\documentclass[11pt]{article}
\usepackage{pdfsync}
\usepackage{hyperref}
\usepackage{array} 
\usepackage{amsfonts}
\usepackage{amsmath}
\usepackage{amssymb}
\usepackage{graphicx}
\usepackage{color}

\newtheorem{thm}{Theorem}[section]
\newtheorem{cor}[thm]{Corollary}
\newtheorem{lem}[thm]{Lemma}
\newtheorem{prop}[thm]{Proposition}

\newtheorem{defn}[thm]{Definition}

\numberwithin{equation}{section}

\newcommand{\RR}{\mathbb{R}}
\newcommand{\NN}{\mathbb{N}}

\def\qed{\,\unskip\kern 6pt \penalty 500d
\raise -2pt\hbox{\vrule \vbox to8pt{\hrule width 6pt
\vfill\hrule}\vrule}\par}

\begin{document}

\title{\textbf{The Fisher-KPP problem \\ [0.7mm]
with doubly nonlinear ``fast'' diffusion.
} \\[3mm]}

\author{{\Large Alessandro Audrito\footnote{Also affiliated with Universit\`a degli Studi di Torino, Italy.} ~and~ Juan Luis V\'azquez}\\ [4pt]
Departamento de Matem\'{a}ticas \\ [4pt] Universidad
Aut\'{o}noma de Madrid
} 
\date{\vspace{-5ex} }

\maketitle

\begin{abstract}
The famous Fisher-KPP reaction diffusion model combines  linear diffusion with the typical Fisher-KPP reaction term, and appears in a number of relevant applications. It is remarkable as a mathematical model since, in the case of linear diffusion, it possesses a family of travelling waves that describe the asymptotic behaviour of a wide class solutions $0\leq u(x,t)\leq 1$  of the problem posed in the real line. The existence of propagation wave with finite speed has been confirmed in the cases of ``slow'' and      ``pseudo-linear'' doubly nonlinear diffusion too, see \cite{AA-JLV:art}. We investigate here the corresponding theory with ``fast'' doubly nonlinear diffusion and we find that general solutions show a non-TW asymptotic behaviour, and exponential propagation in space for large times. Finally, we prove precise bounds for the level sets of general solutions, even when we work in with spacial dimension $N \geq 1$. In particular, we show that location of the level sets is approximately linear for large times, when we take spatial logarithmic scale, finding a strong departure from the linear case, in which appears the famous Bramson logarithmic correction.

\end{abstract}


%
%
%
%
%
%
%
%
%
%
%
\section{Introduction}\label{SECTIONINTRODUCTIONFAST}
In this paper we study the doubly nonlinear (DNL) reaction-diffusion problem posed in the whole Euclidean space
\begin{equation}\label{eq:REACTIONDIFFUSIONEQUATIONPLAPLACIAN}
\begin{cases}
\begin{aligned}
\partial_tu = \Delta_p u^m + f(u) \quad &\text{in } \RR^N\times(0,\infty) \\
u(x,0) = u_0(x) \quad\quad\quad &\text{in } \RR^N.
\end{aligned}
\end{cases}
\end{equation}
We want to describe the asymptotic behaviour of the solution $u = u(x,t)$ for large times and for a specific range of the parameters $m > 0$ and $p > 1$. We recall that the $p$-Laplacian is a nonlinear operator defined for all $1 \leq p < \infty$ by the formula
\[
\Delta_p v := \nabla\cdot(|\nabla v|^{p-2}\nabla v)
\]
and we consider the more general diffusion term
\[
\Delta_p u^m := \Delta_p(u^m) = \nabla\cdot(|\nabla (u^m)|^{p-2}\nabla (u^m)),
\]
called ``doubly nonlinear''operator. Here,  $\nabla$ is the spatial gradient while $\nabla\cdot$ is the spatial divergence. The doubly nonlinear operator (which can be though as the composition of the $m$-th power and the $p$-Laplacian) is much used in the elliptic and parabolic literature (see \cite{C-D-D-S-V:art, DB:book, EstVaz:art, Kal:survey, Lindq:art, V1:book,V2:book} and their references) and allows to recover the Porous Medium operator choosing $p = 2$ or the $p$-Laplacian operator choosing $m = 1$. Of course, choosing $m=1$ and $p = 2$ we obtain the classical Laplacian.

\begin{figure}[h]
  \centering
  \includegraphics[scale=0.45]{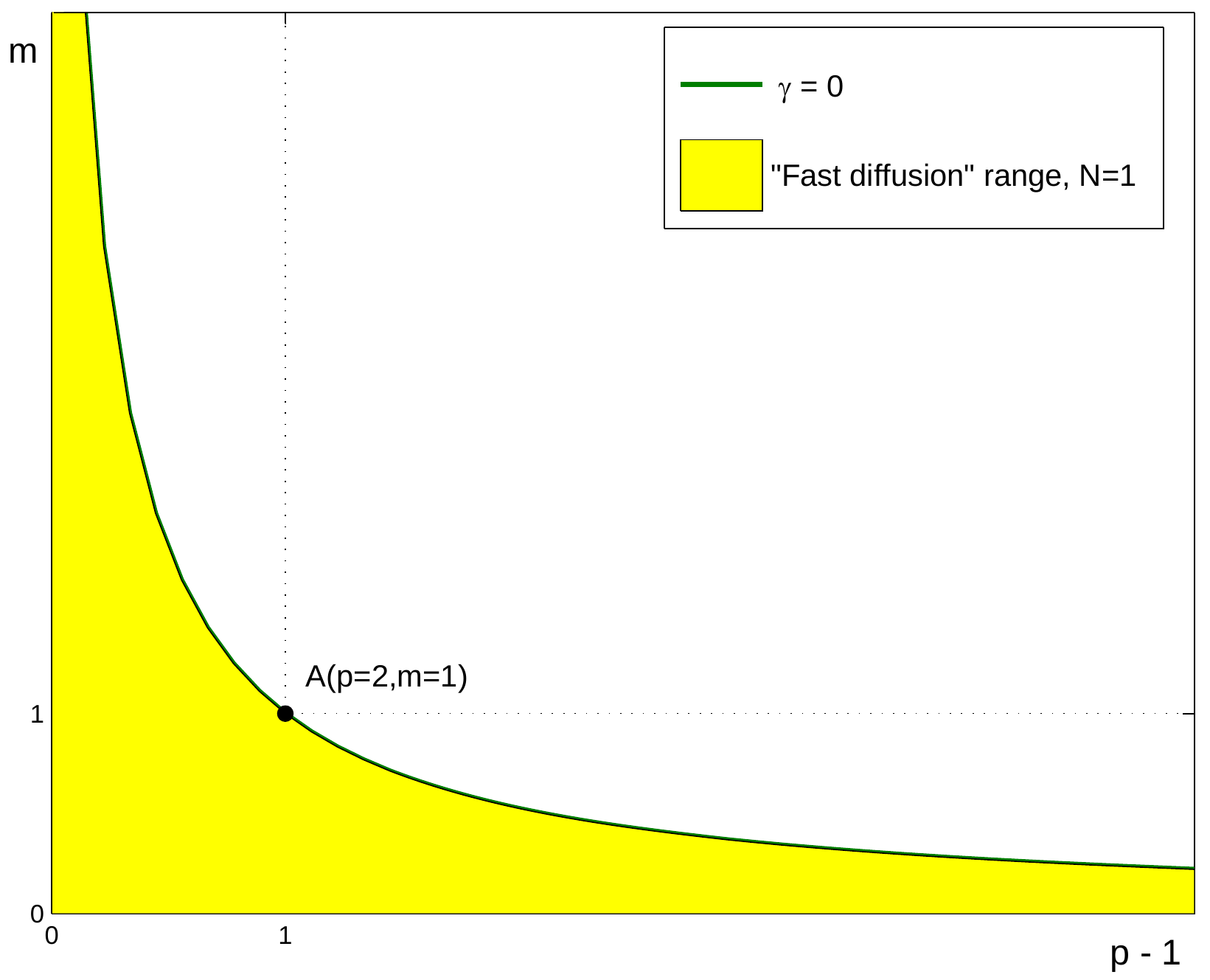} \quad
  \includegraphics[scale=0.45]{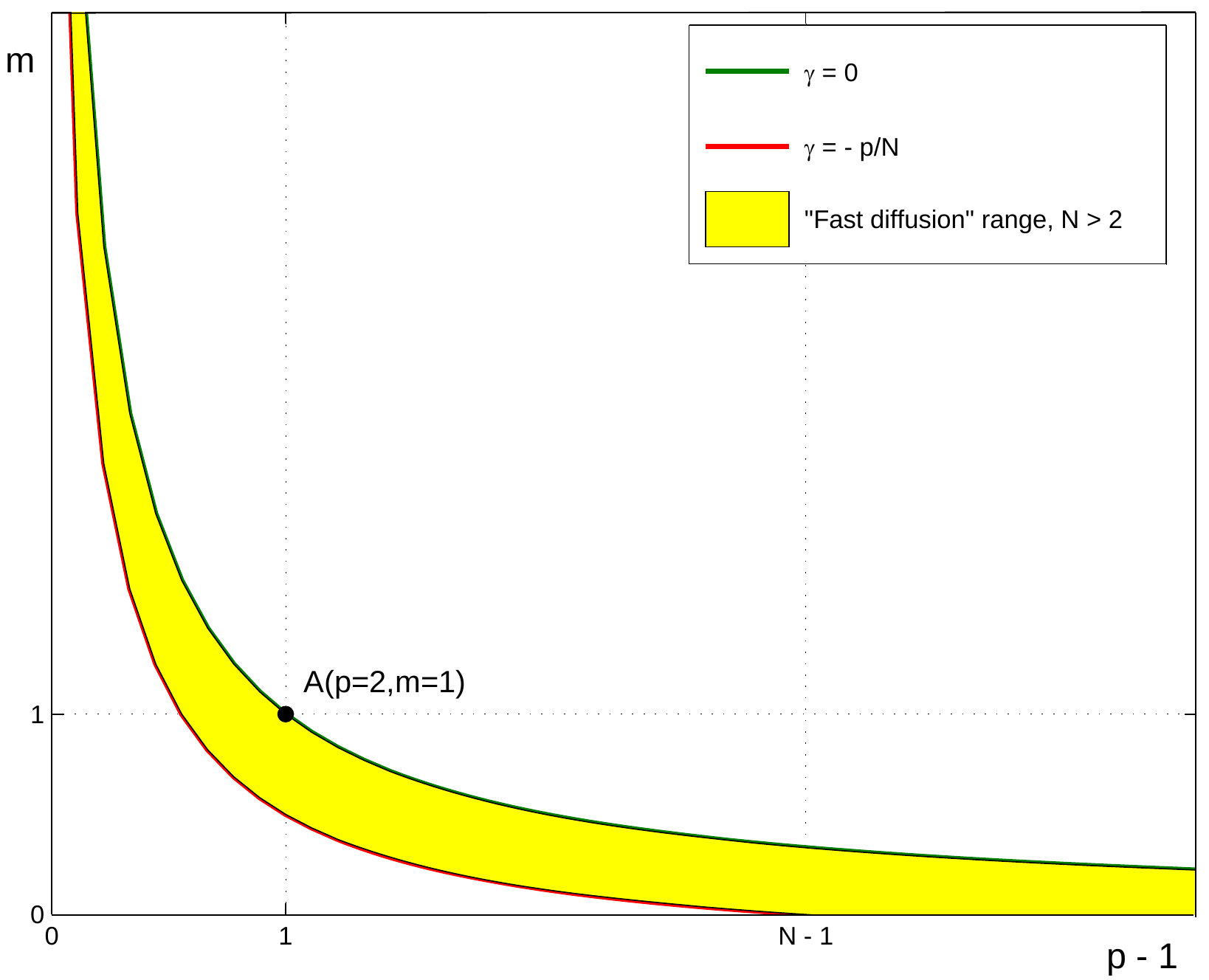}
  \caption{The ``fast diffusion'' range in the $(m,p-1)$-plane.}\label{fig:SIMULFASTCASERANGE}
\end{figure}

\noindent Before proceeding, let us fix some important restrictions and notations. We define the constants
\[
\gamma := m(p-1) - 1 \quad \text{ and } \quad \widehat{\gamma} := -\gamma
\]
and we make the assumption:
\begin{equation}\label{eq:RESTRICTIONONGGFASTDIFF}
-p/N < \gamma < 0 \qquad \text{ i.e. } \qquad 0 < \widehat{\gamma} < p/N
\end{equation}
that we call ``fast diffusion assumption''(cfr. with \cite{V1:book, V2:book}). Note that the shape of the region depends on the dimension $N \geq 1$. Two examples are reported in Figure \ref{fig:SIMULFASTCASERANGE} (note that the region in the case $N=2$ is slightly different respect to the case $N > 2$ and we have not displayed it). We introduce the constant $\widehat{\gamma}$ since its positivity simplifies the reading of the paper and allows us to make the computations simpler to follow.

The case $\gamma \geq 0$, i.e. $\widehat{\gamma} \leq 0$, has been recently studied in \cite{AA-JLV:art}. In this setting, the authors have showed that the equation in \eqref{eq:REACTIONDIFFUSIONEQUATIONPLAPLACIAN} possesses a special class of travelling waves which describe the asymptotic behaviour for large times of more general solutions (see Subsection \ref{SECTIONPREVIOUSRESULTSFAST} for a summary of the results of the case $\gamma \geq 0$). Our main goal is to prove that the case $\gamma < 0$ presents significative departs in the asymptotic properties of the solutions of problem \eqref{eq:REACTIONDIFFUSIONEQUATIONPLAPLACIAN}. In particular, we will see that general solutions do not move with constant speed but with exponential spacial propagation for large times. This fact is the most interesting deviance respect to the classical theory in which TWs play an important role in the study of the asymptotic behaviour.

The function $f(\cdot)$ is a reaction term modeled on the famous references by Fisher \cite{Fisher:art}, and Kolmogorov-Petrovski-Piscounoff \cite{K-P-P:art} in their seminal works on the existence of traveling wave propagation. The classical example is the logistic term $f(u) = u(1-u)$, $0 \leq u \leq 1$. More generally, we will assume that
\begin{equation}\label{eq:ASSUMPTIONSONTHEREACTIONTERM}
\begin{cases}
f : [0,1] \to \RR \text{ and } f \in C^1([0,1]) \\
f(0) = 0 = f(1) \; \text{ and } \; f(u) > 0 \text{ in } (0,1) \\
f \text{ is concave in } [0,1]
\end{cases}
\end{equation}
see  \cite{Aro-Wein1:art, Aro-Wein2:art, Fisher:art, K-P-P:art} for a more complete description of the model. Moreover, we will suppose that the initial datum is a Lebesgue-measurable function and satisfies
\begin{equation}\label{eq:ASSUMPTIONSONTHEINITIALDATUMFAST}
\begin{cases}
u_0(x) \leq C|x|^{-p/\widehat{\gamma}}, \quad \text{for some } C>0 \\
u_0 \not = 0 \; \text{ and } \; 0 \leq u_0 \leq 1.
\end{cases}
\end{equation}
Note that the previous assumption is pretty much general than the more typical continuous with compact support initial data. Moreover, since $0 < \widehat{\gamma} < p/N$, all data satisfying \eqref{eq:ASSUMPTIONSONTHEINITIALDATUMFAST} are automatically integrable, $u_0 \in L^1(\RR^N)$.
\paragraph{Main results and organization of the paper.} The paper is divided in parts as follows:

In Section \ref{SECTIONPREVIOUSRESULTSFAST} we present some known theorems about problem \eqref{eq:REACTIONDIFFUSIONEQUATIONPLAPLACIAN}. Our goal is to give to the reader a quite complete resume on the previous work and related bibliography, to connect it with the new results contained in this paper.

\medskip

In Section \ref{CONVERGENCETOZEROFAST} we begin the study of the asymptotic behaviour of the solutions of problem \eqref{eq:REACTIONDIFFUSIONEQUATIONPLAPLACIAN}-\eqref{eq:ASSUMPTIONSONTHEREACTIONTERM}-
\eqref{eq:ASSUMPTIONSONTHEINITIALDATUMFAST}, with restriction \eqref{eq:RESTRICTIONONGGFASTDIFF}. In particular, we firstly introduce the critical exponent
\[
\sigma_{\ast} := \frac{\widehat{\gamma}}{p}f'(0),
\]
by giving a formal motivation and, later, we prove the following theorem.
\begin{thm}\label{CONVERGENCETOZEROFASTDIFFUSION}
Fix $N \geq 1$. Let $m > 0$ and $p > 1$ such that $0 < \widehat{\gamma} < p/N$. Then for all $\sigma > \sigma_{\ast}$, the solution $u = u(x,t)$ of problem \eqref{eq:REACTIONDIFFUSIONEQUATIONPLAPLACIAN} with initial datum \eqref{eq:ASSUMPTIONSONTHEINITIALDATUMFAST} satisfies
\[
u(x,t) \to 0 \quad \text{ uniformly in } \{|x| \geq e^{\sigma t} \} \; \text{ as } t \to \infty.
\]
\end{thm}
For all $\sigma > \sigma_{\ast}$, we call $\{|x| \geq e^{\sigma t}\}$ ``exponential outer set'' or, simply, ``outer set''. The previous theorem shows that, for large times, the solution $u = u(x,t)$ converges to zero on the ``outer set'' and represents the first step of our asymptotic study.

\medskip

In Section \ref{SECTIONEXPONENTIALEXPANSIONSUPERLEVELSETS} we proceed with the asymptotic analysis, studying the solution of problem \eqref{eq:REACTIONDIFFUSIONEQUATIONPLAPLACIAN} with initial datum
\begin{equation}\label{eq:INITIALDATUMTESTFASTLEVELSETS1}
\widetilde{u}_0(x) :=
\begin{cases}
\begin{aligned}
\widetilde{\varepsilon}  \;\,\qquad\qquad &\text{if }|x| \leq \widetilde{\varrho}_0 \\
a_0|x|^{-p/\widehat{\gamma}} \quad &\text{if } |x| > \widetilde{\varrho}_0,
\end{aligned}
\end{cases}
\end{equation}
where $\widetilde{\varepsilon}$ and $\widetilde{\varrho}_0$ are positive real numbers and $a_0 := \widetilde{\varepsilon}\,\widetilde{\varrho}_0^{\,p/\widehat{\gamma}}$. We show the following crucial proposition.
\begin{prop}\label{EXPANPANSIONOFMINIMALLEVELSETS}
Fix $N \geq 1$. Let $m > 0$ and $p > 1$ such that $0 < \widehat{\gamma} < p/N$ and let $0 < \sigma < \sigma_{\ast}$. Then there exist $t_0 > 0$, $\widetilde{\varepsilon} > 0$ and $\widetilde{\varrho}_0 > 0$ such that the solution $u = u(x,t)$ of problem \eqref{eq:REACTIONDIFFUSIONEQUATIONPLAPLACIAN} with initial datum \eqref{eq:INITIALDATUMTESTFASTLEVELSETS1} satisfies
\[
u(x,t) \geq \widetilde{\varepsilon} \quad \text{in } \{|x| \leq e^{\sigma t}\} \text{ for all } \; t \geq t_0.
\]
\end{prop}
This result asserts that for all initial data \eqref{eq:INITIALDATUMTESTFASTLEVELSETS1} ``small enough'' and for all $\sigma < \sigma_{\ast}$, the solution of problem \eqref{eq:REACTIONDIFFUSIONEQUATIONPLAPLACIAN} is strictly greater than a fixed positive constant on the ``exponential inner sets'' (or ``inner sets'') $\{|x| \leq e^{\sigma t}\}$ for large times. Hence, it proves the non existence of travelling wave solutions (TWs) since ``profiles'' moving with constant speed of propagation cannot describe the asymptotic behaviour of more general solutions (see Section \ref{SECTIONPREVIOUSRESULTSFAST} for the definition of TWs).

\noindent Moreover, this property will be really useful for the construction of sub-solutions of general solutions since, as we will see, it is always possible to place an initial datum with the form \eqref{eq:INITIALDATUMTESTFASTLEVELSETS1} under a general solution of \eqref{eq:REACTIONDIFFUSIONEQUATIONPLAPLACIAN} and applying the Maximum Principle (see Lemma \ref{LEMMAPLACINGBARENBLATTUNDERSOLUTION}).

\medskip

In Section \ref{SECTIONASYMPTOTICBEHAVIOURFAST} we analyze the asymptotic behaviour of the solution of problem \eqref{eq:REACTIONDIFFUSIONEQUATIONPLAPLACIAN}, \eqref{eq:ASSUMPTIONSONTHEINITIALDATUMFAST} in the ``inner sets'' $\{|x| \leq e^{\sigma t}\}$. Along with Theorem \ref{THEOREMBOUNDSFORLEVELSETSFAST} the next theorem is the main result of this paper.
\begin{thm}\label{CONVERGENCETOONEFASTDIFFUSION}
Fix $N \geq 1$. Let $m > 0$ and $p > 1$ such that $0 < \widehat{\gamma} < p/N$. Then for all $\sigma < \sigma_{\ast}$, the solution $u = u(x,t)$ of problem \eqref{eq:REACTIONDIFFUSIONEQUATIONPLAPLACIAN} with initial datum \eqref{eq:ASSUMPTIONSONTHEINITIALDATUMFAST} satisfies
\[
u(x,t) \to 1 \quad \text{ uniformly in } \{|x| \leq e^{\sigma t} \} \; \text{ as } t \to \infty.
\]
\end{thm}
This theorem can be summarized by saying that the function $u(x,t)$ converges to the steady state 1 in the ``inner sets'' for large times. From the point of view of the applications, we can say that the density of population $u = u(x,t)$ invades all the available space propagating exponentially for large times.
\begin{figure}[!ht]
  \centering
  \includegraphics[scale=0.4]{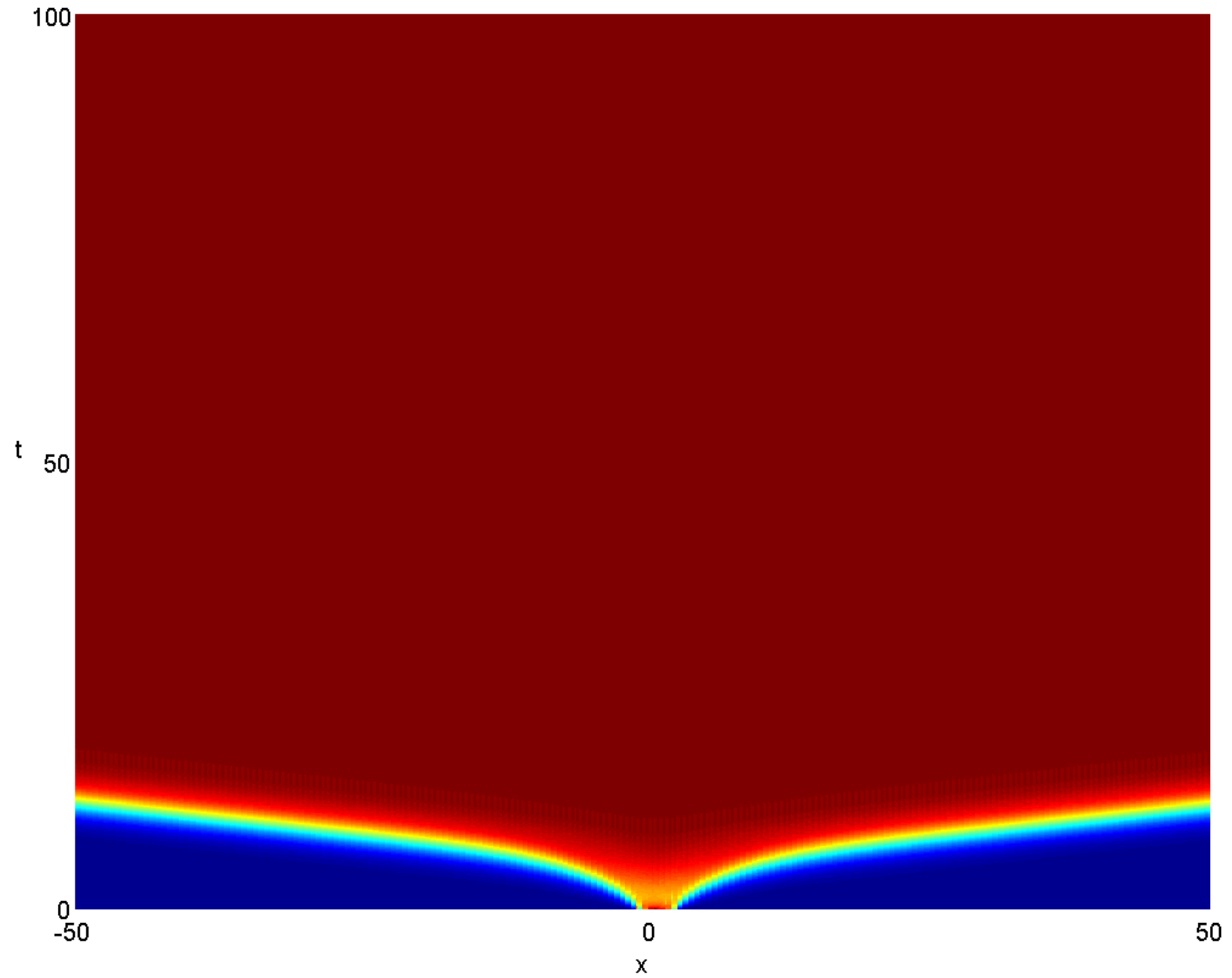} \quad
  \includegraphics[scale=0.45]{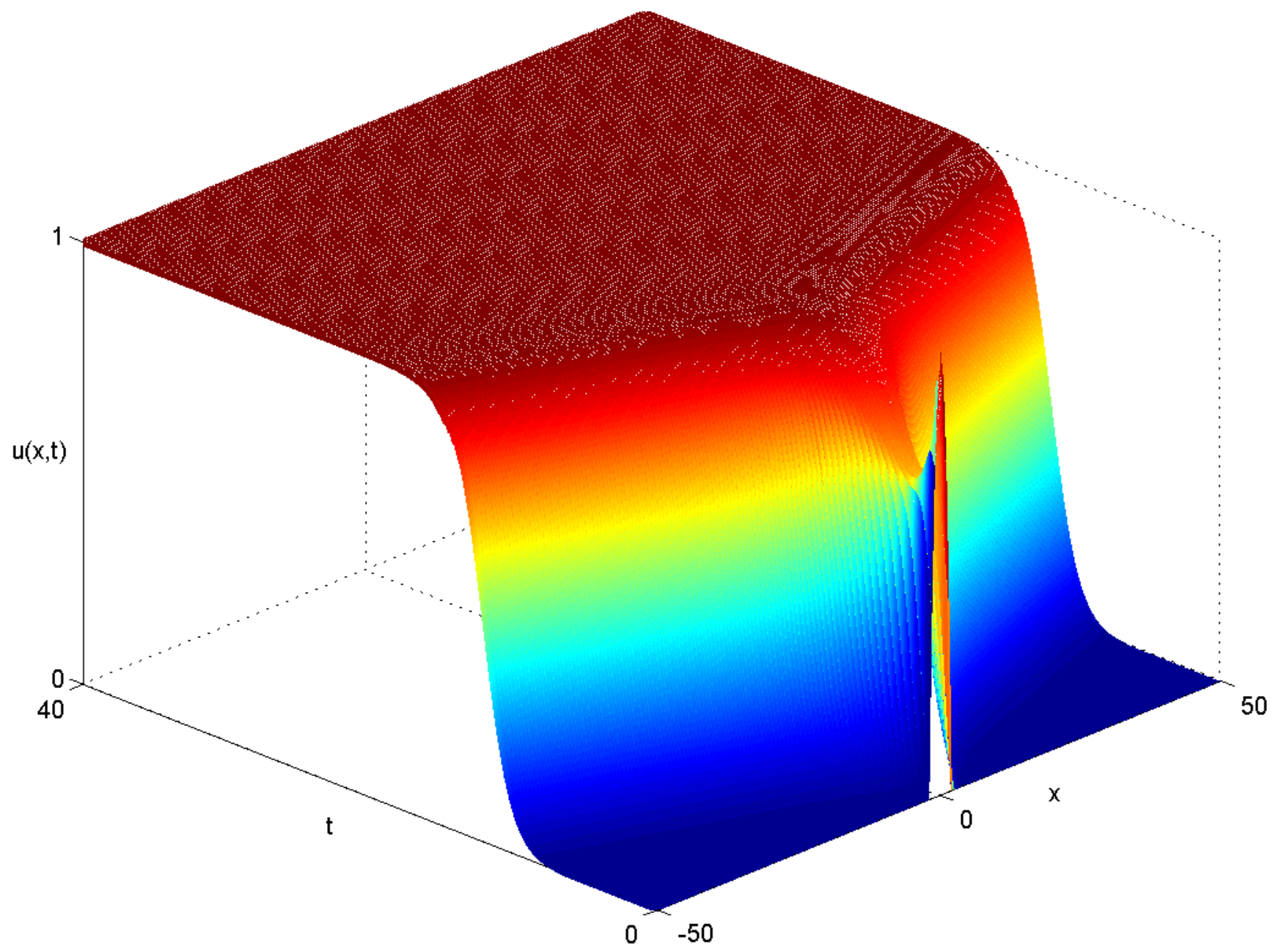}
  \caption{The fast diffusion case: convergence to 1 in the inner sets $\{|x| \leq e^{\sigma t}\}$, for all $\sigma < \sigma_{\ast}$.}\label{fig:SIMULFASTCASEEDPS}
\end{figure}

\medskip

In Section \ref{SECTIONBOUNDSFORLEVELSETSFAST} we consider the classical reaction term $f(u) = u(1-u)$. We find interesting bounds for the level sets of the solution of problem \eqref{eq:REACTIONDIFFUSIONEQUATIONPLAPLACIAN}, \eqref{eq:ASSUMPTIONSONTHEINITIALDATUMFAST}. In particular, we prove that the information on the level sets of the general solutions is contained, up to a multiplicative constant, in the set $|x| = e^{\sigma_{\ast}t}$, for large times.
\begin{thm}\label{THEOREMBOUNDSFORLEVELSETSFAST}
Fix $N \geq 1$. Let $m > 0$ and $p > 1$ such that $0 < \widehat{\gamma} < p/N$, and take $f(u)=u(1-u)$. Then for all $0 < \omega < 1$, there exists a constant $C_{\omega} > 0$ and a time $t_{\omega} > 0$ large enough, such that the solution of problem \eqref{eq:REACTIONDIFFUSIONEQUATIONPLAPLACIAN} with initial datum \eqref{eq:ASSUMPTIONSONTHEINITIALDATUMFAST} and reaction $f(u) = u(1-u)$ satisfies
\begin{equation}\label{eq:STATEMENTLEVELSETSBOUNDSFAST}
\{|x| > C_{\omega}e^{\sigma_{\ast}t}\} \subset \{u(x,t) < \omega \} \quad \text{ and } \quad \{|x| < C_{\omega}^{-1}e^{\sigma_{\ast}t}\} \subset \{u(x,t) > \omega \}
\end{equation}
for all $t \geq t_{\omega}$. In particular, we have
\[
E_{\omega}(t) = \{x \in \RR^N : u(x,t) = \omega \} \subset \{x \in \RR^N : C_{\omega}^{-1}e^{\sigma_{\ast}t} \leq |x| \leq C_{\omega}e^{\sigma_{\ast}t}\} \quad \text{for all } t \geq t_{\omega}.
\]
\end{thm}
An important feature of this result is that for all $0 < \omega < 1$, the set $\{C_{\omega}^{-1}e^{\sigma_{\ast}t} \leq |x| \leq C_{\omega}e^{\sigma_{\ast}t}\}$ does not depend on some $\sigma \not= \sigma_{\ast}$, while in Theorem \ref{CONVERGENCETOZEROFASTDIFFUSION} and Theorem \ref{CONVERGENCETOONEFASTDIFFUSION} the ``outer sets'' and the ``inner sets'' depend on $\sigma > \sigma_{\ast}$ and $\sigma < \sigma_{\ast}$, respectively. Moreover, taking a \emph{spatial logarithmic scale} we can write the estimate
\[
E_{\omega}(t) := \{x \in \RR^N : u(x,t) = \omega \} \subset \{x \in \RR^N : -\ln C_{\omega}\le  \ln|x| - \sigma_{\ast}t \leq   \ln C_{\omega}\},
\]
for $t$ large enough. Actually, this result was not known for ``fast'' nonlinear diffusion neither for the Porous Medium case, nor for the $p$-Laplacian case. However, it was proved by Cabr\'e and Roquejoffre for the fractional Laplacian $(-\Delta)^{1/2}$ in \cite{C-R2:art}, in dimension $N=1$.

In order to fully understand the importance of Theorem \ref{THEOREMBOUNDSFORLEVELSETSFAST}, we need to compare it with the linear case $m=1$ and $p=2$, see formula \eqref{eq:BRAMSONCORRECTION1INTRO}. As we will explain later, in the linear case the location of the level sets is given by a main linear term in $t$ with a logarithmic shift for large times, see \cite{B1:art, B2:art, Hamel-N-R-R:art}. In other words, the propagation of the front is linear ``up to'' a logarithmic correction, for large times. Now, Theorem \ref{THEOREMBOUNDSFORLEVELSETSFAST} asserts that this correction does not occur in the ``fast diffusion'' range. Using the logarithmic scale, we can compare the behaviour of our level sets with the ones of formula \eqref{eq:BRAMSONCORRECTION1INTRO} for linear diffusion, noting that there is no logarithmic deviation, but the location of the level sets is approximately linear for large times (in spatial logarithmic scale, of course), and moreover there is a bounded interval of uncertainty on each level set location.
\normalcolor

\medskip

In Section \ref{SECTIONMAXPRINCCYLDOMAINS} we prove a Maximum Principle for a parabolic equation of $p$-Laplacian type in non-cylindrical domains, see Proposition \ref{MAXPRINNONCYLDOMAINS}. The idea of comparing sub- and super-solutions in non-cylindrical domains comes from \cite{C-R2:art} and it will turn out to be an extremely useful technical tool in the proof of Theorem \ref{CONVERGENCETOONEFASTDIFFUSION}.

\medskip

Section \ref{APPENDIXSELFSIMSOLINCINITDATA} is an appendix. We present some knew results on the existence, uniqueness and regularity for solutions of the ``pure diffusive'' parabolic equation with $p$-Laplacian diffusion and non-integrable initial data. In particular, we focus on radial data $u_0(x) = |x|^{\lambda}$, $\lambda > 0$ and we study some basic properties of the self-similar solutions with datum $u_0 = u_0(x)$. The results of this section are needed for proving Theorem \ref{CONVERGENCETOONEFASTDIFFUSION}.

\medskip

Finally, in Section \ref{SECTIONFINALREMARKSFAST} we conclude the paper with some comments and open problems related to our study. In particular, we focus on the range of parameters
\[
\gamma \leq -p/N  \qquad \text{ i.e. } \qquad \widehat{\gamma} \geq p/N.
\]
The case $\widehat{\gamma} = p/N$ is critical in our study while the range $\widehat{\gamma} > p/N$ is also known in literature as ``very fast'' diffusion range. One of the problems of this range is the lack of basic tools and basic theory (existence, uniqueness, regularity of the solutions and estimates) known for the Porous Medium Equation and for the $p$-Laplacian Equation, but not in the doubly nonlinear setting.
%
%
%
%
%
%
%
\section{Preliminaries and previous results}\label{SECTIONPREVIOUSRESULTSFAST}
In this brief section, for the reader's convenience, we recall some known results about problem \eqref{eq:REACTIONDIFFUSIONEQUATIONPLAPLACIAN} with related bibliography, and we introduce some extremely useful tools such as ``Barenblatt solutions'', we will need through the paper.
%
%
%
%
%
%
\subsection{Previous results}
As we have explained before, we present here the literature and past works linked to our paper, in order to motivate our study. Basically, the goal is to give to the reader a suitable background on the Fisher-KPP theory, so that our new results can be compared and fully understood.
\paragraph{Finite propagation: the doubly nonlinear case.} In \cite{AA-JLV:art}, we studied problem \eqref{eq:REACTIONDIFFUSIONEQUATIONPLAPLACIAN} assuming $\gamma = m(p-1) - 1 \geq 0$ and initial datum satisfying
\begin{equation}\label{eq:ASSUMPTIONSONTHEINITIALDATUM}
\begin{cases}
u_0 : \RR^N \to \RR \text{ is continuous with compact support: } u_0 \in \mathcal{C}_c(\RR^N) \\
u_0 \not \equiv 0 \; \text{ and } \; 0 \leq u_0 \leq 1.
\end{cases}
\end{equation}
Before enunciating the main results we need to introduce the notion of Travelling Waves (TWs). They are special solutions with remarkable applications, and there is a huge mathematical literature devoted to them. Let us review the main concepts and definitions.

\begin{figure}[!ht]
  \centering
  \includegraphics[scale=0.4]{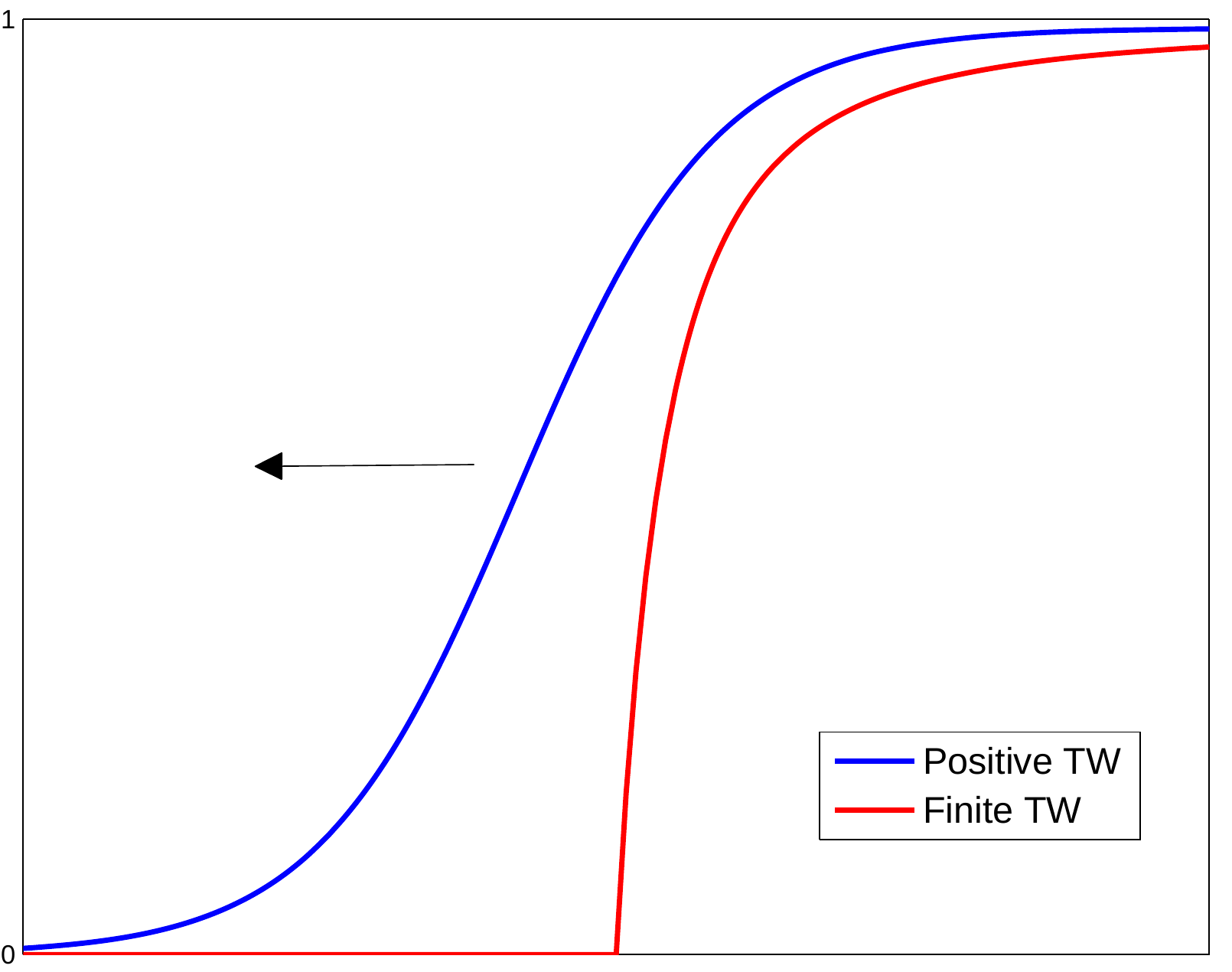}
  \includegraphics[scale=0.4]{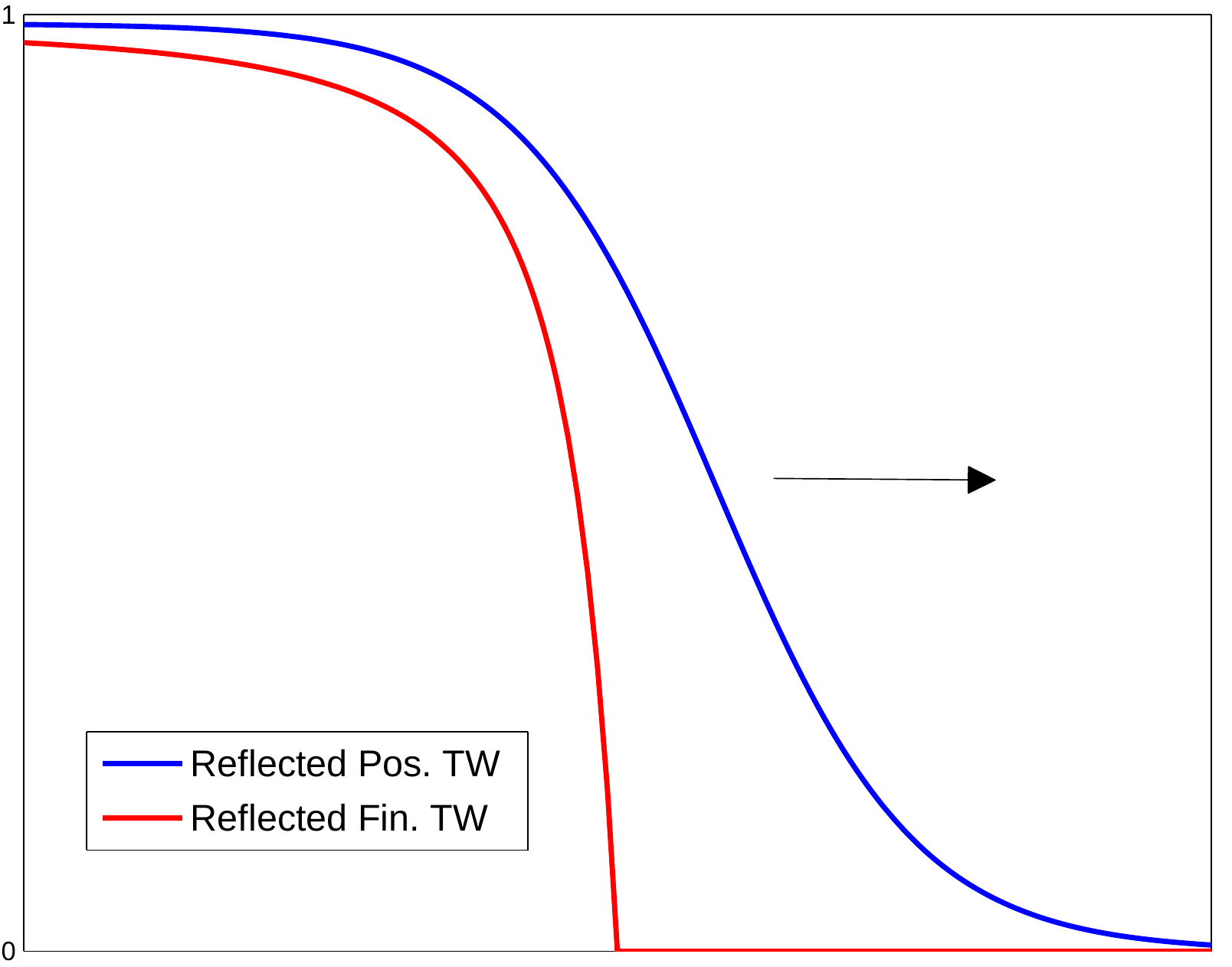}
  \caption{Examples of admissible TWs and their ``reflections'': Finite and Positive types}\label{fig:ADMISSIBLEANDREFLECTEDFINITETW}
\end{figure}

Fix $m > 0$ and $p \geq 1$ and assume that we are in space dimension 1 (note that when $N = 1$, the doubly nonlinear operator has the simpler expression $\Delta_pu^m = \partial_x(|\partial_xu^m|^{p-2}\partial_xu^m)$). A TW solution of the equation
\begin{equation}\label{eq:PLAPLACIANREACTIONDIFFUSIONEQUATIONINTRO}
\partial_tu = \Delta_pu^m + f(u) \quad \text{in }\RR\times(0,\infty)
\end{equation}
is a solution of the form $u(x,t) = \varphi(\xi)$, where $\xi = x + ct$, $c > 0$ and the \emph{profile} $\varphi(\cdot)$ is a real function. In the application to the Fisher-KPP problem, the profile is assumed to satisfy
\begin{equation}\label{eq:CONDITIONONPHIADMISSIBLETWINTRO}
0 \leq \varphi \leq 1, \quad \varphi(-\infty) = 0, \; \varphi(\infty) = 1 \quad \text{and} \quad \varphi' \geq 0.
\end{equation}
In that case we say that $u(x,t) = \varphi(\xi)$ is an \emph{admissible} TW solution. Similarly, one can consider admissible TWs of the form $u(x,t) = \varphi(\xi)$ with $\xi = x - ct$, $\varphi$ decreasing and such that $\varphi(-\infty) = 1$ and $\varphi(\infty) = 0$.  But it is easy to see that  these two options are equivalent, since the the shape of the profile of the second one can be obtained by reflection of the first one, $\varphi_{-c}(\xi)=\varphi_c(-\xi)$, and it moves in the opposite direction of propagation.
\\
Finally, an admissible TW is said \emph{finite} if $\varphi(\xi) = 0$ for $\xi \leq \xi_0$ and/or $\varphi(\xi) = 1$ for $\xi \geq \xi_1$, or \emph{positive} if $\varphi(\xi) > 0$, for all $\xi \in \RR$. The line $x = \xi_0 - ct$ that separates the regions of positivity and vanishing of $u(x,t)$ is then called the \emph{free boundary}. Now, we can proceed.

We proved that the existence of admissible TW solutions depends on the wave's speed of propagation. In particular, we showed the following theorem, cfr. with Theorem 2.1 and Theorem 2.2 of \cite{AA-JLV:art}.
\begin{thm}\label{THEOREMEXISTENCEOFTWS}
Let $m > 0$ and $p > 1$ such that $\gamma \geq 0$. Then there exists a unique $c_{\ast} = c_{\ast}(m,p) > 0$ such that equation \eqref{eq:PLAPLACIANREACTIONDIFFUSIONEQUATIONINTRO} possesses a unique admissible TW for all $c \geq c_{\ast}(m,p)$ and does not have admissible TWs for $0 < c < c_{\ast}(m,p)$. Uniqueness is intended up to reflection or horizontal displacement.

\noindent Moreover, if $\gamma > 0$, the TW corresponding to the value $c = c_{\ast}(m,p)$ is finite (i.e., it vanishes in an infinite half-line), while the TWs corresponding to the values $c > c_{\ast}(m,p)$ are positive everywhere.

\noindent Finally, when $\gamma = 0$, any admissible TW is positive everywhere.
\end{thm}
The concept of admissible TWs and the problem of their existence was firstly introduced in \cite{Fisher:art} and \cite{K-P-P:art}. Then Aronson and Weinberger, see \cite{Aro-Wein1:art, Aro-Wein2:art}, proved Theorem \ref{THEOREMEXISTENCEOFTWS} in the case of the linear diffusion, i.e. $m=1$ and $p=2$ (note that the choice $m=1$ and $p=2$ is a subcase of $\gamma = 0$). Later, the problem of the existence of critical speeds and admissible TWs for the Fisher-KPP equation has been studied for the Porous Medium diffusion ($m > 1$ and $p=2$), see \cite{Aronson1:art,Aronson2:art} and \cite{DP-S:art,DePablo-Vazquez1:art,DePablo-Vazquez2:art}. Recently, see \cite{Eng-Gav-San:art, Gav-San:art}, it has been proved the existence of admissible TWs and admissible speeds of propagation when $m=1$ and $p>2$, i.e. $p$-Laplacian diffusion. In Theorem \ref{THEOREMEXISTENCEOFTWS} we generalized these results when doubly nonlinear diffusion is considered and $\gamma \geq 0$.

\noindent Then we focused on the PDE part in which we studied the asymptotic behaviour of more general solutions, proving the following theorem (Theorem 2.6 of \cite{AA-JLV:art}).
\begin{thm}\label{NTHEOREMCONVERGENCEINNEROUTERSETS}
Fix $N \geq 1$. Let $m > 0$ and $p > 1$ such that $\gamma \geq 0$.

\noindent (i) For all $0 < c < c_{\ast}(m,p)$, the solution $u(x,t)$ of the initial-value problem \eqref{eq:REACTIONDIFFUSIONEQUATIONPLAPLACIAN} with initial datum \eqref{eq:ASSUMPTIONSONTHEINITIALDATUM} satisfies
\[
u(x,t) \to 1 \text{ uniformly in } \{|x| \leq ct\} \;\text{ as } t \to \infty.
\]
\noindent (ii) Moreover, for all $c > c_{\ast}(m,p)$ it satisfies,
\[
u(x,t) \to 0 \text{ uniformly in } \{|x| \geq ct\} \;\text{ as } t \to \infty.
\]
\end{thm}
In the case in which the function $u = u(x,t)$ stands for a density of population, the statement of the previous theorem means that the individuals tend to occupy all the available space and, for large times, they spread with constant speed, see \cite{Aro-Wein1:art, Aro-Wein2:art, Fisher:art, K-P-P:art}. From the mathematical point of view, we can state that the steady state $u = 1$ is asymptotically stable while the null solution $u = 0$ is unstable and, furthermore, the asymptotic stability/instability can be measured in terms of speed of convergence of the solution which, in this case, is asymptotically linear in distance of the front location as function of time.

\noindent Again we recall that for $m=1$ and $p=2$, the previous theorem was showed in \cite{Aro-Wein1:art, Aro-Wein2:art}, while for $m>1$ and $p=2$ in \cite{DP-S:art}. We point out that in this last paper, the authors worked with a slightly different reaction term, they called ``strong reaction'', see also \cite{DePablo-Vazquez1:art,DePablo-Vazquez2:art}.

\noindent In the linear case $m=1$ and $p=2$, the statements of Theorem \ref{NTHEOREMCONVERGENCEINNEROUTERSETS} were improved. Indeed, when $N=1$, Bramson showed an interesting property of the level sets $E_{\omega}(t) = \{x > 0: u(x,t) = \omega \}$, $\omega \in (0,1)$, of the solution $u = u(x,t)$ of equation \eqref{eq:PLAPLACIANREACTIONDIFFUSIONEQUATIONINTRO} (with $m=1$ and $p=2$) with reaction term satisfying \eqref{eq:ASSUMPTIONSONTHEREACTIONTERM}. In particular, in \cite{B1:art} and \cite{B2:art}, it was proved that for all $\omega \in (0,1)$ there exist constants $x_{\omega}$, $a > 0$ and $C_{\omega} > 0$ such that
\begin{equation}\label{eq:BRAMSONCORRECTION1INTRO}
E_{\omega}(t) \subset \Big[c_{\ast}t - \frac{3}{2\omega_{\ast}}\ln t - x_{\omega} - \frac{a}{\sqrt{t}}- \frac{C_{\omega}}{t}, c_{\ast}t - \frac{3}{2\omega_{\ast}}\ln t - x_{\omega} - \frac{a}{\sqrt{t}} + \frac{C_{\omega}}{t}\Big]
\end{equation}
for $t$ large enough, where $\omega_{\ast} = c_{\ast}/2$. The previous formula is interesting since it allows to estimate the ``delay'' of the solution $u = u(x,t)$ from the positive TW with critical speed $c = c_{\ast}$ which, according to \eqref{eq:BRAMSONCORRECTION1INTRO}, grows in time and consists in a logarithmic deviance. Furthermore, he showed that general solutions converge uniformly to the TW with critical speed of propagation (once it is ``shifted'' of a logarithmic factor $3/(2\omega_{\ast})\ln t$), for large times. More recently, similar results have been proved in \cite{Hamel-N-R-R:art, NolRoqRyz2015:art, NolRoqRyz2016:art} with PDEs techniques.
\paragraph{Exponential propagation for Porous Medium fast diffusion.} Let's now resume briefly the results of King and McCabe (\cite{King-M:art}), which have inspired this paper. They considered the Porous Medium case, which is obtained by taking $p=2$ in the equation in \eqref{eq:REACTIONDIFFUSIONEQUATIONPLAPLACIAN}:
\begin{equation}\label{eq:REACTIONDIFFUSIONEQUATIONPORMEDKINGMCCABE}
\begin{cases}
\begin{aligned}
\partial_tu = \nabla\cdot(u^{-(1-m)}\nabla u) + u(1-u) \quad &\text{in } \RR^N\times(0,\infty) \\
u(x,0) = u_0(x) \qquad\qquad\qquad\qquad\quad\; &\text{in } \RR^N,
\end{aligned}
\end{cases}
\end{equation}
in the fast diffusion range, $0 < m < 1$ (note that we absorbed a factor $m$ by using a simple change of variables). They considered non-increasing radial initial data $u_0 \in L^1(\RR^N)$ decaying faster than $r^{-2/(1-m)}$ as $r = |x| \sim \infty$ and studied radial solutions of problem \eqref{eq:REACTIONDIFFUSIONEQUATIONPORMEDKINGMCCABE}.

\noindent They showed that when $(N-2)_+/N := m_c < m < 1$, the radial solutions $u = u(r,t)$ of the previous equation converge pointwise to 1 for large times with exponential rate $r(t) \sim e^{\sigma t}$, for $\sigma < (1-m)/2$, and that the ``main variation in concentration occurs on the scale $O(e^{nt/2})$'', (see \cite{King-M:art} pag. 2544) where $n = 1-m$ in their notation. We will present an adaptation of their methods to our case in Section \ref{CONVERGENCETOZEROFAST}.

\noindent Note that ``our'' critical exponent $\sigma_{\ast}$ generalizes the value $(1-m)/2$ to the case $p > 1$ and to more general reaction terms than $f(u) = u(1-u)$. This is a severe departure from the TW behavior of the standard Fisher-KPP model since there are no TW solutions. Instead, they found that radial solutions of \eqref{eq:REACTIONDIFFUSIONEQUATIONPORMEDKINGMCCABE} have non-TWs form for large times:
\[
u(r,t) \sim \frac{(\kappa r^{-2}e^{n t})^{1/n}}{1 + (\kappa r^{-2}e^{n t})^{1/n}}, \quad t \sim \infty \quad \text{ and } \quad r = O(e^{nt/2}),
\]
where $n = 1-m$ and $\kappa = 2(2-nN)/n^2$. The case $0 < m < m_c$ is studied too and, as we have anticipated, we will discuss this range in the final section with some comments.

\noindent Exponential propagation happens also with fractional diffusion, both linear and nonlinear, see for instance \cite{C-R1:art,C-R2:art,S-V:art} and the references therein.  We will not enter here into the study of the relations of our paper with nonlinear fractional diffusion, though it is an interesting topic.

\noindent Finally, we recall that infinite speed of propagation depends not only on the diffusion operator but also on the initial datum. In particular, in \cite{Hamel-R:art}, Hamel and Roques found that the solutions of the Fisher-KPP problem with linear diffusion i.e., ($m=1$ and $p=2$) propagate exponentially fast for large times if the initial datum has a power-like spatial decay at infinity.

\noindent The scene is set for us to investigate what happens in the presence of a fast doubly nonlinear  diffusion.
%
%
%
%
%
%
\subsection{Preliminaries on doubly nonlinear diffusion.}\label{PRELIMINARIESINTRO}
Now we present some basic results concerning the Barenblatt solutions of the ``pure diffusive'' doubly nonlinear parabolic equation which are essential to develop our study in the next sections (the reference for this issue is \cite{V1:book}). Moreover, we recall some basic facts on existence, uniqueness, regularity and Maximum Principles for the solutions of problem \eqref{eq:REACTIONDIFFUSIONEQUATIONPLAPLACIAN}.
\paragraph{Barenblatt solutions.} Fix $m > 0$ and $p > 1$ such that $0 < \widehat{\gamma} < p/N$ and consider the ``pure diffusive'' doubly nonlinear problem:
\begin{equation}\label{eq:PARABOLICPLAPLACIANEQUATIONINTRO}
\begin{cases}
\begin{aligned}
\partial_tu = \Delta_p u^m \;\quad &\text{in } \RR^N\times(0,\infty) \\
u(t) \to M\delta_0 \quad\;\, &\text{in } \RR^N \text{ as } t \to 0,
\end{aligned}
\end{cases}
\end{equation}
where $M\delta_0(\cdot)$ is the Dirac's function with mass $M>0$ in the origin of $\RR^N$ and the convergence has to be intended in the sense of measures.

\noindent It has been proved (see \cite{V1:book}) that problem \eqref{eq:PARABOLICPLAPLACIANEQUATIONINTRO} admits continuous weak solutions in self-similar form $B_M(x,t) = t^{-\alpha}F_M(xt^{-\alpha/N})$, called Barenblatt solutions, where the \emph{profile} $F_M(\cdot)$ is defined by the formula:
\[
F_M(\xi) = \Big[C_M + k|\xi|^{\frac{p}{p-1}} \Big]^{-\frac{p-1}{\widehat{\gamma}}},
\]
where
\[
\alpha = \frac{1}{p/N - \widehat{\gamma}}, \quad k = \frac{\widehat{\gamma}}{p}\Big(\frac{\alpha}{N}\Big)^{\frac{1}{p-1}},
\]
$C_M>0$ is determined in terms of the mass choosing $M = \int_{\RR^N}B_M(x,t)dx$ (see \cite{V1:book} for a complete treatise). We point out that there is an equivalent formulation (see \cite{V6:art} for the case $p=2$) in which the Barenblatt solutions are written in the form
\begin{equation}\label{eq:SECONDFORMULATIONBARSOLFAST}
B_M(x,t) = R(t)^{-N} \Big[ D + (\widehat{\gamma}/p)\big|x R(t)^{-1} \big|^{\frac{p}{p-1}} \Big]^{-\frac{p-1}{\widehat{\gamma}}}, \qquad R(t) = \big[(N/\alpha)t \big]^{\frac{\alpha}{N}},
\end{equation}
where $D > 0$ is a new constant. It will be useful to keep in mind that we have the formula
\begin{equation}\label{eq:RELATIONMASSESBARENBLATTSOLINTRO}
B_M(x,t) = MB_1(x,M^{-\widehat{\gamma}}t)
\end{equation}
which describes the relationship between the Barenblatt solution of mass $M > 0$ and mass $M = 1$ and the estimates on the profile corresponding to the Barenblatt solution of mass $M>0$:
\begin{equation}\label{eq:ESTIMATEONPROFILEOFMASS1}
K_2 (1 + |\xi|^{p/\widehat{\gamma}})^{-1} \leq F_M(\xi) \leq K_1|\xi|^{-p/\widehat{\gamma}} \quad \text{for all } \xi \in \RR^N
\end{equation}
for suitable positive constants $K_1$ and $K_2$ depending on $M>0$.
\paragraph{Existence, Uniqueness, Regularity and Maximum Principles.} Before presenting the main results of this paper, we briefly discuss the basic properties of the solutions of problem \eqref{eq:REACTIONDIFFUSIONEQUATIONPLAPLACIAN}. Results about existence of weak solutions of the pure diffusive problem and its generalizations, can be found in the survey \cite{Kal:survey} and the large number of  references therein. The problem of uniqueness was studied later (see for instance \cite{DBen-Her1:art, DBen-Her2:art, Tsu:art, V2:book, Wu-Yin-Li:art}). The classical reference for the regularity of nonlinear parabolic equations is \cite{L-S-U:book}, followed by a wide literature. For the Porous Medium case ($p=2$) we refer to \cite{V1:book, V2:book}, while for the $p$-Laplacian case we suggest \cite{DB:book, Lindq:art} and the references therein. Finally, in the doubly nonlinear setting, we refer to \cite{Ag-Blan-Car:art, For-Sos-Ves:art, For-Sos-Ves1:art, Ivan:art, Ves:art}. The results obtained show the H\"{o}lder continuity of the solution of problem \eqref{eq:REACTIONDIFFUSIONEQUATIONPLAPLACIAN}. Finally, we mention \cite{DB:book, V2:book, Wu-Yin-Li:art} for a proof of the Maximum Principle.
%
%
%
%
%
%
%
%
%
%
%
%
%
\section{Convergence to 0 in the ``outer sets''}\label{CONVERGENCETOZEROFAST}
In this section we study the asymptotic behaviour of the solution of the Cauchy problem \eqref{eq:REACTIONDIFFUSIONEQUATIONPLAPLACIAN} with non-trivial initial datum $u_0(\cdot)$ satisfying \eqref{eq:ASSUMPTIONSONTHEINITIALDATUMFAST}:
\[
u_0(x) \leq C|x|^{-\frac{p}{\widehat{\gamma}}} \quad \text{ and } \quad 0 \leq u_0 \leq 1,
\]
for some constant $C > 0$. We recall here the definition of the critical exponent
\[
\sigma_{\ast} = \frac{\widehat{\gamma}}{p}f'(0), \qquad 0 < \widehat{\gamma} < p/N.
\]
Before proceeding, let us see how to formally derive the value of the critical exponent $\sigma_{\ast}$ in the case $0 < \widehat{\gamma} < p/N$ and $f(u) = u(1-u)$ (note that $f'(0) = 1$). We follow the methods used in \cite{King-M:art}.

\noindent First of all, we fix $0 < \widehat{\gamma} < p/N$ and we consider radial solutions of the equation in \eqref{eq:REACTIONDIFFUSIONEQUATIONPLAPLACIAN}, which means
\[
\partial_t u = r^{1-N}\partial_r\big(r^{N-1}|\partial_ru^m|^{p-2}\partial_ru^m\big) + u(1-u), \quad r > 0, \; t > 0.
\]
Note that the authors of \cite{King-M:art} worked with a slightly different equation (they absorbed the multiplicative factor $m^{p-1}$ with a simple change of variables). We linearize the reaction term and we assume that $u = u(r,t)$ satisfies
\begin{equation}\label{eq:LINEARIZEDEQUATIONKINGMCCABE}
\partial_t u = r^{1-N}\partial_r\big(r^{N-1}|\partial_ru^m|^{p-2}\partial_ru^m\big) + u, \quad \text{for } r \sim \infty.
\end{equation}
Now, we look for a solution of \eqref{eq:LINEARIZEDEQUATIONKINGMCCABE} of the form $u(r,t) \sim r^{-p/\widehat{\gamma}}G(t)$ for $r \sim \infty$ which agrees with the assumption \eqref{eq:ASSUMPTIONSONTHEINITIALDATUMFAST} on the initial datum and with the linearization \eqref{eq:LINEARIZEDEQUATIONKINGMCCABE}. It is straightforward to see that for such solution, the function $G = G(t)$ has to solve the equation
\begin{equation}\label{eq:EQUATIONOFGKINGMCCABE}
\frac{dG}{dt} = G + \kappa G^{1-\widehat{\gamma}}, \quad t \geq 0 \quad\quad \kappa:= \frac{(p-\widehat{\gamma} N)(mp)^{p-1}}{\widehat{\gamma}^p}.
\end{equation}
Note that since $0 <\widehat{\gamma} < p/N$, we have that $\kappa$ is well defined and positive, while $1-\widehat{\gamma} = m(p-1) > 0$. Equation \eqref{eq:EQUATIONOFGKINGMCCABE} belongs to the famous Bernoulli class and can be explicitly integrated:
\[
G(t) = \big(ae^{\widehat{\gamma} t} - \kappa\big)^{\frac{1}{\widehat{\gamma}}}, \quad a \geq \kappa, \quad t \geq 0.
\]
Hence, for all fixed $t \geq 0$, we obtain the asymptotic expansion for our solution
\begin{equation}\label{eq:1EXPANSIONSOLUTIONINGMCCABE}
u(r,t) \sim r^{-\frac{p}{\widehat{\gamma}}}\big(ae^{\widehat{\gamma} t} - \kappa\big)^{\frac{1}{\widehat{\gamma}}}, \quad r \sim \infty.
\end{equation}
Now, for all fixed $r > 0$, we consider a solution $\zeta_0 = \zeta_0(r,t)$ of the \emph{logistic} equation
\[
\partial_t\zeta_0 = \zeta_0(1-\zeta_0), \quad t \geq 0,
\]
which describes the state in which there is not diffusion and the dynamics is governed by the reaction term. We assume to have
\[
u(r,t) \sim \zeta_0(r,t) \quad \text{for } t \sim \infty, 
\]
where the leading-order term $\zeta_0 = \zeta_0(r,t)$ satisfies
\begin{equation}\label{eq:SOLUTIONOFLOGISTICKINGMCCABE}
\zeta_0 (r,t) \sim \frac{\phi(r)e^t}{1 + \phi(r)e^t}, \quad \text{for } t \sim \infty, \;\; r \sim \infty,
\end{equation}
for some unknown function $\phi = \phi(r)$, with $\phi(r) \to 0$, as $r \to \infty$. Now, matching \eqref{eq:1EXPANSIONSOLUTIONINGMCCABE} with \eqref{eq:SOLUTIONOFLOGISTICKINGMCCABE} for $t$ large and $r \sim \infty$, we easily deduce
\[
\phi(r) \sim (ar^{-p})^{1/\widehat{\gamma}}, \quad \text{for } r \sim \infty.
\]
Thus, substituting $\phi(r) \sim (a r^{-p})^{1/\widehat{\gamma}}$ in \eqref{eq:SOLUTIONOFLOGISTICKINGMCCABE} and taking $r \sim e^{\widehat{\gamma}/p t}$ for $t \sim \infty$, we have
\[
u(r,t) \sim \frac{(a r^{-p}e^{\widehat{\gamma} t})^{1/\widehat{\gamma}}}{1 + (a r^{-p}e^{\widehat{\gamma} t})^{1/\widehat{\gamma}}} = \frac{\widehat{a} e^t}{r^{p/\widehat{\gamma}} + \widehat{a} e^t} \quad \text{for } t \sim \infty, \;\; r \sim e^{\widehat{\gamma}/p t},
\]
where $\widehat{a} = a^{1/\widehat{\gamma}} \geq \kappa^{1/\widehat{\gamma}}$. The previous formula corresponds to a ``similarity reduction'' (see \cite{King-M:art}, pag. 2533) of the logistic equation with $\zeta_0 = \zeta_0(r/e^{\widehat{\gamma}/p t})$.

\noindent Note that taking $r \geq e^{\sigma t}$ and $\sigma > \widehat{\gamma}/p$, we have $u(r,t) \sim 0$ for $t \sim \infty$ while if $r \leq e^{\sigma t}$ and $\sigma < \widehat{\gamma}/p$ we have $u(r,t) \sim 1$ for $t \sim \infty$. This means that setting $\sigma_{\ast} = \widehat{\gamma}/p$, $r(t) \sim e^{\sigma_{\ast} t}$ is a ``critical'' curve, in the sense that it separates the region in which the solution $u = u(r,t)$ converges to $u=0$ to the one which converges to $u=1$. We will show this property in Theorem \ref{CONVERGENCETOZEROFASTDIFFUSION} and Theorem \ref{CONVERGENCETOONEFASTDIFFUSION}.

In what follows, we prove that the solution of problem \eqref{eq:REACTIONDIFFUSIONEQUATIONPLAPLACIAN} with initial datum \eqref{eq:ASSUMPTIONSONTHEINITIALDATUMFAST} converges uniformly to the trivial solution $u = 0$ in the outer set $\{|x| \geq e^{\sigma t} \}$ as $t \to \infty$ if $\sigma > \sigma_{\ast}$. In the nex sections we will prove that this solution converges uniformly to the equilibrium point $u = 1$ in the inner set $\{|x| \leq e^{\sigma t} \}$ as $t \to \infty$ if $\sigma < \sigma_{\ast}$.
\paragraph{Proof of Theorem \ref{CONVERGENCETOZEROFASTDIFFUSION}.} Fix $N \geq 1$, $0 < \widehat{\gamma} < p/N$, and $\sigma > \sigma_{\ast}$. First of all, we construct a super-solution for problem \eqref{eq:REACTIONDIFFUSIONEQUATIONPLAPLACIAN}, \eqref{eq:ASSUMPTIONSONTHEINITIALDATUMFAST} using the hypothesis on the function $f(\cdot)$. Indeed, since $f(u) \leq f'(0)u$ for all $0 \leq u \leq 1$, the solution of the linearized problem
\[
\begin{cases}
\begin{aligned}
\partial_t\overline{u} = \Delta_p\overline{u}^m + f'(0)\overline{u} \quad &\text{in } \RR^N\times(0,\infty) \\
\overline{u}(x,0) = u_0(x) \;\;\quad\qquad &\text{in } \RR^N,
\end{aligned}
\end{cases}
\]
gives the super-solution we are interested in and, by the Maximum Principle, we deduce $u(x,t) \leq \overline{u}(x,t)$ in $\RR^N\times(0,\infty)$. Now, consider the change of the time variable
\[
\tau(t) = \frac{1}{f'(0)\widehat{\gamma}}\Big[1 - e^{-f'(0)\widehat{\gamma} t} \Big], \quad \text{for } t \geq 0,
\]
with $0 \leq \tau(t) \leq \tau_{\infty} := \frac{1}{f'(0)\widehat{\gamma}}$. Then the function $\overline{v}(x,\tau) = e^{-f'(0)t}\overline{u}(x,t)$ solves the problem
\[
\begin{cases}
\begin{aligned}
\partial_{\tau}\overline{v} = \Delta_p\overline{v}^m \,\qquad &\text{in } \RR^N\times(0,\tau_{\infty}) \\
\overline{v}(x,0) = u_0(x) \quad &\text{in } \RR^N.
\end{aligned}
\end{cases}
\]
From the properties of the profile of the Barenblatt solutions and the hypothesis on the initial datum \eqref{eq:ASSUMPTIONSONTHEINITIALDATUMFAST}, it is evident that there exist positive numbers $M$ and $\theta$ such that $u_0(x) \leq B_M(x,\theta)$ in $\RR^N$ and so, by comparison, we obtain
\[
\overline{v}(x,\tau) \leq B_M(x,\theta + \tau) \quad \text{in } \RR^N\times(0,\tau_{\infty}).
\]
Now, since the profile of Barenblatt solutions satisfies $F_1(\xi) \leq K_1 |\xi|^{-\frac{p}{\widehat{\gamma}}}$ for some constant $K_1 > 0$ and for all $\xi \in \RR^N$ (see \eqref{eq:ESTIMATEONPROFILEOFMASS1}), we can perform the chain of upper estimates
\[
\begin{aligned}
u(x,t) &\leq \overline{u}(x,t) = e^{f'(0)t}\overline{v}(x,\tau) \\
&\leq e^{f'(0)t}B_M(x,\theta + \tau) =  e^{f'(0)t} M^{1 +\alpha\widehat{\gamma}}(\theta + \tau)^{-\alpha} F_1\big(x (M^{-\widehat{\gamma}}(\theta +\tau))^{-\alpha/N}\big) \\
& \leq e^{f'(0)t} K_1 M^{1 +\alpha\widehat{\gamma}}(\theta + \tau)^{-\alpha} (M^{-\widehat{\gamma}}(\theta +\tau))^{\frac{\alpha p}{N\widehat{\gamma} }} |x|^{-p/\widehat{\gamma}} \\
& \leq K e^{f'(0)t}|x|^{-p/\widehat{\gamma}},
\end{aligned}
\]
where we set $K := K_1 (2\tau_{\infty})^{1/\widehat{\gamma}}$ and we used the first relation in \eqref{eq:RELATIONMASSESBARENBLATTSOLINTRO} in the third inequality. Note that we used that
$1 + \alpha = \alpha p/N$, too. Now, supposing $|x| \geq e^{\sigma t}$ in the last inequality, we get
\[
u(x,t) \leq K e^{(f'(0) - p\sigma/\widehat{\gamma} )t} \to 0 \quad \text{in } \{ |x| \geq e^{\sigma t}\} \text{ as } t \to \infty,
\]
since we have chosen $\sigma > \sigma_{\ast}$, completing the proof. $\Box$
%
%
%
%
%
%
%
%
%
%
%
%
\section{Existence of expanding super-level sets}\label{SECTIONEXPONENTIALEXPANSIONSUPERLEVELSETS}
This section is devoted to prove that for all $\sigma < \sigma_{\ast} := \widehat{\gamma}f'(0)/p$ and initial data ``small enough'', the solution of problem \eqref{eq:REACTIONDIFFUSIONEQUATIONPLAPLACIAN} lifts up to a (small) positive constant on the ``inner sets'' $\{|x| \leq e^{\sigma t}\}$ for large times.

\noindent Let $\widetilde{\varepsilon}$ and $\widetilde{\varrho}_0$ be positive real numbers and, for all $0 < \widehat{\gamma} < p/N$, consider the initial datum
\begin{equation}\label{eq:INITIALDATUMTESTFASTLEVELSETS}
\widetilde{u}_0(x) :=
\begin{cases}
\begin{aligned}
\widetilde{\varepsilon}  \;\,\qquad\qquad &\text{if }|x| \leq \widetilde{\varrho}_0 \\
a_0|x|^{-p/\widehat{\gamma}} \quad &\text{if }|x| > \widetilde{\varrho}_0,
\end{aligned}
\end{cases}
\end{equation}
where $a_0 := \widetilde{\varepsilon}\,\widetilde{\varrho}_0^{\,p/\widehat{\gamma}}$. Note that $\widetilde{u}_0(\cdot)$ has ``tails'' which are asymptotic to the profile of the Barenblatt solutions for $|x|$ large (see formula \eqref{eq:ESTIMATEONPROFILEOFMASS1}). The choice \eqref{eq:INITIALDATUMTESTFASTLEVELSETS} will be clear in the next sections, where we will show the convergence of the solution of problem \eqref{eq:REACTIONDIFFUSIONEQUATIONPLAPLACIAN}, \eqref{eq:ASSUMPTIONSONTHEINITIALDATUMFAST} to the steady state $u = 1$. The nice property of the initial datum \eqref{eq:INITIALDATUMTESTFASTLEVELSETS} is that it can be employed as initial ``sub-datum'' as the following Lemma shows.
\begin{lem}\label{LEMMAPLACINGBARENBLATTUNDERSOLUTION}
Fix $N \geq 1$ and let $m > 0$ and $p > 1$ such that $0 < \widehat{\gamma} < p/N$. Then for all $\theta > 0$, there exist $t_1 > \theta$, $\widetilde{\varepsilon} > 0$, and $\widetilde{\varrho}_0 > 0$, such that the solution $u = u(x,t)$ of problem \eqref{eq:REACTIONDIFFUSIONEQUATIONPLAPLACIAN} with nontrivial initial datum $0 \leq u_0 \in L^1(\RR^N)$ satisfies
\[
u(x,t_1) \geq \widetilde{u}_0(x) \quad \text{in } \RR^N
\]
where $\widetilde{u}_0(\cdot)$ is defined in \eqref{eq:INITIALDATUMTESTFASTLEVELSETS}.
\end{lem}
\emph{Proof.} Let $u = u(x,t)$ the solution of problem \eqref{eq:REACTIONDIFFUSIONEQUATIONPLAPLACIAN} with nontrivial initial datum $0 \leq u_0 \in L^1(\RR^N)$ and consider the solution $v = v(x,t)$ of the purely diffusive Cauchy problem:
\[
\begin{cases}
\begin{aligned}
\partial_t v = \Delta_pv^m \;\quad\quad &\text{in } \RR^N\times(0,\infty) \\
v(x,0) = u_0(x) \quad &\text{in } \RR^N.
\end{aligned}
\end{cases}
\]
It satisfies $v(x,t) \leq u(x,t)$ in $\RR^N\times[0,\infty)$ thanks to the Maximum Principle.

Let $\theta > 0$. Since $v(\cdot,\theta)$ is continuous in $\RR^N$ and non identically zero (the mass of the solution is conserved in the ``good'' exponent range $0 < \widehat{\gamma} < p/N$), we have that it is strictly positive in a small ball $\overline{B}_{\varrho}(x_0)$, $x_0 \in \RR^N$ and $\varrho > 0$. Without loss of generality, we may take $x_0 = 0$. So, by continuity, we deduce $v(x,t+\theta) \geq \delta$ in $\overline{B}_{\varrho}\times[0,\tau]$, for some small $\delta > 0$ and $\tau > 0$.

Now, let us consider the function $v_{\theta}(x,t) := v(x,t + \theta)$ and the exterior cylinder
\[
S:= \{|x| \geq \varrho\}\times(0,\tau).
\]
We compare $v_{\theta}(x,t)$ with a ``small'' Barenblatt solution $B_M(x,t)$ (we mean that $M$ is small) at time $t=0$ and on the boundary of $S$. Recall that Barenblatt solutions have the self-similar form
\[
B_M(x,t) = t^{-\alpha}\Big[C_M + k\big|xt^{-\frac{\alpha}{N}}\big|^{\frac{p}{p-1}} \Big]^{-\frac{p-1}{\widehat{\gamma}}},
\]
where $\alpha$ and $k$ are positive constants defined in Subsection \ref{PRELIMINARIESINTRO}, and $C_M > 0$ depends on the mass.

The comparison at time $t = 0$ is immediate since $v_{\theta}(x,0) \geq 0$ and $B_M(x,0) = 0$ for all $|x| \geq \varrho$. Now, let us take $|x| = \varrho$. We want to show that $v_{\theta}(|x|=\varrho,t) \geq B_M(|x| = \varrho,t)$ for all $0 \leq t < \tau$. A simple computation shows that we can rewrite the Barenblatt solution as
\[
B_M(|x| = \varrho,t) = \frac{t^{\frac{1}{\widehat{\gamma}}}}{\Big[C_Mt^{\frac{\alpha p}{N(p-1)}} + k\varrho^{\frac{p}{p-1}}\Big]^{\frac{p-1}{\widehat{\gamma}}}} \leq \bigg( \frac{\tau}{k^{p-1}\varrho^p} \bigg)^{\frac{1}{\widehat{\gamma}}},
\]
since $0 \leq t < \tau$. Thus, since $v_{\theta}(|x|=\varrho,t) \geq \delta$, it is sufficient to have
\[
\bigg( \frac{\tau}{k^{p-1}\varrho^p} \bigg)^{\frac{1}{\widehat{\gamma}}} \leq \delta, \quad \text{ i.e. } \quad \tau \leq \delta^{\widehat{\gamma}}k^{p-1}\varrho^p.
\]
This condition is satisfied taking $\tau > 0$ small enough. We may now use the Maximum Principle to obtain the conclusion $v_{\theta}(x,t) \geq B_M(x,t)$ in the whole of $S$.

In particular, we evaluate the comparison at $t = \tau$, and we have
\[
\begin{cases}
\begin{aligned}
v(x,\tau + \theta) \geq \delta \,\;\qquad\quad\qquad &\text{if } |x| \leq \varrho \\
v(x,\tau + \theta) \geq B_M(x,\tau) \qquad &\text{if } |x| \geq \varrho.
\end{aligned}
\end{cases}
\quad \Leftrightarrow \quad
\begin{cases}
\begin{aligned}
v(x,t_1) \geq \delta \,\;\quad\qquad\qquad\qquad &\text{if } |x| \leq \varrho \\
v(x,t_1) \geq B_M(x,t_1 - \theta) \qquad &\text{if } |x| \geq \varrho,
\end{aligned}
\end{cases}
\]
where we set $t_1 = \tau + \theta$. Let us fix $\widetilde{\varrho}_0 := \varrho$, $0 < \widetilde{\varepsilon} \leq \delta$, and $a_0 := \widetilde{\varepsilon}\, \widetilde{\varrho}_0^{\,p/\widehat{\gamma}}$. By taking $\widetilde{\varepsilon} > 0$ smaller, we can assume $k^{p-1} a_0^{\widehat{\gamma}} < \tau = t_1 -\theta$. Now, we verify that
\[
B_M(x,t_1-\theta) \geq a_0|x|^{-p/\widehat{\gamma}}, \quad \text{ for all } |x| \geq \widetilde{\varrho}_0,
\]
and some suitable constant $C_M > 0$. Writing the expression for the Barenblatt solutions, the previous inequality reads:
\[
C_M \leq \frac{(t_1-\theta)^{\frac{1}{p-1}} - ka_0^{\frac{\widehat{\gamma}}{p-1}} }{\big[a_0^{\widehat{\gamma}}(t_1-\theta)^{\frac{\alpha p}{N}}\big]^{\frac{1}{p-1}}} |x|^{\frac{p}{p-1}} := K|x|^{\frac{p}{p-1}}, \quad \text{ for all } |x| \geq \widetilde{\varrho}_0.
\]
Note that the coefficient $K$ of $|x|^{p/(p-1)}$ is positive thanks to our assumptions on $\widetilde{\varepsilon} > 0$. Now, since $|x| \geq \widetilde{\varrho}_0$, we deduce that a sufficient condition so that the previous inequality is satisfied is $C_M \leq K\widetilde{\varrho}_0^{\,p/(p-1)}$. Consequently, we have shown that for all $\theta > 0$, there exist $t_1 > \theta$, $\widetilde{\varepsilon} > 0$, and $\widetilde{\varrho}_0 > 0$, such that
\[
u(x,t_1) \geq v(x,t_1) \geq \widetilde{u}_0(x), \quad \text{ for all } x \in \RR^N,
\]
which is our thesis. $\Box$

\bigskip

We ask the reader to note that improved global positivity estimates were proved in \cite{Her-Pier:art,  Vascppme} and \cite{Bon-Vaz1:art} for the Porous Medium Equation. Now, with the next crucial lemma, we prove that the expansion of the super-level sets of the solution $u = u(x,t)$ of problem \eqref{eq:REACTIONDIFFUSIONEQUATIONPLAPLACIAN} with initial datum \eqref{eq:INITIALDATUMTESTFASTLEVELSETS} is exponential for all $\sigma < \sigma_{\ast}$ and large times.
\begin{lem}\label{LEMMAEXPANDINGFASTLEVELSETS}
Fix $N \geq 1$. Let $m > 0$ and $p > 1$ such that $0 < \widehat{\gamma} < p/N$, and let $0 < \sigma < \sigma_{\ast}$.
\\
Then there exist $t_0 > 0$ and $0 < \widetilde{\varepsilon}_0 < 1$ which depend only on $m$, $p$, $N$ and $f$, such that the following hold. For all $0 < \widetilde{\varepsilon} \leq \widetilde{\varepsilon}_0$, there exists $\widetilde{\varrho}_0 > 0$ (large enough depending on $\widetilde{\varepsilon} > 0$), such that the solution $u = u(x,t)$ of problem \eqref{eq:REACTIONDIFFUSIONEQUATIONPLAPLACIAN} with initial datum
\eqref{eq:INITIALDATUMTESTFASTLEVELSETS} satisfies
\[
u(x,jt_0) \geq \widetilde{\varepsilon} \quad \text{in } \{|x| \leq \widetilde{\varrho}_0 e^{\sigma j t_0}\}, \text{ for all } \; j \in \NN_+ = \{1,2,\dots\}.
\]
\end{lem}
\emph{Proof.} We prove the assertion of the thesis by induction on $j = 1,2,\dots$. We follow the ideas presented by Cabr\'e and Roquejoffre in \cite{C-R2:art} and, later, in \cite{S-V:art}, for fractional diffusion.

\emph{Step0.} We set $j = 1$, $0 < \sigma < \sigma_{\ast}$ and introduce some basic definitions and quantities we will use during the proof. First of all, let $C_1$ be the constant corresponding to the profile $F_1(\cdot)$ (see Section \ref{PRELIMINARIESINTRO}) and let $K_1$ and $K_2$ be defined as in \eqref{eq:ESTIMATEONPROFILEOFMASS1} with $M = 1$. In order to avoid huge expressions in the following of the proof, we introduce the constants
\begin{equation}\label{eq:DEFINITIONSOFCONSTANTSKFASTLEVELSETS}
\overline{K}:= \big(C_1^{(p-1)/\widehat{\gamma}} K_1^{-\alpha\widehat{\gamma}} \big)^{N/(\alpha p)} \quad \text{ and } \quad \widetilde{K} := \frac{K_2}{2} C_1^{\frac{p-1}{\widehat{\gamma}}}.
\end{equation}
We fix $0 < \delta < 1$ sufficiently small such that
\begin{equation}\label{eq:CONDITIONONDELTAFASTLEVELSETS}
\frac{\widehat{\gamma}}{p}\lambda > \sigma, \qquad\quad \lambda := f(\delta)/\delta.
\end{equation}
Then, we consider $t_0$ sufficiently large such that
\begin{equation}\label{eq:CONDITIONONT0FASTLEVELSETS}
\widetilde{K} e^{\lambda t_0} \geq 2^{\alpha} \quad \text{ and } \quad \frac{K_2}{2K_1}e^{\lambda t_0} \geq e^{\frac{p}{\widehat{\gamma}} \sigma t_0}
\end{equation}
(note that such a $t_0$ exists thanks to \eqref{eq:CONDITIONONDELTAFASTLEVELSETS}) and we define $\widetilde{\varepsilon}_0 := \delta e^{-\lambda t_0}$. Finally, fix $0 < \widetilde{\varepsilon} \leq \widetilde{\varepsilon}_0$ and choose $\widetilde{\varrho}_0$ large enough such that
\begin{equation}\label{eq:CONDITIONONVARRHO0FASTLEVELSETS}
\widetilde{\varrho}_0^{\,p} \geq \frac{K_1^{\widehat{\gamma}}}{\lambda \widehat{\gamma} \widetilde{\varepsilon}^{\,\widehat{\gamma}}}.
\end{equation}
We anticipate that the choice of the subtle conditions \eqref{eq:CONDITIONONDELTAFASTLEVELSETS}, \eqref{eq:CONDITIONONT0FASTLEVELSETS} and \eqref{eq:CONDITIONONVARRHO0FASTLEVELSETS} will be clarified during the proof.

\emph{Step1.} Construction of a sub-solution of problem \eqref{eq:REACTIONDIFFUSIONEQUATIONPLAPLACIAN}, \eqref{eq:INITIALDATUMTESTFASTLEVELSETS} in $\RR^N\times[0,t_0]$. First of all, we construct a Barenblatt solution of the form $B_{M_1}(x,\theta_1)$ such that
\begin{equation}\label{eq:SUBBARENBLATTINITIALDATAFASTLEVELSETS}
B_{M_1}(x,\theta_1) \leq \widetilde{u}_0(x) \quad \text{in } \RR^N.
\end{equation}
Since the profile of the Barenblatt solution is decreasing, we impose $B_{M_1}(0,\theta_1) = \widetilde{\varepsilon}$ in order to satisfy \eqref{eq:SUBBARENBLATTINITIALDATAFASTLEVELSETS} in the set $\{|x| \leq \widetilde{\varrho}_0\}$. Moreover, using \eqref{eq:ESTIMATEONPROFILEOFMASS1} and noting that $1 + \alpha\widehat{\gamma} = \alpha p/N$, it simple to get
\[
B_{M_1}(x,\theta_1) \leq K_1 \theta_1^{\frac{1}{\widehat{\gamma}}} |x|^{-\frac{p}{\widehat{\gamma}}} \quad \text{in } \RR^N
\]
and so, it is sufficient to require $K_1 \theta_1^{\frac{1}{\widehat{\gamma}}} = a_0$, so that \eqref{eq:SUBBARENBLATTINITIALDATAFASTLEVELSETS} is valid in $\{|x| \geq \widetilde{\varrho}_0\}$. Thus, it is simple to obtain the relations
\begin{equation}\label{eq:CHOICEOFBARENBLATTPARAMETERSFASTLEVELSETS}
M_1 = \overline{K} \widetilde{\varrho}_0^{\,N} \widetilde{\varepsilon} \qquad \text{ and } \qquad \theta_1 = K_1^{-\widehat{\gamma}} \widetilde{\varrho}_0^{\,p} \widetilde{\varepsilon}^{\,\widehat{\gamma}}
\end{equation}
Now, consider the linearized problem
\begin{equation}\label{eq:FIRSTLINEARIZEDPROBLEMFASTLEVELSETS}
\begin{cases}
\begin{aligned}
\partial_tw = \Delta_p w^m + \lambda w \quad &\text{in } \RR^N\times(0,\infty)\\
w(x,0) = \widetilde{u}_0(x) \;\;\quad\quad &\text{in } \RR^N
\end{aligned}
\end{cases}
\end{equation}
and the change of variable
\begin{equation}\label{eq:FIRSTCHANGEOFVARIABLEFASTLEVELSETS}
\tau(t) = \frac{1}{\lambda\widehat{\gamma}}\Big[ 1 - e^{-\lambda\widehat{\gamma} t} \Big], \quad \text{for } t \geq 0.
\end{equation}
Note that $0 \leq \tau(t) \leq \tau_{\infty} := \frac{1}{\lambda\widehat{\gamma}}$ and the function $\widetilde{w}(x,\tau) = e^{-\lambda t}w(x,t)$ solves the ``pure diffusive'' problem
\begin{equation}\label{eq:PUREDIFFUSIVEPROBLEMMINLEVELSETFAST}
\begin{cases}
\begin{aligned}
\partial_{\tau}\widetilde{w} = \Delta_p \widetilde{w}^m  \qquad &\text{in } \RR^N\times(0,\tau_{\infty})\\
w(x,0) = \widetilde{u}_0(x) \quad &\text{in } \RR^N.
\end{aligned}
\end{cases}
\end{equation}
Since $B_{M_1}(x,\theta_1) \leq \widetilde{u}_0(x) \leq \widetilde{\varepsilon}$ for all $x \in \RR^N$, from the Maximum Principle we get
\begin{equation}\label{eq:INEQUALITYSUBSOLUTIONFASTLEVELSETS}
B_{M_1}(x,\theta_1 + \tau) \leq \widetilde{w}(x,\tau) \leq \widetilde{\varepsilon} \quad \text{in } \RR^N\times(0,\tau_{\infty}).
\end{equation}
Hence, using the concavity of $f$ and the second inequality in \eqref{eq:INEQUALITYSUBSOLUTIONFASTLEVELSETS} we get
\[
w(x,t) = e^{\lambda t}\widetilde{w}(x,\tau) \leq \widetilde{\varepsilon}_0 e^{\lambda t_0} = \delta, \quad \text{in } \RR^N\times[0,t_0]
\]
and so, since $w \leq \delta$ implies $f(\delta)/\delta \leq f(w)/w$, we have that $w$ is a sub-solution of problem \eqref{eq:REACTIONDIFFUSIONEQUATIONPLAPLACIAN}, \eqref{eq:INITIALDATUMTESTFASTLEVELSETS} in $\RR^N\times[0,t_0]$. Finally, using the first inequality in \eqref{eq:INEQUALITYSUBSOLUTIONFASTLEVELSETS}, we obtain
\begin{equation}\label{eq:FIRSTINEQUALITYFORCHOOSINGEPSFASTLEVELSETS}
u(x,t) \geq e^{\lambda t} \widetilde{w}(x,\tau) \geq e^{\lambda t} B_{M_1}(x,\theta_1 + \tau) \quad \text{in } \RR^N\times[0,t_0].
\end{equation}

\emph{Step2.} In this step, we show that the choices made in  \eqref{eq:CONDITIONONDELTAFASTLEVELSETS}, \eqref{eq:CONDITIONONT0FASTLEVELSETS} and \eqref{eq:CONDITIONONVARRHO0FASTLEVELSETS} allow us to find positive numbers $\widetilde{\varrho}_1$ and $a_1$ such that $u(x,t_0) \geq \widetilde{u}_1(x)$ for all $x \in \RR^N$, where
\[
\widetilde{u}_1(x) :=
\begin{cases}
\begin{aligned}
\widetilde{\varepsilon} = a_1\widetilde{\varrho}_1^{\,-p/\widehat{\gamma}}  \quad &\text{if }|x| \leq \widetilde{\varrho}_1 \\
a_1|x|^{-p/\widehat{\gamma}} \,\;\quad\quad &\text{if }|x| > \widetilde{\varrho}_1
\end{aligned}
\end{cases}
\quad \text{and } \qquad \widetilde{\varrho}_1 \geq \widetilde{\varrho}_0 e^{\sigma t_0},
\]
which implies the thesis for $j = 1$. Now, in order to find $\widetilde{\varrho}_1$ and $a_1$ we proceed with the chain of inequalities in \eqref{eq:FIRSTINEQUALITYFORCHOOSINGEPSFASTLEVELSETS} for the values $t = t_0$, $\tau_0 = \tau(t_0)$ and $|x| = \widetilde{\varrho}_1$. Imposing
\begin{equation}\label{eq:FIRSTCOMPATIBILITYCHOICEFASTLEVELSETS}
\widetilde{\varrho}_1 \big[M_1^{-\widehat{\gamma}}(\theta_1 + \tau_0)\big]^{-\alpha/N} \geq 1,
\end{equation}
using \eqref{eq:ESTIMATEONPROFILEOFMASS1} and observing that $(1+z)^{-1} \geq (2z)^{-1}$ for all $z\geq 1$, we look for $\widetilde{\varrho}_1$ and $a_1$ such that
\[
\begin{aligned}
u(x,t_0)|_{|x| = \widetilde{\varrho}_1} & \geq e^{\lambda t_0} B_{M_1}(x,\theta_1 + \tau_0)|_{|x| = \widetilde{\varrho}_1}
\\
& \geq e^{\lambda t_0} K_2 M_1^{1+\alpha\widehat{\gamma}}(\theta_1 + \tau_0)^{-\alpha}\Big\{1 + \widetilde{\varrho}_1^{\,\frac{p}{\widehat{\gamma}}}\big[M_1^{-\widehat{\gamma}}(\theta_1 + \tau_0)\big]^{-\frac{\alpha p}{\;N\widehat{\gamma}}} \Big\}^{-1} \\
& \geq \frac{K_2}{2}e^{\lambda t_0} M_1^{1+\alpha\widehat{\gamma}}(\theta_1 + \tau_0)^{-\alpha}\Big\{ \widetilde{\varrho}_1 \big[M_1^{-\widehat{\gamma}}(\theta_1 + \tau_0)\big]^{-\frac{\alpha }{N}} \Big\}^{-\frac{p}{\widehat{\gamma}}}
\\
& = \frac{K_2}{2}e^{\lambda t_0} (\theta_1 + \tau_0)^{\frac{1}{\widehat{\gamma}}}\widetilde{\varrho}_1^{-\frac{p}{\widehat{\gamma}}}
\\
&\geq \widetilde{\varepsilon} = a_1\widetilde{\varrho}_1^{-\frac{p}{\widehat{\gamma}}}.
\end{aligned}
\]
Thus, we get $u(x,t_0)|_{|x| = \widetilde{\varrho}_1} \geq \widetilde{\varepsilon}$ taking, for instance, $\widetilde{\varrho}_1 > 0$ such that
\begin{equation}\label{eq:SECONDCOMPATIBILITYCHOICEFASTLEVELSETS}
\frac{K_2}{2}e^{\lambda t_0} (\theta_1 + \tau_0)^{\frac{1}{\widehat{\gamma}}}\widetilde{\varrho}_1^{-\frac{p}{\widehat{\gamma}}} = \widetilde{\varepsilon}.
\end{equation}
Note that this choice of $\widetilde{\varrho}_1 > 0$ performs the equality at the end of the previous chain and the value $a_1 = \frac{K_2}{2}e^{\lambda t_0} (\theta_1 + \tau_0)^{1/\widehat{\gamma}}$ is determined too.

\emph{Remark 1}. Note that the conditions $e^{\lambda t_0}B_{M_1}(x,\theta_1 + \tau_0)|_{|x| = \widetilde{\varrho}_1} \geq \widetilde{\varepsilon}$ and \eqref{eq:FIRSTCOMPATIBILITYCHOICEFASTLEVELSETS} are sufficient to assure $u(x,t_0) \geq u_1(x)$ in $\RR^N$. Indeed, it is guaranteed in the set $\{|x| \leq \widetilde{\varrho}_1\}$ since the profile $F_1(\cdot)$ is non-increasing. On the other hand, if $|x| \geq \widetilde{\varrho}_1$ we have $|x| \big[M_1^{-\widehat{\gamma}}(\theta_1 + \tau_0)\big]^{-\alpha/N} \geq 1$ by \eqref{eq:FIRSTCOMPATIBILITYCHOICEFASTLEVELSETS} and so, following the chain of inequalities as before, we get
\[
u(x,t_0) \geq \frac{K_2}{2}e^{\lambda t_0} (\theta_1 + \tau_0)^{\frac{1}{\widehat{\gamma}}}|x|^{-\frac{p}{\widehat{\gamma}}} = a_1 |x|^{-\frac{p}{\widehat{\gamma}}} = \widetilde{u}_1(x) \quad \text{in } \{|x| \geq \widetilde{\varrho}_1\}.
\]

Now, in order to conclude the proof of the case $j = 1$, we must check that the conditions \eqref{eq:FIRSTCOMPATIBILITYCHOICEFASTLEVELSETS} and \eqref{eq:SECONDCOMPATIBILITYCHOICEFASTLEVELSETS}  actually represent a possible choice and the value of $t_0$, defined at the beginning, performs their compatibility. The compatibility between \eqref{eq:FIRSTCOMPATIBILITYCHOICEFASTLEVELSETS} and \eqref{eq:SECONDCOMPATIBILITYCHOICEFASTLEVELSETS} can be verified imposing
\[
\frac{K_2}{2}e^{\lambda t_0} (\theta_1 + \tau_0)^{\frac{1}{\widehat{\gamma}}} = \widetilde{\varepsilon}\, \widetilde{\varrho}_1^{\,\frac{p}{\widehat{\gamma}}} \geq \widetilde{\varepsilon} \big[M_1^{-\widehat{\gamma}}(\theta_1 + \tau_0)\big]^{\frac{\alpha p}{\;N\widehat{\gamma}}},
\]
which can be rewritten using the definitions \eqref{eq:CHOICEOFBARENBLATTPARAMETERSFASTLEVELSETS} as
\begin{equation}\label{eq:CHOICEOFT0CONDITION1FASTLEVELSETS}
\widetilde{K} e^{\lambda t_0} \geq \Bigg(1 + \frac{\tau_0}{\theta_1}\Bigg)^{\alpha},
\end{equation}
Now, it is simple to verify that condition \eqref{eq:CONDITIONONVARRHO0FASTLEVELSETS} implies $\tau_{\infty} \leq \theta_1$ and so it holds $\tau_0 \leq \theta_1$ too. Hence, a sufficient condition so that \eqref{eq:CHOICEOFT0CONDITION1FASTLEVELSETS} is satisfied and does not depend on $\widetilde{\varepsilon} > 0$ is
\[
\widetilde{K} e^{\lambda t_0} \geq 2^{\alpha},
\]
i.e., our initial choice of $t_0$ in \eqref{eq:CONDITIONONT0FASTLEVELSETS} which proves the compatibility between \eqref{eq:FIRSTCOMPATIBILITYCHOICEFASTLEVELSETS} and \eqref{eq:SECONDCOMPATIBILITYCHOICEFASTLEVELSETS}.

\emph{Remark 2}. Rewriting formula \eqref{eq:SECONDCOMPATIBILITYCHOICEFASTLEVELSETS} using the definition of $\theta_1$, it is simple to deduce
\begin{equation}\label{THIRDCOMPATIBILITYCHOICEFASTLEVELSETS}
\Bigg( \frac{\widetilde{\varrho}_1}{\widetilde{\varrho}_0}\Bigg)^{\frac{p}{\widehat{\gamma}}} = \frac{K_2}{2K_1}e^{\lambda t_0}\Bigg(1 + \frac{\tau_0}{\theta_1}\Bigg)^{\frac{1}{\widehat{\gamma}}}
\end{equation}
and, using the second hypothesis on $t_0$ in \eqref{eq:CONDITIONONT0FASTLEVELSETS}, it is straightforward to obtain $\widetilde{\varrho}_1 \geq \widetilde{\varrho}_0 e^{\sigma t_0}$. In particular, we have shown
\[
u(x,t_0) \geq \widetilde{\varepsilon} \quad \text{in } \{|x| \leq \widetilde{\varrho}_0e^{\sigma t_0} \},
\]
i.e., the thesis for $j = 1$.
\paragraph{Iteration.} Set $t_j := (j+1)t_0$, $\widetilde{\varrho}_j := \widetilde{\varrho}_0e^{\sigma jt_0}$ and $a_j := \widetilde{\varepsilon}\,\widetilde{\varrho}_j^{\,p/\widehat{\gamma}}$ for all $j \in \NN$ and define
\begin{equation}\label{eq:KINITIALDATUMTESTFASTLEVELSETS}
\widetilde{u}_j(x) =
\begin{cases}
\begin{aligned}
\widetilde{\varepsilon}  \;\,\qquad\qquad &\text{if }|x| \leq \widetilde{\varrho}_j \\
a_j |x|^{-p/\widehat{\gamma}} \quad &\text{if }|x| > \widetilde{\varrho}_j.
\end{aligned}
\end{cases}
\end{equation}
We suppose to have proved that the solution of problem \eqref{eq:REACTIONDIFFUSIONEQUATIONPLAPLACIAN}, \eqref{eq:INITIALDATUMTESTFASTLEVELSETS} satisfies
\[
u(x,t_{j-1}) \geq \widetilde{u}_j(x) \quad \text{in } \RR^N, \quad\text{for some }  j \in \NN_+
\]
and we show $u(x,t_{j}) \geq \widetilde{u}_{j+1}(x)$ in $\RR^N$ for the values $\widetilde{\varrho}_{j+1}$ and $a_{j+1}$. From the induction hypothesis, we have that the solution $v(x,t)$ of the problem
\begin{equation}\label{eq:KCAUCHYPROBLEMFASTLEVELSETS}
\begin{cases}
\begin{aligned}
\partial_tv = \Delta_pv^m + f(v) \quad &\text{in } \RR^N\times(t_{j-1},\infty)\\
v(x,t_{j-1}) = \widetilde{u}_j(x) \;\;\quad &\text{in } \RR^N
\end{aligned}
\end{cases}
\end{equation}
is a sub-solution of problem \eqref{eq:REACTIONDIFFUSIONEQUATIONPLAPLACIAN}, \eqref{eq:INITIALDATUMTESTFASTLEVELSETS} in $\RR^N\times[t_{j-1},\infty)$ which implies $u(x,t) \geq v(x,t)$ in $\RR^N\times[t_{j-1},\infty)$ and so, it is sufficient to prove $v(x,t_{j}) \geq \widetilde{u}_{j+1}(x)$ in $\RR^N$. Since we need to repeat almost the same procedure of the case $j = 1$, we only give a brief sketch of the induction step.

\emph{Step1'}. Construction of a sub-solution of problem \eqref{eq:KCAUCHYPROBLEMFASTLEVELSETS}, \eqref{eq:KINITIALDATUMTESTFASTLEVELSETS}, in $\RR^N\times[t_{j-1},t_j]$. With the same techniques used in \emph{Step1}, we construct a Barenblatt solution $B_{M_{j+1}}(x,\theta_{j+1}) \leq \widetilde{u}_j(x)$ in $\RR^N$ with parameters
\begin{equation}\label{eq:CHOICEOFBARENBLATTPARAMETERSKFASTLEVELSETS}
M_{j+1} = \overline{K} \widetilde{\varrho}_j^{\,N} \widetilde{\varepsilon} \qquad \text{ and } \qquad \theta_{j+1} = K_1^{-\widehat{\gamma}} \widetilde{\varrho}_j^{\,p} \widetilde{\varepsilon}^{\,\widehat{\gamma}}
\end{equation}
and a sub-solution of problem \eqref{eq:KCAUCHYPROBLEMFASTLEVELSETS}, \eqref{eq:KINITIALDATUMTESTFASTLEVELSETS}: $w(x,t) = e^{\lambda (t- t_{j-1})} \widetilde{w}(x,\widetilde{\tau})$ in $\RR^N\times[t_{j-1},t_j]$, where
\[
\widetilde{\tau}(t) = \frac{1}{\lambda\widehat{\gamma}}\Big[ 1 - e^{-\lambda\widehat{\gamma}(t-t_{j-1})} \Big], \quad \text{for } t \geq t_{j-1}.
\]
In particular, note that $\theta_{j+1} \geq \theta_j \geq \ldots \geq \theta_1$, $0 \leq \widetilde{\tau}(t) \leq \tau_{\infty}$, $\widetilde{\tau}(t_j) := \widetilde{\tau}_j = \tau_0$ and
\[
v(x,t_j) \geq e^{\lambda t_0} B_{M_{j+1}}(x,\theta_{j+1} + \widetilde{\tau}_j).
\]

\emph{Step2'}. We have to study a chain of inequalities similar to the one carried out in \emph{Step2} verifying that
\[
e^{\lambda t_0} B_{M_{j+1}}(x,\theta_{j+1} + \widetilde{\tau}_j)|_{|x| = \widetilde{\varrho}_{j+1}} \geq \widetilde{\varepsilon}.
\]
Thus, imposing conditions similar to \eqref{eq:FIRSTCOMPATIBILITYCHOICEFASTLEVELSETS} and \eqref{eq:SECONDCOMPATIBILITYCHOICEFASTLEVELSETS} and requiring their compatibility, we have to check the validity of the inequality
\[
\widetilde{K} e^{\lambda t_0} \geq \Bigg(1 + \frac{\tau_0}{\theta_{j+1}}\Bigg)^{\alpha}.
\]
Since $\theta_1 \leq \theta_{j+1}$, we have $\tau_0 \leq \theta_{j+1}$ and so, a sufficient condition so that the previous inequality is satisfied is $\widetilde{K} e^{\lambda t_0} \geq 2^{\alpha}$, which is guaranteed by the initial choice of $t_0$. Finally, following the reasonings of the case $j = 1$ it is simple to obtain the relation
\[
\Bigg( \frac{\widetilde{\varrho}_{j+1}}{\widetilde{\varrho}_{j}}\Bigg)^{\frac{p}{\widehat{\gamma}}} \geq \frac{K_2}{2K_1}e^{\lambda t_0}\Bigg(1 + \frac{\tau_0}{\theta_{j+1}}\Bigg)^{\frac{1}{\widehat{\gamma}}}
\]
which implies
\[
\widetilde{\varrho}_{j+1} \geq \widetilde{\varrho}_j e^{\sigma t_0} \geq \ldots \geq \widetilde{\varrho}_0 e^{\sigma j t_0},
\]
and we complete the proof. $\Box$

\bigskip

\noindent {\bf Proof of Proposition \ref{EXPANPANSIONOFMINIMALLEVELSETS}.} The previous lemma proves that for the sequence of times $t_j = (jt_0)_{j \in \NN_+}$ and for any choice of the parameter $0 < \sigma < \sigma_{\ast}$, the solution of problem \eqref{eq:REACTIONDIFFUSIONEQUATIONPLAPLACIAN}, \eqref{eq:INITIALDATUMTESTFASTLEVELSETS} reaches a positive value $\widetilde{\varepsilon}$ in the sequence of sets $ \{|x| \leq \widetilde{\varrho}_0e^{\sigma j t_0}\}$ where $\widetilde{\varrho}_0 > 0$ is chosen large enough (in particular, we can assume $\widetilde{\varrho}_0 \geq 1$).

\noindent Actually, we obtained a more useful result. First of all, note that, for all $0 < \sigma < \sigma_{\ast}$, Lemma \ref{LEMMAEXPANDINGFASTLEVELSETS} implies
\[
u(x,jt_0) \geq \widetilde{\varepsilon} \quad \text{in } \{|x| \leq e^{\sigma j t_0}\}, \text{ for all } \; j \in \NN_+,
\]
for all $0 < \widetilde{\varepsilon} \leq \widetilde{\varepsilon}_0 = \delta e^{-f'(0)t_0}$.
Moreover, since conditions \eqref{eq:CONDITIONONT0FASTLEVELSETS} are satisfied for all $t_0 \leq t_1 \leq 2t_0$, we can repeat the same proof of Lemma \ref{LEMMAEXPANDINGFASTLEVELSETS}, modifying the value of $\widetilde{\varepsilon}_0$ and choosing a different value $\underline{\widetilde{\varepsilon}}_0 = \delta e^{-2f'(0)t_0} > 0$, which is smaller but strictly positive for all $t_0 \leq t_1 \leq 2t_0$. Hence, it turns out that for all $0 < \widetilde{\varepsilon} \leq \underline{\widetilde{\varepsilon}}_0$, it holds
\[
u(x,t) \geq \widetilde{\varepsilon} \quad \text{in } \{ |x| \leq e^{\sigma t}\},  \text{ for all } \; t_0 \leq t \leq 2t_0.
\]
Now, iterating this procedure as in the proof of Lemma \ref{LEMMAEXPANDINGFASTLEVELSETS}, it is clear that we do not have to change the value of $\underline{\widetilde{\varepsilon}}_0$ when $j \in \NN_+$ grows and so, for all $0 < \widetilde{\varepsilon} \leq \underline{\widetilde{\varepsilon}}_0$, we obtain
\[
u(x,t) \geq \widetilde{\varepsilon} \quad \text{in } \{|x| \leq e^{\sigma t}\}, \text{ for all } \; j \in \NN_+ \;\text{ and for all }\;\; j t_0 \leq t \leq (j+1)t_0.
\]
Then, using the arbitrariness of $j \in \NN_+$, we complete the proof. $\Box$
\paragraph{Remark.} Note that, to be precise, in the proof of Proposition \ref{EXPANPANSIONOFMINIMALLEVELSETS}, we have to combine Lemma \ref{LEMMAPLACINGBARENBLATTUNDERSOLUTION} with Lemma \ref{LEMMAEXPANDINGFASTLEVELSETS} as follows. Let $u = u(x,t)$ the solution of problem \eqref{eq:REACTIONDIFFUSIONEQUATIONPLAPLACIAN} with initial datum \eqref{eq:ASSUMPTIONSONTHEINITIALDATUMFAST}. We wait a time $t_1 > 0$ given by Lemma \ref{LEMMAPLACINGBARENBLATTUNDERSOLUTION}, in order to have
\[
u(x,t_1) \geq \widetilde{u}_0(x) \quad \text{in } \RR^N,
\]
for all $\widetilde{\varrho}_0 > 0$ and some $\widetilde{\varepsilon} > 0$ depending on $t_1$. Now, thanks to the Maximum Principle, we deduce $u(x,t + t_1) \geq \widetilde{u}(x,t)$ in $\RR^N\times[0,\infty)$, where we indicate with $\widetilde{u} = \widetilde{u}(x,t)$ the solution of problem \eqref{eq:REACTIONDIFFUSIONEQUATIONPLAPLACIAN} with initial datum $\widetilde{u}_0 = \widetilde{u}_0(x)$. In this way, we deduce the statement of Lemma \ref{LEMMAEXPANDINGFASTLEVELSETS} for more general initial data satisfying \eqref{eq:ASSUMPTIONSONTHEINITIALDATUMFAST} and we can prove Proposition \ref{EXPANPANSIONOFMINIMALLEVELSETS}.
%
%
%
%
%
%
%
%
%
%
%
%
\section{Convergence to 1 in the ``inner sets''}\label{SECTIONASYMPTOTICBEHAVIOURFAST}
As mentioned in the introduction, we now address to the problem of showing the convergence of a general solution of problem \eqref{eq:REACTIONDIFFUSIONEQUATIONPLAPLACIAN}, \eqref{eq:ASSUMPTIONSONTHEINITIALDATUMFAST} to the steady state $u = 1$. As anticipated, we find that the convergence to 1 is exponential for large times, with exponent $\sigma < \sigma_{\ast}$. This fact represents an interesting deviance, respect to the case $\gamma \geq 0$ (i.e. $\widehat{\gamma} \leq 0$) in which the solutions converge with constant speed for large times and show a TW asymptotic behaviour.
\paragraph{Proof of Theorem \ref{CONVERGENCETOONEFASTDIFFUSION}.} Fix $0 < \widehat{\gamma} < p/N$, $0 < \sigma < \sigma_{\ast}$ and set $\underline{w} := 1 - u^m$ in $\RR^N\times(0,\infty)$. We will prove that for all $\varepsilon > 0$, there exists $t_{\varepsilon} > 0$ such that
\[
\underline{w}(x,t) \leq \varepsilon \quad \text{ in } \{|x|\leq e^{\sigma t}, t \geq t_{\varepsilon}\},
\]
which is equivalent to the assertion of the thesis.

\emph{Step1.} Fix $\sigma < \nu < \sigma_{\ast}$ and consider the inner set $\Omega_I := \{|x|\leq e^{\nu t}, t \geq t_1\}$, where $t_1 > 0$ is initially arbitrary. We recall that Proposition \ref{EXPANPANSIONOFMINIMALLEVELSETS} assures the existence of $\widetilde{\varepsilon} > 0$ and $t_0 > 0$ such that $u \geq \widetilde{\varepsilon}$ in the set $\{|x|\leq e^{\nu t}, t \geq t_0\}$. In particular, for all $t_1 \geq t_0$, we have that $u = u(x,t)$ is bounded from below and above in the inner set:
\begin{equation}\label{eq:APPLICATIONOFPROPLEVELSETS}
\widetilde{\varepsilon} \leq u \leq 1 \quad \text{ in } \Omega_I.
\end{equation}
Moreover, it is not difficult to see that, setting $a(x,t) = (1/m)u^{1-m}$ and $c(x,t) = f(u)/\underline{w}$, the function $\underline{w} = 1 - u^m$ solves the problem
\begin{equation}\label{eq:PROBLEMSOLVEDBYV}
\begin{cases}
\begin{aligned}
a(x,t)\partial_t\underline{w} - \Delta_p\underline{w} + c(x,t)\underline{w} = 0 \quad &\text{in } \RR^N\times(t_1,\infty) \\
\underline{w}(x,t_1) = 1 - [u(x,t_1)]^m \qquad\qquad &\text{in } \RR^N.
\end{aligned}
\end{cases}
\end{equation}
Using \eqref{eq:APPLICATIONOFPROPLEVELSETS}, it is simple to see that
\[
a_0 \leq a(x,t) \leq a_1 \quad \text{in } \Omega_I
\]
where
\[
a_0 :=
\begin{cases}
\begin{aligned}
(1/m)\widetilde{\varepsilon}^{1-m} \quad &\text{if } 0 < m < 1 \\
1/m \qquad\quad\;\;\; &\text{if } m \geq 1
\end{aligned}
\end{cases}
\quad
a_1 :=
\begin{cases}
\begin{aligned}
1/m \qquad\quad\;\;\; &\text{if } 0 < m < 1 \\
(1/m)\widetilde{\varepsilon}^{1-m} \quad &\text{if } m \geq 1.
\end{aligned}
\end{cases}
\]
For what concerns $c(x,t)$, it is bounded from below in $\Omega_I$:
\[
c(x,t) \geq c_0 \quad \text{in } \Omega_I,
\]
where $c_0 > 0$ and depends on $\widetilde{\varepsilon}$ and $m$. Indeed, if $0 < m < 1$ we have that $c(x,t) = f(u)/(1-u^m) \geq f(u)/(1-u)$ for all $0 \leq u \leq 1$. Hence, we get our bound from below recalling \eqref{eq:APPLICATIONOFPROPLEVELSETS} and noting that
\[
\frac{f(u)}{1-u} \sim -f'(1) > 0  \quad \text{as } u \sim 1.
\]
If $m \geq 1$, we have the formula
\[
c(x,t) = \frac{f(u)}{1-u^m} = \frac{f(u)}{(1-u)(1 + u + \ldots + u^{m-1})},
\]
and so, since $u \leq 1$ and arguing as in the case $0 < m < 1$, we deduce
\[
c(x,t) \geq (1/m)\frac{f(u)}{1-u} \geq c_0 \quad \text{in } \Omega_I
\]
for some $c_0 > 0$ depending on $\widetilde{\varepsilon}$ and $m$. In particular, it follows that $\underline{w} = \underline{w}(x,t)$ satisfies
\begin{equation}\label{eq:DIFFINEQSUBSOLFASTCONV1}
a(x,t)\partial_t\underline{w} - \Delta_p\underline{w} + c_0 \underline{w} \leq 0 \quad \text{in } \Omega_I,
\end{equation}
i.e., $\underline{w} = \underline{w}(x,t)$ is a sub-solution for the equation in problem \eqref{eq:PROBLEMSOLVEDBYV} in the set $\Omega_I$.

\emph{Step2.} In this step, we look for a super-solution $\overline{w} = \overline{w}(x,t)$ of problem \eqref{eq:PROBLEMSOLVEDBYV} with $\partial_t\overline{w} \leq 0$ in $\RR^N\times(t_1,\infty)$. We consider the solution of the problem
\begin{equation}\label{eq:CAUCHYPROBLEMINITIALDATUMINCREASINGSUPERSOLFAST}
\begin{cases}
\begin{aligned}
a_1\partial_t\overline{w} - \Delta_p\overline{w} + c_0\overline{w} = 0 \quad &\text{in } \RR^N\times(t_1,\infty) \\
\overline{w}(x,t_1) = 1 + |x|^{\lambda} \;\,\quad\qquad &\text{in } \RR^N.
\end{aligned}
\end{cases}
\end{equation}
According to the resume presented in Section \ref{APPENDIXSELFSIMSOLINCINITDATA}, problem \eqref{eq:CAUCHYPROBLEMINITIALDATUMINCREASINGSUPERSOLFAST} is well posed if $0 < \lambda < p/(p-2)$ when $p > 2$. Further assumptions are not needed when $1 < p \leq 2$. Furthermore, since $c_0 > 0$ can be chosen smaller and $a_1 > 0$ larger, we make the additional assumption
\begin{equation}\label{eq:ASSUMPIONONC0ANDA1}
\frac{c_0}{a_1} = \nu\lambda.
\end{equation}
Now, we define the function
\[
\tau(t) :=
\begin{cases}
\begin{aligned}
\frac{1}{c_0(2-p)}\Big[e^{(c_0/a_1)(2-p)(t-t_1)} - 1 \Big] \;\;\quad &\text{if } 1 < p < 2 \\
\frac{1}{a_1}(t-t_1) \qquad\qquad\qquad\qquad\qquad\quad &\text{if } p = 2 \\
\frac{1}{c_0(p-2)}\Big[1 - e^{-(c_0/a_1)(p-2)(t-t_1)} \Big] \quad &\text{if } p > 2.
\end{aligned}
\end{cases}
\]
Note that $\tau = \tau(t)$ is increasing with $\tau(t_1) = 0$ for all $p > 1$. Moreover, we define the limit of $\tau(t)$ as $t \to \infty$ with the formula
\[
\tau_{\infty} :=
\begin{cases}
\begin{aligned}
\infty \qquad\qquad\quad\; &\text{if } 1 < p \leq 2 \\
[c_0(p-2)]^{-1} \quad &\text{if } p > 2.
\end{aligned}
\end{cases}
\]
Then, the function $\widetilde{w}(x,\tau) := e^{(c_0/a_1)(t-t_1)}\overline{w}(x,t)$ (with $\tau = \tau(t)$) solves the ``pure diffusive'' problem
\[
\begin{cases}
\begin{aligned}
\partial_{\tau}\widetilde{w} = \Delta_p\widetilde{w} \qquad\qquad\quad\;\; &\text{in } \RR^N\times(0,\tau_{\infty}) \\
\widetilde{w}(x,0) = 1 + |x|^{\lambda} \quad\qquad &\text{in } \RR^N.
\end{aligned}
\end{cases}
\]
As we explained in Section \ref{APPENDIXSELFSIMSOLINCINITDATA}, for all $\tau_1 \geq 0$ the problem
\begin{equation}\label{eq:DEFINITIONOFALPHABETALAMBDA}
\begin{cases}
\begin{aligned}
\partial_{\tau}U = \Delta_pU \qquad\, &\text{in } \RR^N\times(\tau_1,\infty) \\
U(x,\tau_1) =|x|^{\lambda} \quad &\text{in } \RR^N
\end{aligned}
\end{cases}
\end{equation}
admits self-similar solutions $U(x,\tau+\tau_1) = \tau^{-\alpha_{\lambda}}F(|x|(\tau+\tau_1)^{-\beta_{\lambda}})$, with self-similar exponents
\[
\alpha_{\lambda} = -\frac{\lambda}{(1-\lambda)p + 2\lambda} \quad \text{ and } \quad \beta_{\lambda} =  \frac{1}{(1-\lambda)p + 2\lambda},
\]
and profile $F(\xi) \geq 0$ with $F'(\xi) \geq 0$ for all $\xi \geq 0$, where we set $\xi = |x|(\tau+\tau_1)^{-\beta_{\lambda}}$. Note that since we assumed $0 < \lambda < p/(p-2)$ when $p > 2$, the self-similar exponents are well defined with $\alpha_{\lambda} < 0$ and $\beta_{\lambda} > 0$ for all $p > 1$. Finally, recall that it is possible to describe the spacial ``decay'' of the self-similar solutions for large values of the variable $\xi = |x|(\tau+\tau_1)^{-\beta_{\lambda}}$, with the bounds
\begin{equation}\label{eq:ESTIMATEONULARGEXI}
H_2 |x|^{\lambda} \leq U(x,\tau+\tau_1) \leq H_1 |x|^{\lambda}, \quad \text{ for all } |x|(\tau+\tau_1)^{-\beta_{\lambda}} \geq h
\end{equation}
for a constant $h \gg 0$ large enough, see formula \eqref{eq:ESTIMATEONULARGEXI1}. Now, it is not difficult to see that $\widetilde{w}(x,\tau) = 1 + U(x,\tau + \tau_1)$, for all fixed delays $\tau_1 \geq 0$. Moreover, we compute the time derivative:
\[
\begin{aligned}
\partial_t \overline{w}(x,t) & = \partial_t \big\{ e^{-\frac{c_0}{a_1}(t-t_1)}[1 + U(x,\tau + \tau_1)] \big\} \\
& = -(\tau + \tau_1)^{-\alpha_{\lambda}-1}e^{-\frac{c_0}{a_1}(t-t_1)}
\bigg\{ \frac{c_0}{a_1}(\tau + \tau_1)^{\alpha_{\lambda}+1} + \bigg[\frac{c_0}{a_1}(\tau + \tau_1) + \alpha_{\lambda} \tau' \bigg]F(\xi) + \beta_{\lambda} \tau' \xi F'(\xi) \bigg\},
\end{aligned}
\]
where $\xi = |x|(\tau + \tau_1 )^{-\beta_{\lambda}}$ and $\tau'$ stands for the derivative respect with the variable $t\geq0$. Let's set
\[
Q(t) := (c_0/a_1)(\tau + \tau_1) + \alpha_{\lambda} \tau'.
\]
Since, $F(\cdot)$, $F'(\cdot)$, and $\tau'(\cdot)$ are non-negative and $\beta_{\lambda} > 0$, in order to have $\partial_t \overline{w}(x,t) \leq 0$, it is sufficient to show $Q(t) \geq 0$ for all $t \geq t_1$ and a suitable choice of $\tau_1 > 0$.

\noindent If $p = 2$, this follows from a direct and immediate computation, choosing $\tau_1 > 0$ large enough. If $1 < p < 2$, we may proceed similarly. It is simple to see that condition $Q(t) \geq 0$ for $t \geq t_1$ reads
\[
[1 + \alpha_{\lambda}(2-p)]e^{(c_0/a_1)(2-p)(t-t_1)} \geq 1 - \frac{\tau_1}{\tau_{\infty}}.
\]
Consequently, since $1 + \alpha_{\lambda}(2-p) \geq 0$, it is sufficient to choose $\tau_1 \geq \tau_{\infty}$. Finally, when $p > 2$ it holds $\tau'(t) \leq 1/a_1$ for all $t \geq t_1$. Hence, it is simple to see that the choice $\tau_1 \geq -\alpha_{\lambda}/c_0$ is a sufficient condition so that $Q(t) \geq 0$ for all $t \geq t_1$. We stress that the choice of $\tau_1 > 0$ is independent of $t_1 > 0$.

\noindent Now, using the fact that $\partial_t\overline{w}(x,t) \leq 0$ in $\RR^N\times(t_1,\infty)$ and that $0 \leq a(x,t) \leq a_1$ in $\Omega_I$, it is straightforward to see that
\begin{equation}\label{eq:DIFFINEQSUPSOLFASTCONV1}
a(x,t)\partial_t\overline{w} - \Delta_p\overline{w} + c_0\overline{w} \geq 0 \quad \text{in } \Omega_I.
\end{equation}

\emph{Step3.} Now we compare the functions $\underline{w}$ and $\overline{w}$, applying the Maximum Principle of Section \ref{SECTIONMAXPRINCCYLDOMAINS}. Hence, we have to check that the assumptions in Proposition \ref{MAXPRINNONCYLDOMAINS} are satisfied.

(A1). It is simple to see that it holds $\underline{w}(x,t_1) \leq \overline{w}(x,t_1)$ in $\RR^N$. Indeed, we have $\overline{w}(x,t_1) \geq 1$ while $\underline{w}(x,t_1) = 1 - [u(x,t_1)]^m \leq 1$.

(A2). We have to check that $\underline{w} \leq \overline{w}$ on the boundary of $\Omega_I$, i.e., on the set $\{|x| = e^{\nu t}, t \geq t_1\}$. We use the first estimate in \eqref{eq:ESTIMATEONULARGEXI}:
\[
\begin{aligned}
\overline{w} &= e^{-(c_0/a_1)(t-t_1)}\widetilde{w}(x,\tau) =  e^{-(c_0/a_1)(t-t_1)}[1 + U(x,\tau + \tau_1)] \\
&\geq e^{-(c_0/a_1)(t-t_1)} (1 + H_2|x|^{\lambda}) = e^{-(c_0/a_1)(t-t_1)} (1 + H_2e^{\nu\lambda t})  \\
& = e^{(c_0/a_1)t_1}(H_2 + e^{-(c_0/a_1)t}) \geq 1 \geq \underline{w} \quad \text{ in } \{|x| = e^{\nu t}, t \geq t_1\}.
\end{aligned}
\]
First of all, we point out that the last equality in the preceding chain is satisfied thanks to the first assumption in \eqref{eq:ASSUMPIONONC0ANDA1}, i.e., $c_0/a_1 = \nu \lambda$.

\noindent Secondly, we note that the first inequality holds only if $|x|(\tau + \tau_1)^{-\beta_{\lambda}} \geq h$, which means
\begin{equation}\label{eq:CONDITIONONT1CONVERGENCETO1FAST}
e^{\nu t} \geq h (\tau + \tau_1)^{\beta_{\lambda}}.
\end{equation}
As the reader can easily check, when $p=2$, \eqref{eq:CONDITIONONT1CONVERGENCETO1FAST} is satisfied by taking $t_1 \geq t_0$ so that  $e^{\nu t_1} \geq h\sqrt{\tau_1}$. If $p > 2$, it is sufficient to fix $t_1 \geq t_0$ to have $e^{\nu t_1} \geq h(\tau_{\infty} + \tau_1)^{\beta_{\lambda}}$.

\noindent The case $1 < p < 2$ is a little bit subtle. First of all, set $b_{\lambda}:= 1/\beta_{\lambda} = (1-\lambda)p + 2\lambda$ and note that, thanks to assumption \eqref{eq:ASSUMPIONONC0ANDA1}, we have that \eqref{eq:CONDITIONONT1CONVERGENCETO1FAST} is automatically satisfied if
\[
e^{\nu b_{\lambda}t} \geq h^{b_{\lambda}}\big\{ [c_0(2-p)]^{-1} e^{\nu\lambda(2-p)(t-t_1)} +\tau_1 \big\},
\]
which, since $b_{\lambda} = p + \lambda(2-p)$, is equivalent to
\[
e^{\nu\lambda(2-p)t}\big\{e^{\nu p t} - [c_0(2-p)]^{-1}h^{b_{\lambda}} \big\} \geq h^{b_{\lambda}}\big\{[c_0(2-p)]^{-1}e^{-\nu \lambda(2-p)t_1} + \tau_1 \big\}.
\]
Finally, it is straightforward to see that the last inequality is satisfied for all $t \geq t_1 \geq t_0$ so that
\[
e^{\nu\lambda p t_1} \geq h^{b_{\lambda}}\big\{ 2 [c_0(2-p)]^{-1} + \tau_1 \big\}.
\]
Hence, we have that condition \eqref{eq:CONDITIONONT1CONVERGENCETO1FAST} is satisfied for all $p > 1$ when $t_1 \geq t_0$ is taken large enough.

(A3). To check this third assumption it is sufficient to combine \eqref{eq:DIFFINEQSUBSOLFASTCONV1} and \eqref{eq:DIFFINEQSUPSOLFASTCONV1}, restricting they validity to the set $\{|x| \leq e^{\nu t}, t \geq t_1\}$.

\noindent Hence, we deduce $\underline{w} \leq \overline{w}$ in $\{|x| \leq e^{\nu t}, t \geq t_1\}$ by applying Proposition \ref{MAXPRINNONCYLDOMAINS}. So we have
\[
\begin{aligned}
\underline{w}(x,t) &\leq \overline{w}(x,t) = e^{-\frac{c_0}{a_1}(t-t_1)}[1 + U(x,\tau + \tau_1)] \leq e^{-\frac{c_0}{a_1}(t-t_1)}[1 + U(e^{\nu t},\tau + \tau_1)] \\
&\leq e^{-\frac{c_0}{a_1}(t-t_1)}[1 + H_1|x|^{\lambda}] \leq e^{-\frac{c_0}{a_1}(t-t_1)}[1 + H_1e^{\nu\lambda t}]
\end{aligned}
\]
in the set $\{|x| \leq e^{\nu t}, t \geq t_1 \}$, thanks to \eqref{eq:CONDITIONONT1CONVERGENCETO1FAST}.

\noindent Now, let us fix $\varepsilon > 0$ and take a time $t_{\varepsilon}' > 0$, and a constant $H_{\varepsilon} > 0$ such that
\[
t_{\varepsilon}' \geq t_1 - \frac{a_1\ln (\varepsilon/2)}{c_0} \qquad  \text{ and }  \qquad H_{\varepsilon}^{\lambda} \leq \frac{\varepsilon}{2H_1}e^{-\frac{c_0}{a_1}t_1}.
\]
These choices combined with the previous chain of inequalities give us
\[
\underline{w}(x,t) \leq e^{-\frac{c_0}{a_1}(t_{\varepsilon}-t_1)} + H_{\varepsilon}^{\lambda}H_1e^{+\frac{c_0}{a_1}t_1} \leq \frac{\varepsilon}{2} + \frac{\varepsilon}{2} \leq \varepsilon,
\]
in the set $\{|x| \leq H_{\varepsilon}e^{\nu t}, t \geq t_{\varepsilon}'\}$. Finally, noting that $\{|x| \leq e^{\sigma t}, t \geq t_{\varepsilon}\} \subset \{|x| \leq H_{\varepsilon}e^{\nu t}, t \geq t_{\varepsilon}\}$ for all $\sigma < \nu$ and for some $t_{\varepsilon} > 0$ large enough, we complete the proof of the theorem using the arbitrariness of $\sigma < \nu < \sigma_{\ast}$. $\Box$
\paragraph{Remarks.} We end this section with two remarks. First of all, we ask the reader to note that at the beginning of the previous proof, we have made the change of variable $\underline{w} = 1 - u^m$ in order to obtain problem \eqref{eq:PROBLEMSOLVEDBYV}, which has non-constant coefficients, but the diffusion operator simplifies to a $p$-Laplacian. This is a considerable advantage since we can employ the well known theory of $p$-Laplacian diffusion and non-integrable initial data (see Section \ref{APPENDIXSELFSIMSOLINCINITDATA} for more details and references) to construct the super-solution given by problem \eqref{eq:CAUCHYPROBLEMINITIALDATUMINCREASINGSUPERSOLFAST}. A different approach could be studying the existence, uniqueness and regularity of solutions for the doubly nonlinear equation and non-integrable initial data in the fast diffusion range $0 < \widehat{\gamma} < p/N$. Up to our knowledge, this theory has not been developed yet.

\noindent Secondly, we point out that in the previous proof we showed a slightly different result too stated in the following corollary. It will be very useful in Section \ref{SECTIONBOUNDSFORLEVELSETSFAST}, where we will study the behaviour of the solution $u = u(x,t)$ on the set $|x| = e^{\sigma_{\ast}t}$.
\begin{cor}\label{CONJECTUREFORLEVELSETSSIGMAASTFAST}
Let $m > 0$ and $p > 1$ such that $0 < \widehat{\gamma} < p/N$ and let $u = u(x,t)$ be the solution of the problem \eqref{eq:REACTIONDIFFUSIONEQUATIONPLAPLACIAN} with initial datum \eqref{eq:ASSUMPTIONSONTHEINITIALDATUMFAST}. Suppose that there exist $\nu > 0$, $\overline{\varrho} > 0$, $\overline{\varepsilon} > 0$ and $t_0 >0$ such that
\[
u(x,t) \geq \overline{\varepsilon} \quad \text{in } \{|x| \leq \overline{\varrho}e^{\nu t}\} \text{ for all } \; t \geq t_0.
\]
Then, for all $0 < \omega < 1$, there exist $C_{\omega} > 0$ large enough and $t_{\omega} \geq t_0$ such that
\[
u(x,t) \geq \omega \quad \text{in } \{|x| \leq C_{\omega}^{-1}e^{\nu t}\} \text{ for all } \; t \geq t_{\omega}.
\]
\end{cor}
The proof coincides with the one of Theorem \ref{CONVERGENCETOONEFASTDIFFUSION}. Notice indeed that we have begun by assuming that $u \geq \overline{\varepsilon}:=\widetilde{\varepsilon}$ in $\{|x|\leq e^{\nu t}, t \geq t_1\}$ and for all $\varepsilon > 0$, we proved the existence of $t_{\varepsilon}'>0$ and $H_{\varepsilon}>0$ such that
\[
u(x,t) \geq 1 - \varepsilon, \quad \text{in } \{|x|\leq H_{\varepsilon}e^{\nu t}, t \geq t_{\varepsilon}'\}.
\]
We point out that in the previous statement the value of the exponent $\nu > 0$ does not change. In the proof of Theorem \ref{CONVERGENCETOONEFASTDIFFUSION}, we need to take $\sigma < \nu$ to obtain a ``convergence inner set'' not depending on $H_{\varepsilon}$. Indeed, in the second one we prove the convergence of the solution $u = u(x,t)$ to the steady state 1 in the set $\{|x| \leq e^{\sigma t}\}$ for all $\sigma < \sigma_{\ast}$, while now the exponent $\nu > 0$ is arbitrary.
%
%
%
%
%
%
%
%
%
%
%
%
\section{\texorpdfstring{\boldmath}{}Case \texorpdfstring{$f(u) = u(1-u)$}{}. Precise bounds for level sets}\label{SECTIONBOUNDSFORLEVELSETSFAST}
We devote this section to the proof of Theorem \ref{THEOREMBOUNDSFORLEVELSETSFAST}. A similar result was showed in \cite{C-R2:art} for the Fisher-KPP equation with fractional diffusion. In particular, they studied the case of the fractional Laplacian $(-\Delta)^{1/2}$  and worked in dimension $N=1$, see Theorem 1.6 of \cite{C-R2:art} for more details. In our setting, we consider the classical reaction term $f(u) = u(1-u)$ and the problem
\begin{equation}\label{eq:REACTIONDIFFUSIONEQUATIONPLAPLACIANCLASSICAL}
\begin{cases}
\partial_t u = \Delta_pu^m + u(1-u) \quad \text{in } \RR^N\times(0,\infty) \\
u(x,0) = u_0(x) \;\;\,\qquad\qquad \text{in } \RR^N
\end{cases}
\end{equation}
where $u_0(\cdot)$ satisfies \eqref{eq:ASSUMPTIONSONTHEINITIALDATUMFAST}, i.e., $u_0(x) \leq C|x|^{-p/\widehat{\gamma}}$ and $0 \leq u_0 \leq 1$, for some constant $C>0$. As always, we do not pose restrictions on the dimension $N \geq 1$, and we will work in the ``fast diffusion'' range $0 < \widehat{\gamma} < p/N$.

We divide the proof in two main parts. The first one is devoted to prove the ``upper bound'' i.e., the first inclusion \eqref{eq:STATEMENTLEVELSETSBOUNDSFAST}. In the second one, we show the ``lower bound'' (the second inclusion \eqref{eq:STATEMENTLEVELSETSBOUNDSFAST}) which is the most difficult part. We point out that in this part, we have to give to separate proofs for the ranges $\widehat{\gamma} \leq p-1$ and $\widehat{\gamma} > p-1$.
\paragraph{Proof of Theorem \ref{THEOREMBOUNDSFORLEVELSETSFAST}.} The proof is divided in some steps. We construct a super-solution and sub-solutions for problem  \eqref{eq:REACTIONDIFFUSIONEQUATIONPLAPLACIANCLASSICAL} and we use these sub- and super-solutions to prove the level sets bounds \eqref{eq:STATEMENTLEVELSETSBOUNDSFAST}. In particular, we employ the super-solution to prove the upper bound, whilst the sub-solutions for the lower bound.
\paragraph{Upper bound.} We begin to prove the first inclusion in \eqref{eq:STATEMENTLEVELSETSBOUNDSFAST}, i.e., for all $0 < \omega < 1$, there exists a constant $C_{\omega} > 0$ large enough such that it holds
\[
\{x \in \RR^N : |x| > C_{\omega}e^{\sigma_{\ast}t}\} \subset \{x \in \RR^N : u(x,t) < \omega \} \quad \text{for all } t \geq 0,
\]
where $\sigma_{\ast} = \widehat{\gamma}/p$ (recall that in this setting $f'(0) = 1$).

As we mentioned before, we need a special super-solution. Thus, we basically repeat the computations carried out at beginning of Section \ref{CONVERGENCETOZEROFAST}, by looking for radial solutions of the equation
\[
\partial_t\overline{u} = r^{1-N}\partial_r\big(r^{N-1} |\partial_r\overline{u}^m|^{p-2}\partial_r\overline{u}^m\big) + \overline{u}, \quad r = |x| > 0,
\]
with separate variables
\[
\overline{u}(r,t) = r^{-p/\widehat{\gamma}}G(t), \quad r > 0, \, t \geq 0.
\]
We have seen that for such solution, the function $G = G(t)$ has to solve the equation
\[
\frac{dG}{dt} = G + \kappa G^{1-\widehat{\gamma}}, \quad t \geq 0 \quad\quad \kappa:= \frac{(p-\widehat{\gamma} N)(mp)^{p-1}}{\widehat{\gamma}^p}.
\]
Note that since $0 <\widehat{\gamma} < p/N$ we have that $\kappa$ is well defined, positive and $1-\widehat{\gamma} = m(p-1) > 0$. As we have mentioned, we deal with a Bernoulli equation and we can integrate it:
\[
G(t) = \big(ae^{\widehat{\gamma} t} - \kappa\big)^{\frac{1}{\widehat{\gamma}}}, \quad a > \kappa, \quad t \geq 0.
\]
Hence, we found that for all $a > \kappa$, the function
\[
\overline{u}(x,t) = |x|^{-\frac{p}{\widehat{\gamma}}}\big(ae^{\widehat{\gamma} t} - \kappa\big)^{\frac{1}{\widehat{\gamma}}}
\]
is a super-solution of the equation in \eqref{eq:REACTIONDIFFUSIONEQUATIONPLAPLACIANCLASSICAL}.

Now, in order to apply the Maximum Principle, we show $\overline{u}(x,0) \geq u_0(x)$ in $\RR^N$.  We consider the function
\[
\overline{u}_0(x) =
\begin{cases}
\begin{aligned}
1 \qquad\qquad &\text{if } |x| \leq C^{\widehat{\gamma}/p} \\
C|x|^{-p/\widehat{\gamma}} \quad &\text{if } |x| \geq C^{\widehat{\gamma}/p},
\end{aligned}
\end{cases}
\]
where $C > 0$ is taken as in \eqref{eq:ASSUMPTIONSONTHEINITIALDATUMFAST}, see also the beginning of this
section. It is simple to see that we have $\overline{u}_0(x) \geq u_0(x)$ for all $x \in \RR^N$. So, taking
$a \geq k + C^{\widehat{\gamma}}$, we have
\[
\begin{aligned}
\overline{u}(x,0) = (a - \kappa)^{\frac{1}{\widehat{\gamma}}} |x|^{-\frac{p}{\widehat{\gamma}}} \geq
C |x|^{-\frac{p}{\widehat{\gamma}}} = \overline{u}_0(x) \quad &\text{ if } |x| \geq C^{\widehat{\gamma}/p}, \\
\overline{u}(x,0) = (a - \kappa)^{\frac{1}{\widehat{\gamma}}} |x|^{-\frac{p}{\widehat{\gamma}}} \geq
1 = \overline{u}_0(x) \qquad\quad\;\; &\text{ if } |x| \leq C^{\widehat{\gamma}/p}.
\end{aligned}
\]
Consequently, we have $\overline{u}(x,0) \geq u_0(x)$ in $\RR^N$ and we obtain $\overline{u}(x,t) \geq u(x,t)$ in $\RR^N\times[0,\infty)$, by applying the Maximum Principle.

\noindent We are ready to prove the first inclusion in \eqref{eq:STATEMENTLEVELSETSBOUNDSFAST}. Thus, let us fix $0 < \omega < 1$ and consider $x \in \RR^N$ satisfying $u(x,t) > \omega$. Using the super-solution $\overline{u}(x,t)$ constructed before, we have
\[
\omega < \overline{u}(x,t) = |x|^{-\frac{p}{\widehat{\gamma}}}\big(ae^{\widehat{\gamma} t} - \kappa\big)^{\frac{1}{\widehat{\gamma}}} \quad \Rightarrow \quad |x|^{\frac{p}{\widehat{\gamma}}} < \frac{1}{\omega}\big(ae^{\widehat{\gamma} t} - \kappa\big)^{\frac{1}{\widehat{\gamma}}} < \frac{a^{1/\widehat{\gamma}}}{\omega}e^t, \quad t \geq 0,
\]
which, setting $C_{\omega} = \big(a^{1/\widehat{\gamma}}/\omega\big)^{\widehat{\gamma}/p}$, implies
\[
|x| < C_{\omega}e^{\frac{\widehat{\gamma}}{p}t} = C_{\omega}e^{\sigma_{\ast}t}, \quad t \geq 0,
\]
and we conclude the proof of the first inclusion in \eqref{eq:STATEMENTLEVELSETSBOUNDSFAST}.
\paragraph{Lower bound.} We want to show the second inclusion in \eqref{eq:STATEMENTLEVELSETSBOUNDSFAST}: for all $0 < \omega < 1$, there exist $C_{\omega} > 0$ and $t_{\omega}$ large enough, such that
\[
\{x \in \RR^N : |x| < C_{\omega}^{-1}e^{\sigma_{\ast}t}\} \subset \{x \in \RR^N : u(x,t) > \omega \} \quad \text{for all } t \geq t_{\omega}.
\]
The idea consists in constructing a sub-solution $\underline{u} = \underline{u}(x,t)$ that will act as a barrier from below and then employ Corollary \ref{CONJECTUREFORLEVELSETSSIGMAASTFAST}. Even though this idea was firstly used in \cite{C-R2:art}, we stress that our construction is completely independent by the previous one, and the comparison is not done in the whole space, but only on a sub-region. Moreover, we will divide the range $0 < \widehat{\gamma} < p/N$ in two sub-ranges: $\widehat{\gamma} \leq p-1$ and $\widehat{\gamma} > p-1$. Of course, these two new ranges are always accompanied with the ``fast diffusion assumption'' $0 < \widehat{\gamma} < p/N$. We will specify when we will need to distinguish between these different cases. We divide the proof in some steps.

\emph{Step1.} In this first step, we fix some important notations and we define the candidate sub-solution. We know the following two facts:

\noindent $\bullet$ First, for all $\varepsilon > 0$ (small) and for all $r_0  > 0$ (large), there exists $t_{\varepsilon,r_0} > 0$ (large enough) such that
\[
u(x,t) \geq 1 - \varepsilon \quad \text{in } \{|x| \leq r_0 \}, \quad \text{for all } t \geq t_{\varepsilon,r_0}.
\]
This is a direct consequence of Theorem \ref{CONVERGENCETOONEFASTDIFFUSION}.

\noindent $\bullet$ Secondly, by applying Lemma \ref{LEMMAPLACINGBARENBLATTUNDERSOLUTION}, we have that for all $\theta > 0$, there exist $t_1 > \theta$, $\widetilde{\varepsilon} > 0$, and $\widetilde{\varrho}_0 > 0$ satisfying
\[
u(x,t_1)
\geq
\widetilde{u}_0(x) :=
\begin{cases}
\begin{aligned}
\widetilde{\varepsilon}  \;\,\qquad\qquad &\text{if }|x| \leq \widetilde{\varrho}_0 \\
a_0|x|^{-p/\widehat{\gamma}} \quad &\text{if }|x| > \widetilde{\varrho}_0
\end{aligned}
\end{cases}
\quad \text{in } \RR^N,
\]
where $a_0 := \widetilde{\varepsilon}\,\widetilde{\varrho}_0^{\,p/\widehat{\gamma}}$, see \eqref{eq:INITIALDATUMTESTFASTLEVELSETS}.

\noindent Now, for all $m > 0$ and $p > 1$ satisfying $0 < \widehat{\gamma} < p/N$, and $N \geq 1$, we consider the constants
\begin{equation}\label{eq:DEFINITIONDCOSTANTFASTDIFFCASE}
d_1 := p/\widehat{\gamma}^2(p-\widehat{\gamma}N), \quad d_2 := p/\widehat{\gamma}^2[(p-1)(p-\widehat{\gamma}) + \widehat{\gamma}(N-1)], \quad d_3 := (p/\widehat{\gamma})^{2-p}m^{1-p}.
\end{equation}
Note that the assumption $0 < \widehat{\gamma} < p/N$ guarantees that $d_i$ are positive for all $i =1,2,3$. The importance of these constants will be clear later. Now, they are simply needed to choose $r_0 > 0$. We take $r_0$ large (depending only on $m$, $p$, $N$ and $\varepsilon$), satisfying
\begin{equation}\label{eq:CHOICEOFR0ANDT0FASTDIFFCASE}
r_0^p \geq \frac{d_2^{\widehat{\gamma}+1}(1 - \varepsilon)^{-\widehat{\gamma}}}{d_1^{\widehat{\gamma}}d_3(\widehat{\gamma}+1)}\varepsilon^{-1}
\quad
\text{ and }
\quad
r_0^p \geq \frac{d_2^2(d_1 + d_1/d_2)^{\widehat{\gamma}-(p-1)}}{pd_1d_3}\varepsilon^{-1/\widehat{\gamma}}. \end{equation}
The first assumption will be needed in the range $\widehat{\gamma} \leq p-1$, while the second when $\widehat{\gamma} > p-1$. Then we take $\theta := t_{\varepsilon,r_0}$ and $t_0 := t_1 > t_{\varepsilon,r_0}$ and we fix $\widetilde{\varepsilon} > 0$ and $\widetilde{\varrho}_0 > 0$ corresponding to the value of $t_0 > 0$.

Let us define the candidate sub-solution. We consider the function
\[
\underline{u}(r,t) = \frac{e^t}{b\psi(r) + ce^t}, \quad r \geq r_0, \, t \geq t_0, \; \text{ and } \; \psi(r) = r^{p/\widehat{\gamma}},
\]
where, of course, $r = |x|$. We fix
\begin{equation}\label{eq:CHOICEOFPARAMETERSAC}
c = (1-\varepsilon)^{-1} > 1.
\end{equation}
The parameter $b > 0$ will be chosen in the next step, independently from $c$. We ask the reader to note that $c$ depends only on $\varepsilon > 0$. The fact that $c > 1$ will be important later.

\emph{Step2.} Now, we consider the region $\mathcal{R}_0 := \{|x| \geq r_0\} \times [t_0,\infty)$ and, we show that
\[
\underline{u}(x,t_0) \leq u(x,t_0) \quad \text{ for } |x| \geq r_0 \qquad \text{and} \qquad \underline{u}(|x|=r_0,t) \leq u(|x|=r_0,t) \quad \text{ for } t \geq t_0,
\]
in order to assure that $\underline{u}(x,t)$ and $u(x,t)$ are well-ordered at time $t = t_0$ and on the boundary of $\mathcal{R}_0$.

\noindent $\bullet$ Comparison in $\{|x| \geq r_0\}$ at time $t = t_0$. Since both the parameters $b$ and $c$ are positive, we have
\[
\underline{u}(x,t_0) \leq
\begin{cases}
\begin{aligned}
b^{-1}e^{t_0}|x|^{-p/\widehat{\gamma}} \quad &\text{for all } x \in \RR^N \\
b^{-1}e^{t_0}r_0^{-p/\widehat{\gamma}} \;\;\quad &\text{for all } |x| \geq r_0,
\end{aligned}
\end{cases}
\]
and so, comparing with $\widetilde{u}_0 = \widetilde{u}_0(x)$, we obtain
\[
u(x,t_0) \geq \widetilde{u}_0(x) \geq \underline{u}(x,t_0) \quad \text{ for all } |x| \geq r_0,
\]
by taking the parameter $b > 0$ large enough depending on $r_0 > 0$ (but not on $c$):
\[
b \geq (e^{t_0}/\widetilde{\varepsilon})\max\{\widetilde{\varrho}_0^{\,-p/\widehat{\gamma}}, r_0^{-p/\widehat{\gamma}}\}.
\]
$\bullet$ Comparison on the boundary $\{|x|=r_0\}\times\{t \geq t_0\}$. This part is simpler. Indeed, we have
\[
\underline{u}(|x|=r_0,t) \leq 1/c = 1 - \varepsilon \leq u(|x|=r_0,t), \quad \text{for all } t \geq t_0,
\]
thanks to our assumptions on $t_0 > 0$ and $c > 0$, see \eqref{eq:CHOICEOFPARAMETERSAC}.

\emph{Step3.} In this step we prove that $\underline{u} = \underline{u}(r,t)$ is a sub-solution of the equation in \eqref{eq:REACTIONDIFFUSIONEQUATIONPLAPLACIANCLASSICAL} in the region $\mathcal{R}_0$. For the reader convenience, we introduce the expression
\[
A(r,t) := b\psi(r) + ce^t \quad \Rightarrow \quad \underline{u}(r,t) = e^t A(r,t)^{-1}.
\]
We proceed by carrying out some computations. We have
\begin{equation}\label{eq:REACTIONANDTIMEDERIVSUBSOLUTIONBOUNDLEVELSETS}
\partial_t \underline{u} = be^t\psi(r)A(r,t)^{-2}, \qquad
-\underline{u}(1-\underline{u}) = -e^t\big[b\psi(r) + (c-1)e^t\big]A(r,t)^{-2}.
\end{equation}
Now, we need to compute the radial $p$-Laplacian of $\underline{u}^m$, which is given by the formula
\[
-\Delta_{p,r}\underline{u}^m := -r^{1-N}\partial_r\big(r^{N-1} |\partial_r\overline{u}^m|^{p-2}\partial_r\overline{u}^m\big).
\]
First of all, setting $B(t) := (mb)^{p-1}e^{(1-\widehat{\gamma})t} > 0$ and using the fact that $(m+1)(p-1) = p-\widehat{\gamma}$, it is not difficult to obtain
\[
|\partial_r\underline{u}^m|^{p-2}\partial_r\underline{u}^m = -B(t)|\psi'(r)|^{p-2}\psi'(r)A(r,t)^{\widehat{\gamma}-p},
\]
where $\psi' = d\psi/dr$. Consequently, we have
\[
\begin{aligned}
-\Delta_{p,r}\underline{u}^m &= B(t)r^{1-N}\partial_r\big[ r^{N-1} |\psi'(r)|^{p-2}\psi'(r) A(r,t)^{\widehat{\gamma}-p}\big] \\
&= B(t)|\psi'(r)|^{p-2}A(r,t)^{\widehat{\gamma}-p-1}\bigg\{\bigg[\frac{N-1}{r}\psi'(r) + (p-1)\psi''(r)\bigg] A(r,t) - b(p-\widehat{\gamma})(\psi'(r))^2 \bigg\}.
\end{aligned}
\]
Combining the last quantity with the ones in \eqref{eq:REACTIONANDTIMEDERIVSUBSOLUTIONBOUNDLEVELSETS} and multiplying by $B(t)^{-1}|\psi'(r)|^{2-p}A(r,t)^{1+p-\widehat{\gamma}}$, we obtain
\[
\begin{aligned}
B(t)^{-1}&|\psi'(r)|^{2-p}A(r,t)^{1+p-\widehat{\gamma}}\big[\partial_t\underline{u} -\Delta_{p,r}\underline{u}^m - \underline{u}(1-\underline{u})\big] \\
& = \frac{c-1}{(mb)^{p-1}}e^{(1+\widehat{\gamma})t}|\psi'(r)|^{2-p}A(r,t)^{p-1-\widehat{\gamma}}
+ \Big[(p-1)\psi''(r) + \frac{N-1}{r}\psi'(r) \Big] A(r,t) \, + \\
& - b(p - \widehat{\gamma})(\psi'(r))^2.
\end{aligned}
\]
Let us take $\psi(r) = r^{p/\widehat{\gamma}}$ with $\psi'(r) = (p/\widehat{\gamma})r^{\frac{p}{\widehat{\gamma}}-1}$ and $\psi''(r) = (p/\widehat{\gamma})\big(p/\widehat{\gamma}-1\big)r^{\frac{p}{\widehat{\gamma}}-2}$. Since
\[
\begin{aligned}
(p-1)\psi''(r) + &\frac{N-1}{r}\psi'(r) = d_2 \, r^{p/\widehat{\gamma} - 2}, \\
(\psi'(r))^2 &= (p/\widehat{\gamma})^2r^{2(p/\widehat{\gamma} - 1)}, \qquad
|\psi'(r)|^{2-p} = (p/\widehat{\gamma})^{2-p}r^{(2-p)(p/\widehat{\gamma} - 1)}
\end{aligned}
\]
and recalling that $A(r,t) = b r^{p/\widehat{\gamma}} + c e^t$, we substitute in the previous equation
deducing
\[
\begin{aligned}
B(t)^{-1}&|\psi'(r)|^{2-p}A(r,t)^{1+p-\widehat{\gamma}}\big[\partial_t\underline{u} -\Delta_{p,r}\underline{u}^m - \underline{u}(1-\underline{u})\big] \\
& = -d_3\frac{c-1}{b^{p-1}}e^{(1+\widehat{\gamma})t}r^{(2-p)(p/\widehat{\gamma} - 1)} \big(br^{p/\widehat{\gamma}} + ce^t\big)^{p-1-\widehat{\gamma}}  + d_2ce^tr^{p/\widehat{\gamma} - 2} - bd_1r^{2(p/\widehat{\gamma}-1)},
\end{aligned}
\]
where $d_i > 0$, $i =1,2,3$ are chosen as in \eqref{eq:DEFINITIONDCOSTANTFASTDIFFCASE}:
\[
d_1 := p/\widehat{\gamma}^2(p-\widehat{\gamma}N), \quad
d_2 := p/\widehat{\gamma}^2[(p-1)(p-\widehat{\gamma}) + \widehat{\gamma}(N-1)], \quad
d_3 := (p/\widehat{\gamma})^{2-p}m^{1-p}.
\]
Now, multiplying by $r^{-2(p/\widehat{\gamma}-1)}$ and setting $\xi = e^tr^{-p/\widehat{\gamma}} > 0$, it is not difficult to obtain
\[
\begin{aligned}
B(t)^{-1}&|\psi'(r)|^{2-p}r^{-2(p/\widehat{\gamma}-1)}A(r,t)^{1+p-\widehat{\gamma}}\big[\partial_t\underline{u} -\Delta_{p,r}\underline{u}^m - \underline{u}(1-\underline{u})\big] \\
& = -d_3\frac{c-1}{b^{p-1}} \,r^p\, \xi^{1+\widehat{\gamma}} \big(b + c\xi\big)^{p-1-\widehat{\gamma}}  + d_2c\xi - bd_1 \\
& \leq -d_3\frac{c-1}{b^{p-1}} \,r_0^p\, \xi^{1+\widehat{\gamma}} \big(b + c\xi\big)^{p-1-\widehat{\gamma}}  + d_2c\xi - bd_1 := -C_{r_0}(\xi),
\end{aligned}
\]
for all $r \geq r_0$.
\paragraph{Case $\boldsymbol{\widehat{\gamma} \leq p-1}$.} To prove that $\underline{u} = \underline{u}(r,t)$ is a sub-solution, it is sufficient to check that
\begin{equation}\label{eq:SUBSOLUTIONCONDITIONCMINUSFAST}
C_{r_0}(\xi) = d_3\frac{c-1}{b^{p-1}} \,r_0^p\, \xi^{1+\widehat{\gamma}} \big(b + c\xi\big)^{p-1-\widehat{\gamma}} - d_2c\xi + bd_1 \geq 0,
\end{equation}
for all $\xi > 0$. We will prove the previous inequality in two separate intervals $0 \leq \xi \leq \xi_0$ and $\xi \geq \xi_0$, where $\xi_0 > 0$ will be suitably chosen.

\noindent Suppose $0 \leq \xi \leq \xi_0$. In this interval we have $C_{r_0} (\xi) \geq - d_2c\xi_0 + bd_1$ and so, a sufficient condition so that \eqref{eq:SUBSOLUTIONCONDITIONCPLUSFAST} is satisfied (for $0 \leq \xi \leq \xi_0$) is
\begin{equation}\label{eq:SUFFICIENTCONDITIONSMALLINTERVALCMINUS}
c \leq (d_1/d_2)\, b \, \xi_0^{-1}.
\end{equation}

\noindent Suppose $\xi \geq \xi_0$ and assume \eqref{eq:SUFFICIENTCONDITIONSMALLINTERVALCMINUS} to be true. Since we are in the range $\widehat{\gamma} \leq p-1$, we have $(b + c\xi)^{p-1-\widehat{\gamma}} \geq b^{p-1-\widehat{\gamma}}$, and so
\[
C_{r_0} (\xi) \geq C_{1,r_0} (\xi) := d_3\frac{c-1}{b^{\widehat{\gamma}}} \,r_0^p\, \xi^{1+\widehat{\gamma}} - d_2c\xi + bd_1.
\]
Now, we note that condition \eqref{eq:SUFFICIENTCONDITIONSMALLINTERVALCMINUS} not only implies $C_{r_0}(\xi) \geq 0$, but also  $C_{1,r_0}(\xi) \geq 0$ for all $0 \leq \xi \leq \xi_0$. Hence, in order to prove that $C_{1,r_0}(\xi) \geq 0$ for all $\xi \geq \xi_0$, it is sufficient to show that the minimum point of $C_{1,r_0}(\cdot)$ is attained for some $0 < \xi_m \leq \xi_0$. It is straightforward to compute the minimum point $\xi_m$ of $C_{1,r_0}(\cdot)$:
\[
\xi_m^{\widehat{\gamma}} = \frac{d_2}{d_3 (1 + \widehat{\gamma})r_0^p} \frac{c \, b^{\widehat{\gamma}}}{(c-1)}.
\]
For our purpose we may choose
\[
\xi_0^{\widehat{\gamma}} = \xi_m^{\widehat{\gamma}} = \frac{d_2}{d_3 (1 + \widehat{\gamma})r_0^p} \frac{c \, b^{\widehat{\gamma}}}{(c-1)}.
\]
Now, since $\xi_0$ depends on $c$, we need to check that our choice of $\xi_0 > 0$ is compatible with \eqref{eq:SUFFICIENTCONDITIONSMALLINTERVALCMINUS}, which we have assumed to be true. Thus, substituting the value of $\xi_0$ in \eqref{eq:SUFFICIENTCONDITIONSMALLINTERVALCMINUS}, we obtain that the parameter $c$ has to satisfy the inequality
\[
\frac{d_2^{\widehat{\gamma}+1}}{d_1^{\widehat{\gamma}} d_3(\widehat{\gamma} + 1) r_0^p} \, c^{1+\widehat{\gamma}} \leq c - 1.
\]
The crucial fact is that the previous expressions do not depend on $b$. Indeed, taking $c = (1-\varepsilon)^{-1}$ as in \eqref{eq:CHOICEOFPARAMETERSAC}, we can rewrite the previous inequality as
\[
r_0^p \geq \frac{d_2^{\widehat{\gamma}+1}(1 - \varepsilon)^{-\widehat{\gamma}}}{d_1^{\widehat{\gamma}}d_3(\widehat{\gamma}+1)}\varepsilon^{-1},
\]
which exactly our first assumption in \eqref{eq:CHOICEOFR0ANDT0FASTDIFFCASE} on $r_0 > 0$. This proves that for all $\widehat{\gamma} \leq p-1$ and $0 < \widehat{\gamma} < p/N$, the function $\underline{u} = \underline{u}(r,t)$ is a sub-solution for the equation in \eqref{eq:REACTIONDIFFUSIONEQUATIONPLAPLACIANCLASSICAL} in the region $\mathcal{R}_0 = \{|x| \geq r_0\} \times [t_0,\infty)$.
\paragraph{Case $\boldsymbol{\widehat{\gamma} > p-1}$.} In this range the proof is similar, but there are some technical changes that have to be highlighted. We rewrite $C_{r_0}(\cdot)$ as
\[
\begin{aligned}
-C_{r_0}(\xi) :&= -d_3\frac{c-1}{b^{p-1}} \,r_0^p\, \xi^{1+\widehat{\gamma}} \big(b + c\xi\big)^{p-1-\widehat{\gamma}}  + d_2c\xi - bd_1  \\
&= -d_3\frac{c-1}{b^{p-1}} \,r_0^p\, \bigg(\frac{\xi}{b + c\xi}\bigg)^{\widehat{\gamma}-(p-1)} \xi^p  + d_2c\xi - bd_1.
\end{aligned}
\]
So, in order to show that $\underline{u} = \underline{u}(r,t)$ is a sub-solution, we can verify that
\begin{equation}\label{eq:SUBSOLUTIONCONDITIONCPLUSFAST}
C_{r_0}(\xi) = d_3\frac{c-1}{b^{p-1}} \,r_0^p\, \bigg(\frac{\xi}{b + c\xi} \bigg)^{\widehat{\gamma}-(p-1)} \xi^p  - d_2c\xi + bd_1\geq 0,
\end{equation}
for all $\xi > 0$. Again we will pick a ``good'' $\xi_0 > 0$ and prove \eqref{eq:SUBSOLUTIONCONDITIONCPLUSFAST} in the intervals $0 \leq \xi \leq \xi_0$ and $\xi \geq \xi_0$.

\noindent Suppose $0 \leq \xi \leq \xi_0$. As before, in this interval we have $C_{r_0} (\xi) \geq - d_2c\xi_0 + bd_1$ and so, taking again $c$ as in \eqref{eq:SUFFICIENTCONDITIONSMALLINTERVALCMINUS}, i.e.
\[
c \leq (d_1/d_2)\, b \, \xi_0^{-1},
\]
then \eqref{eq:SUBSOLUTIONCONDITIONCPLUSFAST} is automatically satisfied (for $0 \leq \xi \leq \xi_0$).

\noindent Now, suppose $\xi \geq \xi_0$ and assume again \eqref{eq:SUFFICIENTCONDITIONSMALLINTERVALCMINUS} to be true. Since we are in the range $\widehat{\gamma} > p-1$, the function
\[
\xi \to \bigg(\frac{\xi}{b + c\xi} \bigg)^{\widehat{\gamma}-(p-1)}
\]
is increasing (in $\xi$), we have
\[
\begin{aligned}
C_{r_0} (\xi) & \geq d_3\frac{c-1}{b^{p-1}} \,r_0^p\, \bigg(\frac{\xi_0}{b + c\xi_0} \bigg)^{\widehat{\gamma}-(p-1)} \xi^p  - d_2c\xi + bd_1 \\
& \geq d_3\frac{c-1}{b^{p-1}} \,r_0^p\, \bigg(\frac{\xi_0}{b + (d_1/d_2)b} \bigg)^{\widehat{\gamma}-(p-1)} \xi^p  - d_2c\xi + bd_1 \\
& = \frac{d_3(c-1)}{(1 + d_1/d_2)^{\widehat{\gamma}-(p-1)}b^{p-1}} \xi_0^{\widehat{\gamma}-(p-1)} \,r_0^p\, \xi^p  - d_2c\xi + bd_1 := \widetilde{C}_{1,r_0}(\xi),
\end{aligned}
\]
where we used \eqref{eq:SUFFICIENTCONDITIONSMALLINTERVALCMINUS} in the second inequality. Exactly as in the previous case, condition \eqref{eq:SUFFICIENTCONDITIONSMALLINTERVALCMINUS} implies both $C_{r_0}(\xi) \geq 0$ and  $\widetilde{C}_{1,r_0}(\xi) \geq 0$ for all $0 \leq \xi \leq \xi_0$. Hence, we show that the minimum point of $\widetilde{C}_{1,r_0}(\cdot)$ is attained for $\xi_m = \xi_0$ and this gives us $\widetilde{C}_{1,r_0}(\xi) \geq 0$ for all $\xi \geq \xi_0$. The minimum point $\xi_m$ of $\widetilde{C}_{1,r_0}(\cdot)$ is given by the formula:
\[
\xi_m^{p-1} = \frac{d_2(1 + d_1/d_2)^{\widehat{\gamma}-(p-1)}}{pd_3 \, \xi_0^{\widehat{\gamma}-(p-1)} r_0^p} \frac{c \, b^{\widehat{\gamma}}}{c-1}.
\]
So we ask $\xi_m = \xi_0$, i.e.:
\[
\xi_0^{\widehat{\gamma}} = \frac{d_2(1 + d_1/d_2)^{\widehat{\gamma}-(p-1)}}{pd_3 \, r_0^p} \frac{c \, b^{\widehat{\gamma}}}{c-1}.
\]
Again we must check the compatibility between our choice of $\xi_0 > 0$ and \eqref{eq:SUFFICIENTCONDITIONSMALLINTERVALCMINUS}. So, we substitute the value of $\xi_0$ in \eqref{eq:SUFFICIENTCONDITIONSMALLINTERVALCMINUS} and we obtain the inequality
\[
\frac{d_2^2(1 + d_1/d_2)^{\widehat{\gamma}-(p-1)}}{pd_1d_3 \, r_0^p}c^{\frac{1+\widehat{\gamma}}{\widehat{\gamma}}} \leq  (c - 1)^{1/\widehat{\gamma}},
\]
in the parameter $c$. Also in this this case is is really important that the previous expressions do not depend on $b$. We take $c = (1-\varepsilon)^{-1}$ as in \eqref{eq:CHOICEOFPARAMETERSAC} and we rewrite the previous inequality as
\[
r_0^p \geq \frac{d_2^2(1 + d_1/d_2)^{\widehat{\gamma}-(p-1)}}{pd_1d_3} \varepsilon^{-1/\widehat{\gamma}} (1 - \varepsilon)^{(1 + 2\widehat{\gamma})/[\widehat{\gamma}(1 + \widehat{\gamma})]}.
\]
Since $1 - \varepsilon \leq 1$ the last inequality is satisfied thanks to the assumption on $r_0 > 0$ in \eqref{eq:CHOICEOFR0ANDT0FASTDIFFCASE}. Hence, we have showed that $\underline{u} = \underline{u}(r,t)$ is a sub-solution for the equation in \eqref{eq:REACTIONDIFFUSIONEQUATIONPLAPLACIANCLASSICAL} in the region $\mathcal{R}_0 = \{|x| \geq r_0\} \times [t_0,\infty)$, for the range $\widehat{\gamma} > p-1$, too.

Consequently, for all $0 < \widehat{\gamma} < p/N$, we obtain
\[
u(x,t) \geq \underline{u}(x,t) \quad \text{in } \{|x| \geq r_0\} \times [t_0,\infty),
\]
thanks to the comparison at time $t = t_0$ and on the boundary of $\mathcal{R}_0$ done in \emph{Step2}. Note that the Maximum Principle can be applied since $0 \leq \underline{u}(r,t) \leq 1/c = 1-\varepsilon$ in $\RR^N\times[0,\infty)$ and $f(u) = u(1-u)$ can be re-defined outside $[0,1-\varepsilon]$ to be Lipschitz continuous.

\emph{Step4.} In this last step, we conclude the proof. The following procedure holds for all $0 < \widehat{\gamma} < p/N$ (see also \cite{C-R2:art}). Thanks to Corollary \ref{CONJECTUREFORLEVELSETSSIGMAASTFAST}, to deduce the second inclusion in \eqref{eq:STATEMENTLEVELSETSBOUNDSFAST}:
\[
\forall \, 0 < \omega < 1, \;\; \exists \, t_{\omega},\, C_{\omega} \gg 0: \;\; \{|x| < C_{\omega}^{-1}e^{\sigma_{\ast}t}\} \subset \{u(x,t) > \omega \}, \;\; \forall t \geq t_{\omega},
\]
it is sufficient to prove $u(x,t) \geq \overline{\varepsilon}$ in $\{|x| \leq e^{\sigma_{\ast}t} = e^{\widehat{\gamma}t/p}\}\times[t_0,\infty)$, for some $\overline{\varepsilon} > 0$. So, in the set $\{r_0 \leq |x| \leq e^{\sigma_{\ast}t}\}\times[t_0,\infty)$, we have
\[
u(x,t) \geq \underline{u}(x,t) = \frac{e^t}{b|x|^{p/\widehat{\gamma}} + ce^t} \geq \frac{1}{b+c} := \overline{\varepsilon}.
\]
Note that the bound $u(x,t) \geq \overline{\varepsilon}$ can be extended to the region $\{|x| \leq r_0\}\times[t_0,\infty)$, thanks to our assumption on $t_0$ and Theorem \ref{CONVERGENCETOONEFASTDIFFUSION}. Consequently, applying Corollary \ref{CONJECTUREFORLEVELSETSSIGMAASTFAST} with $\nu = \sigma_{\ast}$, $\overline{\varepsilon} = 1/(b+c)$ and $\overline{\varrho} = 1$, we end the proof of the theorem. $\Box$
%
%
%
%
%
%
%
%
%
%
%
\section{A Maximum Principle in non-cylindrical domains}\label{SECTIONMAXPRINCCYLDOMAINS}
In this brief section, we give the proof of a Maximum Principle for a certain class of parabolic equations with $p$-Laplacian diffusion. As mentioned in the introduction, a similar result have been introduced in \cite{C-R2:art}, but proved with different techniques. This comparison principle is crucial in the study of the asymptotic behaviour of the general solutions of the Fisher-KPP problem, see Theorem \ref{CONVERGENCETOONEFASTDIFFUSION}.

\noindent Before proceeding we need to introduce some definitions. First of all, let $r \in C^1([0,\infty);\RR)$ be a positive and non-decreasing function, and consider the ``inner-sets''
\[
\Omega_I^T := \{(x,t) \in \RR^N\times[0,T) : |x| \leq r(t) \}, \quad 0 < T \leq \infty, \qquad \text{with }\; \Omega_I^{\infty}:=\Omega_I.
\]
Now, for all $p > 1$, we consider the equation
\begin{equation}\label{eq:PROBLEMPLAPFORMAXPRIN}
a(x,t)\partial_tu - \Delta_pu + c_0u = 0 \quad \text{in } \RR^N\times(0,\infty) \\
\end{equation}
where $a = a(x,t)$ is a continuous function in $\RR^N\times(0,\infty)$, with $0 < a_0 \leq a(x,t) \leq a_1 < \infty$ in $\Omega_I$, $c_0 > 0$ and $u_0 \in L^1(\RR^N)$. The next definition is given following \cite{Bon-Vaz:art, DB:book}. See also \cite{V2:book}, Chapter 8 for the Porous Medium setting.

\begin{defn}
A nonnegative function $\overline{u} = \overline{u}(x,t)$ is said to be a ``local strong'' super-solution of equation \eqref{eq:PROBLEMPLAPFORMAXPRIN} in $\Omega_I^T$ if

\noindent (i) $\overline{u} \in C_{loc}(0,T:L_{loc}^2(\RR^N)) \cap L_{loc}^p(0,T:W_{loc}^{1,p}(\RR^N))$, and $\partial_t\overline{u} \in L_{loc}^2(\RR^N\times(0,\infty))$;

\noindent (ii) $\overline{u} = \overline{u}(x,t)$ satisfies
\[
\int_{\Omega_I^T}[a(x,t)\partial_t\overline{u} + c_0\overline{u}]\eta + |\nabla \overline{u}|^{p-2}\nabla\overline{u}\,\nabla\eta \geq 0,
\]
for all test function $\eta \in C_c^1(\Omega_I^T)$, $\eta \geq 0$.

\medskip

\noindent A nonnegative function $\underline{u} = \underline{u}(x,t)$ is said to be a ``local strong'' sub-solution of equation \eqref{eq:PROBLEMPLAPFORMAXPRIN} in $\Omega_I^T$ if

\noindent (i) $\underline{u} \in C_{loc}(0,T:L_{loc}^2(\RR^N)) \cap L_{loc}^p(0,T:W_{loc}^{1,p}(\RR^N))$, and $\partial_t\underline{u} \in L_{loc}^2(\RR^N\times(0,\infty))$;

\noindent (ii) $\underline{u} = \underline{u}(x,t)$ satisfies
\[
\int_{\Omega_I^T}[a(x,t)\partial_t\underline{u} + c_0\underline{u}]\eta + |\nabla \underline{u}|^{p-2}\nabla\underline{u}\,\nabla\eta \leq 0,
\]
for all test function $\eta \in C_c^1(\Omega_I^T)$, $\eta \geq 0$.
\end{defn}
\begin{prop}\label{MAXPRINNONCYLDOMAINS}
Consider two functions $\overline{u} = \overline{u}(x,t)$ and $\underline{u} = \underline{u}(x,t)$ defined and continuous in $\RR^N\times(0,\infty)$. Assume that:

\noindent (A1) \; $\overline{u}(x,0) \geq \underline{u}(x,0)$ in $\RR^N$.

\noindent (A2) \; $\overline{u}(x,t) \geq \underline{u}(x,t)$ in $\partial\Omega_I = \{(x,t) \in \RR^N\times(0,\infty) : |x| = r(t)\}$.

\noindent (A3) Finally, assume that $\overline{u} = \overline{u}(x,t)$ is a ``local strong'' super-solution and $\underline{u} = \underline{u}(x,t)$ is a ``local strong'' sub-solution of equation  \eqref{eq:PROBLEMPLAPFORMAXPRIN} in $\Omega_I^T$.

\medskip

\noindent Then $\overline{u} \geq \underline{u}$ in $\Omega_I^T$.
\end{prop}
\emph{Proof.} Let's fix $0 < T \leq \infty$. For all $0 < t < T$, we define the subset of $\RR^N$
\[
\Omega_{I,t} := \{x \in \RR^N: |x|\leq r(t)\}.
\]
We show that for all $t > 0$, it holds
\begin{equation}\label{eq:AUXILIARYTHESISMAXPRINCFAST}
\big\|[\underline{u}(t) - \overline{u}(t)]_+\big\|_{L^1(\Omega_{I,t})} \leq \big\|[\underline{u}(0) - \overline{u}(0)]_+\big\|_{L^1(\RR^N)},
\end{equation}
where $[\cdot]_+$ stands for the positive part. Consequently, we deduce the thesis thanks to assumption (A1). We proceed with a standard argument, see for instance Chapter 8 of \cite{V2:book} for the Porous Medium equation.

\noindent Let's consider a function $p \in C^1(\RR)$ such that
\[
0 \leq p \leq 1, \qquad p(s) = 0 \; \text{ for } s \leq 0, \qquad p'(s) > 0 \; \text{ for } s > 0,
\]
and a sequence $w_j \in C^1(\Omega_I^T)$ such that $w_j \to \underline{u} - \overline{u}$ as $j \to \infty$ in $L_{loc}^p(0,T:W_{loc}^{1,p}(\RR^N))$. Note that we can suppose
\[
w_j \leq 0 \; \text{ on }\; \partial\Omega_I^T = \{(x,t) \in \RR^N\times(0,T): |x| = r(t)\}
\]
thanks to assumption (A2). Hence, if $h \in C_0^1([0,T])$ with $0 \leq h \leq 1$, we can take as test function
\[
\eta_j = p(w_j)h(t), \quad j = 1,2,\ldots
\]
Thus, by the definition of sub- and super-solutions, it is simple to deduce
\[
\int_{\Omega_I^T}\big[a(x,t)\partial_t(\underline{u}-\overline{u}) + c_0(\underline{u} - \overline{u})\big]p(w_j)h + \big<\,|\nabla\underline{u}|^{p-2}\nabla\underline{u} - |\nabla\overline{u}|^{p-2}\nabla\overline{u},\nabla w_j\,\big> p'(w_j)h \,dxdt \leq 0.
\]
The second integral converges to
\[
\int_{\Omega_I^T} \big<\,|\nabla\underline{u}|^{p-2}\nabla\underline{u} - |\nabla\overline{u}|^{p-2}\nabla\overline{u},\nabla \underline{u} - \nabla\overline{u}\,\big> p'(\underline{u} -\overline{u})h \,dxdt \geq 0,
\]
thanks to the fact that $\big<|b|^{p-2}b - |a|^{p-2}a,b-a\big> \geq 0$ for all $a,b \in \RR^N$ and $p > 1$, see the last section of \cite{Lindq:art}. Hence, taking the limit in the second integral we deduce
\[
\int_{\Omega_I^T}\big[a(x,t)\partial_t(\underline{u}-\overline{u}) + c_0(\underline{u} - \overline{u})\big]p(\underline{u}-\overline{u})h \,dxdt\leq 0,
\]
and, letting $p(\cdot) \to sign_+(\cdot):=[sign]_+(\cdot)$, we obtain
\[
\int_{\Omega_I^T}a(x,t)\partial_t(\underline{u}-\overline{u})sign_+(\underline{u}-\overline{u})h + c_0(\underline{u}-\overline{u})sign_+(\underline{u}-\overline{u})h \,dxdt \leq 0.
\]
Now, we have
\[
\frac{d}{dt} [\,\underline{u} - \overline{u}\,]_+ = \partial_t(\underline{u}-\overline{u})sign_+(\underline{u}-\overline{u}),
\]
and, since $[s]_+ = s\cdot sign_+(s) \geq 0$, $a(x,t) \geq a_0 > 0$, and $c_0 > 0$ we easily get
\[
\int_0^T \bigg(\int_{\Omega_{I,t}} \partial_t [\,\underline{u}(t) - \overline{u}(t)\,]_+ dx\bigg) h(t) \,dt \leq 0 \quad \text{ for all }\; h \in C_c^1([0,T]), \; 0 \leq h \leq 1.
\]
Thus, thanks to arbitrariness of $h$, we deduce that
\[
\int_{\Omega_{I,t}} \partial_t [\,\underline{u}(t) - \overline{u}(t)\,]_+ dx \leq 0,
\]
for all $t > 0$. Using assumption (A2) again, it is not difficult to deduce
\[
\frac{d}{dt}\bigg(\int_{\Omega_{I,t}} [\,\underline{u}(t) - \overline{u}(t)\,]_+ dx\bigg) \leq 0,
\]
which implies
\[
\big\|[\underline{u}(t) - \overline{u}(t)]_+\big\|_{L^1(\Omega_{I,t})} \leq \big\|[\underline{u}(0) - \overline{u}(0)]_+\big\|_{L^1(\Omega_{I,0})} \leq \big\|[\underline{u}(0) - \overline{u}(0)]_+\big\|_{L^1(\RR^N)},
\]
i.e., the thesis. $\Box$
\paragraph{Remark.} We point out that the functions we use in the proof of Theorem
\ref{CONVERGENCETOONEFASTDIFFUSION} satisfy the assumptions of regularity required in the statement of
Proposition \ref{MAXPRINNONCYLDOMAINS}, as we have remarked in the introduction. See also the bibliography
reported in the next section.
%
%
%
%
%
%
%
%
%
%
%

\section{Appendix: Self-similar solutions for increasing initial data}\label{APPENDIXSELFSIMSOLINCINITDATA}
In this section, we recall some basic facts about the existence of Barenblatt solutions for the Cauchy problem
\begin{equation}\label{eq:CAUCHYPROBLEMPLAPLACIANINCREASINGINITIALDATA}
\begin{cases}
\begin{aligned}
\partial_tu = \Delta_p u \;\;\;\qquad &\text{in } \RR^N\times(0,\infty) \\
u(x,0) = u_0(x) \quad &\text{in } \RR^N,
\end{aligned}
\end{cases}
\end{equation}
where $p > 1$. In particular, we focus on the specific initial datum
\begin{equation}\label{eq:ASSUMPTIONINDATUMPOWERTYPE}
u_0(x) = |x|^{\lambda}, \quad \lambda > 0.
\end{equation}
A more complete analysis of the self-similarity of the $p$-Laplacian Equation can be found in \cite{Iag-San-Vaz}.
We have decided to dedicate an entire appendix to this topic since solutions of problem
\eqref{eq:CAUCHYPROBLEMPLAPLACIANINCREASINGINITIALDATA} play a main role in the proof of Theorem
\ref{CONVERGENCETOONEFASTDIFFUSION}. Moreover, we think it facilitates the reading and gives us the occasion to
present the related bibliography. Before proceeding with our analysis, we need to recall some important
 properties about problem \eqref{eq:CAUCHYPROBLEMPLAPLACIANINCREASINGINITIALDATA}.
\paragraph{Case $\boldsymbol{p = 2}$.} The existence and uniqueness of solutions for the Heat Equation for
continuous non-integrable initial has been largely studied, see Tychonov \cite{Tyc:art} and the references
therein. In particular, he  proved that if the initial datum satisfies
\begin{equation}\label{eq:ASSUMPTIONINDATUMP=2}
|u_0(x)| \leq b\exp\big(a|x|^2\big), \quad \text{for } |x| \sim \infty,
\end{equation}
for some positive $a$ and $b$, then problem \eqref{eq:CAUCHYPROBLEMPLAPLACIANINCREASINGINITIALDATA},
\eqref{eq:ASSUMPTIONINDATUMP=2} admits a unique (classical) solution defined in $\RR^N\times(0,1/(4a))$.
 More work on this issue can be found in \cite{Wid:art}.
\paragraph{Case $\boldsymbol{p > 2}$.} This range was studied in \cite{DBen-Her2:art}, by
DiBenedetto and Herrero. The authors showed that, under the assumptions
\begin{equation}\label{eq:ASSUMPTIONINDATUMPGEQ2DIBHER}
u_0 \in L_{loc}^1(\RR^N) \quad \text{ and } \quad u_0(x) \leq C|x|^{\lambda}, \text{ as } |x| \to \infty
\end{equation}
for some $C > 0$ and $\lambda < p/(p-2)$, there exists a unique weak solution of problem \eqref{eq:CAUCHYPROBLEMPLAPLACIANINCREASINGINITIALDATA}, \eqref{eq:ASSUMPTIONINDATUMPGEQ2DIBHER} defined in $\RR^N\times(0,\infty)$ (see Theorem 1, Theorem 2, and Theorem 4 of \cite{DBen-Her2:art}). Furthermore, they proved that $\partial_t u \in L_{loc}^2(\RR^N\times(0,\infty))$ (i.e. $u$ is a ``local strong solution'') and the function $(x,t) \to \nabla u(x,t)$ is locally H\"{o}lder continuous in $\RR^N\times(0,\infty)$ (see also \cite{DB-Fried:art}).
\paragraph{Case $\boldsymbol{1 < p < 2}$.} The same authors (see \cite{DBen-Her1:art}) considered problem \eqref{eq:CAUCHYPROBLEMPLAPLACIANINCREASINGINITIALDATA} with $1 < p < 2$ and nonnegative initial data
\begin{equation}\label{eq:ASSUMPTIONINDATUMP1LEQ2DIBHER}
u_0 \in L^1_{loc}(\RR^N) \quad \text{ and } \quad u_0(x) \geq 0 \; \text{ in } \RR^N
\end{equation}
without any assumption on the decay at infinity of $u_0(\cdot)$. First of all, they show existence and the uniqueness of weak solutions of problem \eqref{eq:CAUCHYPROBLEMPLAPLACIANINCREASINGINITIALDATA}, \eqref{eq:ASSUMPTIONINDATUMP1LEQ2DIBHER} by using the Benilan-Crandall regularizing effect, see \cite{Ben-Crand:art}. Then they posed their attention on the regularity of these solutions when the initial datum is a non-negative $\sigma$-finite Borel measure in $\RR^N$, in the range $2N/(N+1):=p_c < p < 2$. In particular, they showed the existence and the uniqueness of a locally H\"{o}lder continuous weak solution in $\RR^N\times(0,\infty)$, with $\partial_tu \in L_{loc}^2(\RR^N\times(0,\infty))$ (i.e. they are ``local strong solutions''), with $(x,t) \to \nabla u(x,t)$ locally H\"{o}lder continuous in $\RR^N\times(0,\infty)$.

\noindent The sub-critical range $1 < p \leq p_c:=2N/(N+1)$ was studied later by Bonforte, Iagar and V\'azquez in \cite{Bon-Vaz:art}. They proved new local smoothing effects when the initial datum is taken in $L^r_{loc}(\RR^N)$ and $p$ sub-critical, and special energy inequalities which are employed to show that bounded local weak solutions are indeed ``local strong solutions'', more precisely $\partial_tu \in L^2_{loc}(\RR^N)$. Then, thanks to the mentioned smoothing effect and known regularity theory (\cite{DB:book} and \cite{DBen-Urb-Ves:art}) they found that the local strong solutions are locally H\"{o}lder continuous.
\paragraph{Barenblatt solutions for problem (\ref{eq:CAUCHYPROBLEMPLAPLACIANINCREASINGINITIALDATA}), (\ref{eq:ASSUMPTIONINDATUMPOWERTYPE}).} From now on we take  $U_0(x) = |x|^{\lambda}$, $\lambda > 0$. We do not make any other assumptions on $\lambda > 0$ if $1 < p \leq 2$, whilst when $p > 2$ we assume $0 < \lambda < p/(p-2)$, according to the theory developed in \cite{DBen-Her2:art}, and presented before. As mentioned before, the assumptions on the parameter $\lambda$ guarantees the existence, the uniqueness and the H\"{o}lder regularity of the solution of problem \eqref{eq:CAUCHYPROBLEMPLAPLACIANINCREASINGINITIALDATA}, \eqref{eq:ASSUMPTIONINDATUMPOWERTYPE}, for all $p > 1$.

\noindent We look for solutions in \emph{self-similar} form
\[
U(x,t) = t^{-\alpha_{\lambda}} F(|x|t^{-\beta_{\lambda}}),
\]
where $\alpha_{\lambda}$ and $\beta_{\lambda}$ are real numbers and $F(\cdot)$ is called profile of the solution. Let $\xi = |x|t^{-\beta_{\lambda}}$ and write $F' = dF/d\xi$. It is not difficult to compute
\[
\partial_t U  = -t^{-\alpha_{\lambda}-1}(\alpha_{\lambda}F(\xi) + \beta_{\lambda}\xi F'(\xi)), \qquad
\Delta_p U = t^{-(\alpha_{\lambda} + \beta_{\lambda})(p-1) - \beta_{\lambda}}\xi^{1-N}\big(\xi^{N-1}|F'(\xi)|^{p-2}F'(\xi)\big)'
\]
and, by taking
\begin{equation}\label{eq:FIRSTCONDONALPHABETAINCREASININDATUM}
2\alpha_{\lambda} + 1 = (\alpha_{\lambda} + \beta_{\lambda})p,
\end{equation}
we have $\alpha_{\lambda} + 1 = (\alpha_{\lambda} + \beta_{\lambda})(p-1) + \beta_{\lambda}$, and so we obtain the equation of the profile
\[
\xi^{1-N}\big(\xi^{N-1}|F'(\xi)|^{p-2}F'(\xi)\big)' + \beta_{\lambda}\xi F'(\xi) + \alpha_{\lambda}F(\xi) = 0.
\]
Furthermore, since \eqref{eq:FIRSTCONDONALPHABETAINCREASININDATUM} guarantees that the equation in \eqref{eq:CAUCHYPROBLEMPLAPLACIANINCREASINGINITIALDATA} is invariant under the transformation $U_k(x,t) = k^{\alpha_{\lambda}}U(k^{\beta_{\lambda}}x,kt)$, $k > 0$, we use the uniqueness of the solution of problem \eqref{eq:CAUCHYPROBLEMPLAPLACIANINCREASINGINITIALDATA}, \eqref{eq:ASSUMPTIONINDATUMPOWERTYPE} to deduce
\[
k^{\alpha_{\lambda} + \lambda\beta_{\lambda}}|x|^{\lambda} = U_k(x,0) = U(x,0) = |x|^{\lambda}, \quad \text{for all } k > 0.
\]
Hence, we get $\alpha_{\lambda} + \lambda\beta_{\lambda} = 0$ and, combining it with \eqref{eq:CAUCHYPROBLEMPLAPLACIANINCREASINGINITIALDATA}, we obtain the precise expressions for the self-similar exponents
\[
\alpha_{\lambda} = -\frac{\lambda}{(1-\lambda)p + 2\lambda}, \qquad \beta_{\lambda} = \frac{1}{(1-\lambda)p + 2\lambda}.
\]
We point out that, thanks to the assumption $0 < \lambda < p/(p-2)$ when $p > 2$, we have $(1-\lambda)p + 2\lambda > 0$ for all $p>1$, and so $\alpha_{\lambda} < 0$ while $\beta_{\lambda} > 0$.
\paragraph{Properties of the Barenblatt solutions.} We are going to prove that the profile $F(\cdot)$ of the Barenblatt solutions is positive and monotone non-decreasing by applying the Aleksandrov's Symmetry Principle. Later, we show some asymptotic properties of the profile $F(\cdot)$.

\noindent Let $U_0(x) = |x|^{\lambda}$, with $0 < \lambda < p/(p-2)$ and, for all $j \in \NN$, consider the approximating sequence of initial data
\[
U_{0j}(x) :=
\begin{cases}
\begin{aligned}
|x|^{\lambda} \quad & \text{if } \; |x| \leq j \\
j^{\lambda}   \quad\;\, & \text{if } \; |x|\geq j.
\end{aligned}
\end{cases}
\]
Note that $U_{0j}(\cdot)$ are both radial non-decreasing and bounded in $\RR^N$. Now, consider the sequence of initial data
\[
v_{0j}(x) := j^{\lambda} - U_{0j}(x) \in C_c(\RR^N) \; \text{ and radial non-increasing},
\]
and the sequence of solutions $v_j(x,t)$ of problem \eqref{eq:CAUCHYPROBLEMPLAPLACIANINCREASINGINITIALDATA} with initial data $v_{0j}(\cdot)$, for all $j \in \NN$. Hence, by applying the Aleksandrov's Symmetry Principle, we deduce that for all times $t>0$, the solutions $v_j(\cdot,t)$ are radially non-increasing in space too. Finally, we define the sequence $U_j(x,t) = j^{\gamma} - v_j(x,t)$ which are radially non-decreasing in space and solve problem \eqref{eq:CAUCHYPROBLEMPLAPLACIANINCREASINGINITIALDATA} with initial data $U_{0j}$, for all $j \in \NN$. Hence, passing to the limit as $j \to \infty$, we have $U_j(x,t) \to U(x,t)$ and the limit $U(x,t)$, solution of problem \eqref{eq:CAUCHYPROBLEMPLAPLACIANINCREASINGINITIALDATA} with initial datum $U_0(\cdot)$, inherits the same radial properties of the sequence $U_j(x,t)$.

\noindent Now, we show the existence of two constants $0 < H_2 < H_1$ such that the following asymptotic bounds hold
\begin{equation}\label{eq:ESTIMATEONULARGEXI1}
H_2 |x|^{\lambda} \leq U(x,t) \leq H_1 |x|^{\lambda}, \quad \text{for } |x|t^{-\beta_{\lambda}} \sim \infty.
\end{equation}
Estimates \eqref{eq:ESTIMATEONULARGEXI1} follow directly from that fact that $U(x,t) \to |x|^{\lambda}$ as $t \to 0$. Indeed, for all fixed $0 \not= x \in \RR^N$, we have that
\[
\big|U(x,t) - |x|^{\lambda}\big| = t^{-\alpha_{\lambda}}\big|F(\xi) - \xi^{\lambda}\big| = |x|^{\lambda} \bigg| \frac{F(\xi)}{\xi^{\lambda}} - 1 \bigg|, \quad \text{ where } \; \xi = |x|t^{-\beta_{\lambda}}.
\]
Since, the left expression converges to 0 as $t \to 0$, we deduce that $F(\xi)/\xi^{\lambda} \to 1$, as $\xi \to \infty$ and, from this limit, we get \eqref{eq:ESTIMATEONULARGEXI1}. $\Box$
\paragraph{Aleksandrov's Symmetry Principle.} The Aleksandrov-Serrin symmetry method was firstly introduced in \cite{Alek:art} and \cite{Ser:art} to show monotonicity of solutions of both (eventually nonlinear) elliptic and parabolic equations. Here, following \cite{V2:book}, we give a short proof for the case of the ``pure diffusive'' $p$-Laplacian equation in \eqref{eq:CAUCHYPROBLEMPLAPLACIANINCREASINGINITIALDATA}, for all $p>1$.

\noindent Before proceeding with the statement, we fix some notations. Let $H$ be an hyperplane in $\RR^N$, $\Omega_1$ and $\Omega_2$ the two half-spaces ``generated'' by $H$, and $\Pi:\Omega_1 \to \Omega_2$ the reflection with respect to the hyperplane $H$.
\begin{thm}
Let $u\geq0$ be a solution of the initial-value problem \eqref{eq:CAUCHYPROBLEMPLAPLACIANINCREASINGINITIALDATA} with initial datum $u_0 \in L^1(\RR^N)$. Suppose that
\[
u_0(x) \geq u_0(\Pi(x)) \quad \text{for all } x \text{ in } \Omega_1.
\]
Then, for all times $t > 0$ it holds
\[
u(x,t) \geq u(\Pi(x),t) \quad \text{for all } x \text{ in } \Omega_1.
\]
In particular, radial initial data generate radial solutions.
\end{thm}
\emph{Proof}. First of all, thanks to the rotation invariance of the equation in \eqref{eq:CAUCHYPROBLEMPLAPLACIANINCREASINGINITIALDATA}, we can assume $H = \{x \in \RR^N : x_1 = 0\}$ and $\Pi(x_1,x_2,\ldots,x_N) = (-x_1,x_2,\ldots,x_N)$. Moreover, it follows that $\widehat{u}(x,t) = u(\Pi(x),t)$ solves problem \eqref{eq:CAUCHYPROBLEMPLAPLACIANINCREASINGINITIALDATA} in $\RR^N\times(0,\infty)$ with initial datum $\widehat{u}_0(x) = u(\Pi(x),0)$.

\noindent Now, we have $u_0(x) \geq u_0(\Pi(x))$ in $\Omega_1$ and $u(x,t) = \widehat{u}(x,t)$ in $H\times(0,\infty) = \partial\Omega_1\times(0,\infty)$. Hence, since the solution is continuous, we get the thesis by applying the Maximum Principle. Note that, to be precise, we should consider solutions of the Cauchy-Dirichlet problem posed in the ball $B_R(0)$ with zero boundary data. These solutions approximate $u = u(x,t)$ and $\widehat{u} = \widehat{u}(x,t)$. Consequently, we can apply the Maximum Principle to these approximate solutions and, finally, pass to the limit as $R \to \infty$. See Chapter 9 of \cite{V2:book} for more details.

\noindent If $u_0(\cdot)$ is radial, we can apply the statement for all hyperplane $H$ passing through the origin of $\RR^N$ and deducing that for all times $t > 0$, the solution $u(\cdot,t)$ is radial respect with the spacial variable too. $\Box$
%
%
%
%
%
%
%
%
%
%
%
\section{Comments and open problems}\label{SECTIONFINALREMARKSFAST}
We end the paper by discussing some open problems. Moreover, we present some final comments and remarks to supplement our work.

\noindent As we have mentioned in the introduction, nonlinear evolution processes give birth to a wide variety of phenomena. Indeed we have seen that solutions of problem \eqref{eq:REACTIONDIFFUSIONEQUATIONPLAPLACIAN} exhibit a travelling wave behaviour for large times when $\gamma \geq 0$, i.e. $\widehat{\gamma} \leq 0$, while infinite speed of propagation when $0 < \widehat{\gamma} < p/N$. It is natural to ask ourselves what happens in the range of parameters $\widehat{\gamma} \geq p/N$ that we call ``very fast diffusion assumption''.

\noindent However, respect to the Porous Medium and the $p$-Laplacian case, we have to face the problem of lack of literature and previous works related to the doubly nonlinear operator (in this range of parameters). For this reason, in the next paragraphs we will briefly discuss what is known for the Porous Medium and the $p$-Laplacian case, trying to guess what could happen in the presence of the doubly nonlinear operator. We stress that our approach is quite \emph{formal}, but can be interesting since it gives a more complete vision of the fast diffusion range, and allows us to explain what are (or could be) the main differences respect to the range $0 < \widehat{\gamma} < p/N$.
\paragraph{The critical case $\boldsymbol{\widehat{\gamma} = p/N}$.} This critical case was firstly studied by King \cite{King:art} and later in \cite{Galak-Pel-Vaz:art} for the Porous Medium setting, i.e. $p=2$ and $m = m_c := (N-2)_+/N$, with $N \geq 3$. When $N = 1,2$ it follows $m = 0$, choice of parameter which goes out of our range and we avoid it. King studied the asymptotic behaviour of radial solutions of the pure diffusive equation
\begin{equation}\label{eq:PUREDIFFUSIVEDNLFAST}
\partial_t \underline{u} = r^{1-N}\partial_r\big(r^{N-1}|\partial_r\underline{u}^m|^{p-2} \partial\underline{u}^m\big) \quad \text{ in } \RR^+\times [0,\infty)
\end{equation}
with $p = 2$, $0 < m \leq m_c$, and $N \geq 3$. Actually, he considered a slightly different equation absorbing a factor $m^{p-1}$ in the time varible and he studied the cases $N=2$ and $m=0$, too. Note that the choice $\widehat{\gamma} = p/N$ corresponds to $m = m_c$ when $p=2$.

\noindent In \cite{King:art}, the author described the asymptotic behaviour of radial solutions of equation \eqref{eq:PUREDIFFUSIVEDNLFAST}, given by the formula
\[
\underline{u}(r,t) \sim \bigg(\frac{(N-2)t}{r^2\ln r}\bigg)^{\frac{N}{2}}, \quad \text{as } t \sim \infty \quad \text{and} \quad t^{-N/(N-2)}\ln r \geq \eta_0,
\]
where $\eta_0$ is a constant depending on $N$ and on the initial datum (see formula (2.34) of \cite{King:art}). In particular, it follows that the solutions of \eqref{eq:PUREDIFFUSIVEDNLFAST} have spacial power like decay $r^{-N}$ ``corrected'' by a logarithmic term for $r \sim \infty$. We are interested in seeing that an analogue decay holds when $p > 1$ and $\widehat{\gamma} = p/N$ in the doubly nonlinear setting. We proceed as in \cite{King:art}, see Section 2.

\emph{Asymptotic behaviour for large $r$.} Let's take for a moment $0 < \widehat{\gamma} < p/N$. Seeking solutions $\underline{u} = \underline{u}(r,t)$ of equation \eqref{eq:PUREDIFFUSIVEDNLFAST} in separate form as in \emph{Step1} of Theorem \ref{THEOREMBOUNDSFORLEVELSETSFAST}, it is simple to see that if the initial datum satisfies \eqref{eq:ASSUMPTIONSONTHEINITIALDATUMFAST}, then for all $t > 0$ we have
\begin{equation}\label{eq:ASYMPTOTICEXPANSIONRLARGEFAST}
\underline{u}(r,t) \sim a\,t^{\frac{1}{\widehat{\gamma}}} \, r^{-\frac{p}{\widehat{\gamma}}}, \quad \text{for } r \sim \infty,
\end{equation}
for some suitable constant $a > 0$. Note that it corresponds to fix $t > 0$ and take the limit as $r \to \infty$ in the formula of the Barenblatt solutions, see Subsection \ref{PRELIMINARIESINTRO}.

\noindent Now, motivated by the previous analysis, we fix $N > p$ (in order to remain in the ranges $m>0$ and $p > 1$), $\widehat{\gamma} = p/N$, and we look for solutions of equation \eqref{eq:PUREDIFFUSIVEDNLFAST} in the form
\[
\underline{u}(r,t) \sim a\,t^{N/p}\,r^{-N}F(r), \quad \text{for } r \sim \infty,
\]
for some correction function $0 \leq F(r) \to 0$ as $r \to \infty$ and some constant $a > 0$. In what follows we ask $rF'(r) = o(F(r))$ as $r \to \infty$, too. It is simple to compute
\[
\begin{aligned}
\partial_t \underline{u} &\sim (aN/p)t^{\frac{N}{p}-1}r^{-N}F(r) \\
\partial_r\underline{u}^m &= ma^mt^{\frac{Nm}{p}}r^{-Nm-1}F(r)^{m-1}(-NF(r) + rF'(r)) \sim  -Nma^m\,t^{\frac{Nm}{p}}r^{-Nm-1}F(r)^m \\
|\partial_r\underline{u}^m|^{p-2} \partial\underline{u}^m &\sim -(Nm)^{p-1}a^{\frac{N-p}{N}}\,t^{\frac{N-p}{p}}r^{1-N}F(r)^{\frac{N-p}{N}}
\end{aligned}
\]
as $r \sim \infty$, where we have used the fact that $m(p-1) = 1 - \widehat{\gamma} = 1- p/N$. Hence, it is simple to see that $\underline{u} = \underline{u}(r,t)$ solves \eqref{eq:PUREDIFFUSIVEDNLFAST} if and only if
\[
(mN)^{p-1}a^{\frac{N-p}{N}}\,r\,\Big(F(r)^{\frac{N-p}{N}}\Big)' + a(N/p)F(r) = 0.
\]
Now, it is clear that a possible choice is $F(r) = (\ln r)^{-b}$, for some $b > 0$, and a straightforward computation shows that the previous equation is satisfied by taking
\[
a^{\frac{p}{N}} = m^{p-1}(N-p)N^{p-2} \quad \text{ and } \quad b = \frac{N}{p},
\]
so that for all $t > 0$, we obtain
\begin{equation}\label{eq:ASYMPTOTICEXPANSIONSOLPUREDIFFRLARGEFASTCRIT}
\underline{u}(r,t) \sim \bigg(\frac{a^{p/N}t}{r^p\ln r}\bigg)^{\frac{N}{p}}, \quad \text{for } r \sim \infty,
\end{equation}
which generalizes the case $p = 2$ and $m = m_c$.

\emph{Barenblatt solutions for $\widehat{\gamma} = p/N$.} As was observed in \cite{King:art} (see pag. 346), \eqref{eq:ASYMPTOTICEXPANSIONSOLPUREDIFFRLARGEFASTCRIT} does not respect the self-similarity reduction of equation \eqref{eq:PUREDIFFUSIVEDNLFAST}. Indeed, it admits ``pseudo-Barenblatt'' solutions which, following the notation of \eqref{eq:SECONDFORMULATIONBARSOLFAST}, can be written in the form $B_D(x,t) = R(t)^{-N}F_D(xR(t)^{-1})$ where
\[
F_D(\xi) = \Big[D + (1/N)|\xi|^{\frac{p}{p-1}} \Big]^{-\frac{(p-1)}{p}N},  \qquad \xi = xR(t)^{-1}, \qquad R(t) = e^t,
\]
and $D > 0$ is a free parameter (cfr. with \cite{V6:art} for the case $p=2$ and with formula \eqref{eq:SECONDFORMULATIONBARSOLFAST} for the range $0 < \widehat{\gamma} < p/N$). We point out that the profile $F_D(\cdot)$ satisfies the inequalities in \eqref{eq:ESTIMATEONPROFILEOFMASS1} with $\widehat{\gamma} = p/N$. However, these self-similar solutions (also called ``of Type III'', see \cite{V1:book}) are quite different from the ones in the range $0 < \widehat{\gamma} < p/N$. In particular, they are eternal, i.e. defined for all $t \in \RR$ and they do not converge to a Dirac Delta as $t \to 0$ (see also \cite{Br-Fr:art}). Finally, for all fixed $t \in \RR$, these self-similar solutions are \emph{not integrable} respect with to the spacial variable and show the spacial decay
\[
B_D(x,t) \sim N^{\frac{(p-1)N}{p}} |x|^{-N}, \quad \text{for } |x| \sim \infty.
\]
Taking into account these facts, when $\widehat{\gamma} = p/N$ it seems reasonable to study problem \eqref{eq:REACTIONDIFFUSIONEQUATIONPLAPLACIAN} with nontrivial initial datum satisfying
\[
0 \leq u_0(x) \leq 1 \quad \text{ and } \quad u_0(x) \leq C\big(|x|^p\ln|x|\big)^{-\frac{N}{p}} \quad \text{for } |x| \sim \infty,
\]
for some constant $C > 0$, and trying to extend the techniques used for the range $0 < \widehat{\gamma} < p/N$, to the this critical case.

\noindent First of all, we can define $\sigma_{\ast} := f'(0)/N$ by continuity. Thus, it is possible to repeat the proof of Theorem \ref{CONVERGENCETOZEROFASTDIFFUSION} by using ``pseudo-Barenblatt'' solutions instead of the usual ones. In this way, for all $\sigma > \sigma_{\ast}$, we show the convergence of the solutions to 0 in the ``outer sets'' $\{|x| \geq e^{\sigma t}\}$, as $t \to \infty$.

\noindent Moreover, thanks to the asymptotic expansion \eqref{eq:ASYMPTOTICEXPANSIONSOLPUREDIFFRLARGEFASTCRIT} it should be possible to prove a version of Lemma \ref{LEMMAPLACINGBARENBLATTUNDERSOLUTION} with
\begin{equation}\label{eq:CRITICALINITIALDATUMFAST}
\widetilde{u}_0(x) :=
\begin{cases}
\begin{aligned}
\widetilde{\varepsilon}  \qquad\qquad\qquad\quad\;\; &\text{if }|x| \leq \widetilde{\varrho}_0 \\
a_0\big(|x|^p\ln|x|\big)^{-\frac{N}{p}} \quad &\text{if }|x| > \widetilde{\varrho}_0,
\end{aligned}
\end{cases}
\end{equation}
for $a_0 := \widetilde{\varepsilon}\,\big(\widetilde{\varrho}_0^{\,p}\ln\widetilde{\varrho}_0\big)^{N/p}$ and some $0 < \widetilde{\varepsilon} < 1$ and $\widetilde{\varrho}_0 > 1$.

\noindent However, it is clear that the methods employed for showing Proposition
\ref{EXPANPANSIONOFMINIMALLEVELSETS} cannot be used in this case too. Indeed, in the range
$0 < \widehat{\gamma} < p/N$, this crucial proposition has been proved by constructing barriers from below with
Barenblatt solutions. This has been possible since the initial datum $\widetilde{u}_0 = \widetilde{u}_0(x)$ in
\eqref{eq:INITIALDATUMTESTFASTLEVELSETS1} shares the same spacial decay of these self-similar solutions. In the
critical case $\widehat{\gamma} = p/N$, this property would not be preserved as
\eqref{eq:ASYMPTOTICEXPANSIONSOLPUREDIFFRLARGEFASTCRIT} suggests. In particular, ``pseudo-Barenblatt'' solutionscannot be placed under an initial datum satisfying \eqref{eq:CRITICALINITIALDATUMFAST} and so the
validity of Proposition \ref{EXPANPANSIONOFMINIMALLEVELSETS} in this critical case remains an open problem.
\paragraph{The range $\boldsymbol{\widehat{\gamma} > p/N}$.} Before discussing the doubly nonlinear diffusion,
let us recall what is known in the Porous Medium setting in the corresponding range of parameters,
$0 < m < m_c:=(N-2)_+/N$, $p=2$, and $N \geq 3$. Consider the Porous Medium Equation
\[
\begin{cases}
\begin{aligned}
\partial_tv = \Delta v^m \qquad\qquad\;\;\, &\text{in } \RR^N\times(0,\infty) \\
v(x,0) = v_0(x) \quad\qquad &\text{in } \RR^N,
\end{aligned}
\end{cases}
\]
where $v_0 \in L^1(\RR^N)\cap L^{\infty}(\RR^N)$. It has been proved that the corresponding solution $v = v(x,t)$ extinguishes in finite time (see for instance \cite{Ben-Crand1:art,King:art,V1:book} and the references therein). In other words, there exists a critical ``extinction time'' $0 < t_c < \infty$ such that $v(\cdot,t) = 0$, for all $t \geq t_c$. Again, the cases $m = 0$, $N = 1$ and $m = 0$, $N = 2$ are critical and we refer to \cite{V1:book}, Chapters 5 to 8.

\emph{Barenblatt solutions for $\widehat{\gamma} > p/N$.} So, even though there is not literature on the subject (at least to our knowledge), it seems reasonable to conjecture that the doubly nonlinear diffusion shows a similar property in the range $\widehat{\gamma} > p/N$, with $N > p$. In particular, also in this case we have ``pseudo-Barenblatt'' solutions written in the form
\begin{equation}\label{eq:PSEUDOBARSOLVERYFAST}
B_D(x,t) = R(t)^{-N} \Big[ D + (\widehat{\gamma}/p)\big|x R(t)^{-1} \big|^{\frac{p}{p-1}} \Big]^{-\frac{p-1}{\widehat{\gamma}}},
\end{equation}
where $D \geq 0$ and, with a strong departure from \eqref{eq:SECONDFORMULATIONBARSOLFAST},
\[
R(t) = \big[(N/|\alpha|)(t_c-t) \big]^{-\frac{|\alpha|}{N}},
\]
where $t_c > 0$ is fixed and stands for the ``extinction time''(cfr. with \cite{V1:book} pag. 194 or  \cite{V6:art} for the case $p=2$, and with formula \eqref{eq:SECONDFORMULATIONBARSOLFAST} for the range $0 < \widehat{\gamma} < p/N$). The existence of this kind of self-similar solutions (also said in \cite{V1:book} ``of Type II'') strengthen the idea that a larger class of solutions have an extinction time, i.e. they vanish in finite time.

\emph{Application to the Fisher-KKP equation.} In Section \ref{CONVERGENCETOZEROFAST} we have seen that the linearized problem
\[
\begin{cases}
\begin{aligned}
\partial_t\overline{u} = \Delta_p\overline{u}^m + f'(0)\overline{u} \quad &\text{in } \RR^N\times(0,\infty) \\
\overline{u}(x,0) = u_0(x) \;\;\quad\qquad &\text{in } \RR^N,
\end{aligned}
\end{cases}
\]
gives a super-solution for the Fisher-KPP problem \eqref{eq:REACTIONDIFFUSIONEQUATIONPLAPLACIAN} with nontrivial initial datum $u_0 \in L^1(\RR^N)$, $0 \leq u_0 \leq 1$. Again, with the change of variable
\[
\tau(t) = \frac{1}{f'(0)\widehat{\gamma}}\Big[1 - e^{-f'(0)\widehat{\gamma} t} \Big], \quad \text{for } t \geq 0,
\]
we deduce that the function $\overline{v}(x,\tau) = e^{-f'(0)t}\overline{u}(x,t)$ solves the problem
\begin{equation}\label{eq:PUREDIFFUSIVEDNLEQUATFINITETIME}
\begin{cases}
\begin{aligned}
\partial_{\tau}\overline{v} = \Delta_p\overline{v}^m \,\qquad &\text{in } \RR^N\times(0,\tau_{\infty}) \\
\overline{v}(x,0) = u_0(x) \quad &\text{in } \RR^N.
\end{aligned}
\end{cases}
\end{equation}
Now, set $\tau_{\infty} := \frac{1}{f'(0)\widehat{\gamma}}$ and note that $0 \leq \tau(t) \leq \tau_{\infty}$. Now, let $\tau_c > 0$ be the ``extinction time'' of the solution of problem \eqref{eq:PUREDIFFUSIVEDNLEQUATFINITETIME}. Thus, we deduce $v(\cdot,\tau) = 0$, for all $\tau \geq \tau_c$ and, if $\tau_c < \tau_{\infty}$, it follows
\[
0 \leq u(\cdot,t) \leq\overline{u}(\cdot,t) = e^{f'(0)t}\overline{v}(\cdot,\tau) = 0, \quad \text{for all } \tau \geq \tau_c,
\]
which implies $u(\cdot,t) = 0$ for all $t \geq \tau_{\infty}\ln\big[ \tau_{\infty}/(\tau_{\infty} - \tau_c) \big]$, and so the solution $u = u(x,t)$ of the Fisher-KPP problem \eqref{eq:REACTIONDIFFUSIONEQUATIONPLAPLACIAN} with initial datum $u_0$ extinguishes in finite time, too. This conclusion holds under the assumption $\tau_c < \tau_{\infty}$, which should be guaranteed if the the initial datum is ``small enough'' (in terms of the mass), see \cite{V1:book} Chapter 5, for the Porous Medium setting. The analysis of the case in which the initial mass is infinite is an interesting open problem.
%
%
%
%
%
%
%
%
%
%
%

\vskip 1cm


\noindent {\textbf{\large \sc Acknowledgments.}} Both authors have been partially funded by Projects MTM2011-24696 and  MTM2014-52240-P (Spain). Work partially supported by the ERC Advanced Grant 2013 n.~339958 ``Complex Patterns for Strongly Interacting Dynamical Systems - COMPAT''. We thank F\'elix del Teso for providing us with the numerical simulations that we display.

%
%
%
\vskip 1cm

\

2000 \textit{Mathematics Subject Classification.}
35K57,  
35K65, 
35C07,  	
35K55, 

\medskip

\textit{Keywords and phrases.}  Fisher-KPP equation, Doubly nonlinear diffusion, Propagation of level sets.

\

\end{document}